\documentclass[12pt,a4paper,reqno]{amsart}
\usepackage{amsmath}
\usepackage{amsfonts}
\usepackage{mathrsfs}
\usepackage{amssymb}
\usepackage{amsthm}
\usepackage{enumerate}
\usepackage[all]{xy}
\usepackage{graphicx}
\usepackage{fullpage}
\usepackage{mathtools}
\usepackage{float}
\usepackage[bookmarksnumbered=true]{hyperref}

\DeclareMathOperator{\diam}{diam}

\DeclareMathOperator{\wind}{wind}

\theoremstyle{plain}
  \newtheorem{theorem}{Theorem}[section]
  \newtheorem{lemma}[theorem]{Lemma}
  
  \newtheorem{proposition}[theorem]{Proposition}

\theoremstyle{remark}
    \newtheorem{definition}[theorem]{Definition}
    \newtheorem{remark}[theorem]{Remark}

\title{The loop-ensemble scaling limits of near-critical percolation in the plane}
\author{Chang-Long Yao}\thanks{State Key Laboratory of Mathematical Sciences, Academy of Mathematics and Systems Science, Chinese Academy of Sciences, Beijing 100190, China. E-mail: deducemath@126.com}

\begin{document}
\maketitle

\begin{abstract}
We prove that the collection of all cluster boundaries of near-critical site percolation (with a fixed parameter) on the triangular lattice has a unique scaling limit in the loop-ensemble space, as the mesh size tends to 0.  This answers a problem proposed by Camia et al. \cite{CFN06}.  Using this, we also prove that near-critical percolation exploration path in a disk converges in distribution to a continuous path as the mesh size goes to 0.  The main idea of our proof is to construct the cut-off loop ensembles under the near-critical coupling, based on the enhanced networks introduced in \cite{GPS18b} and the quad-crossing scaling limits of near-critical percolation established in \cite{GPS18a}; the proof also relies on the study of the geometry of cluster chains.  Some properties of the limiting loop ensemble, including rotation invariance, scale covariance and finite chaining, are established.  Moreover, a few applications of our scaling-limit result and proof strategy are discussed.

\noindent\textbf{MSC2020 subject classifications:} 60K35; 82B43
\end{abstract}

\maketitle

\section{Introduction}
\subsection{Background and motivation}\label{bm}
Percolation is a central model in probability theory and statistical physics for studying phase transitions; see the classical books \cite{BR06,Gri99} and a recent survey \cite{DC18} for general background.  Critical planar percolation is particularly interesting since it has a conformally invariant scaling limit.  The breakthrough works of Schramm \cite{Sch00} and Smirnov \cite{Smi01} allows us to describe the scaling limit for the interfaces of critical site percolation on the triangular lattice $\textsf{T}$ by Schramm-Loewner Evolution with parameter 6 (SLE$_6$); see the surveys \cite{BD13,Sun11,Wer09}.  Based on this, Camia and Newman \cite{CN06} proved that as the mesh size $\eta$ goes to zero, the collection of all cluster boundaries of critical percolation on $\eta\textsf{T}$ converges to a collection of SLE$_6$-type loops, known as the Conformal Loop Ensemble with parameter 6 (CLE$_6$).

General CLE$_{\kappa}$ with $\kappa\in(8/3,8)$ was introduced by Sheffield in \cite{She09}, which is a conformally invariant countable collection of SLE$_{\kappa}$-like random loops (it is a single space-filling loop at the extreme $\kappa=8$ and is an empty set at $\kappa=8/3$).  Note that SLEs describe the scaling limit of a single interface arising in discrete statistical physics models at criticality, while CLEs describe the ``full scaling limit" of the collection of all interfaces.  Besides the CLE$_6$ result for percolation, similar convergence results have been proved for the Ising model \cite{BH19} ($\kappa=3$), the single random current \cite{CX} ($\kappa=3$), the double random current \cite{DCLQ25,DCLQ} ($\kappa=4$) and the FK-Ising model \cite{KS19} ($\kappa=16/3$).

The behavior of near-critical percolation (with the parameter $p$ being close to the critical value $p_c$) is closely related to that of critical percolation.
The scaling relations derived by Kesten \cite{Kes87}, together with SLE techniques and the scaling limit result for critical percolation, leads to the derivation of important properties of both critical and near-critical percolation; see the surveys \cite{Nol08,Wer09}.  For example, it was proved in \cite{SW01} that $\theta(p)=(p-p_c)^{5/36+o(1)}$ as $p\downarrow p_c$, where $\theta(p)$ is the probability that the origin belongs to the infinite open cluster.  Moreover, with the help of scaling relations, some basic geometric properties of subsequential near-critical scaling limits for the collection of boundary loops were obtained in \cite{CJM09}.  In \cite{CFN06,CFN06b}, Camia, Fontes and Newman proposed a conceptual framework to construct a one-parameter family of loop-ensemble scaling limits for near-critical percolation with $p=p_c+\lambda\eta^{3/4}$, where $\eta$ is the mesh size and $\lambda\in\mathbb{R}$.  Their framework combines CLE$_6$ with a ``Poissonian marking" of the double (i.e., touching) points of CLE$_6$.  Then, Garban, Pete and Schramm \cite{GPS13,GPS18a} realized a modified version of this framework, using the quad-crossing topology introduced in \cite{SS11}.  More precisely, they proved that for any fixed $\lambda\in\mathbb{R}$, the near-critical site percolation on $\eta\textsf{T}$ with
\begin{equation}\label{e5}
p=p_c+\lambda r(\eta)
\end{equation}
has a unique quad-crossing scaling limit as $\eta\rightarrow 0$, where $r(\eta):=\eta^2/\alpha_4(\eta,1)$ and $\alpha_4(\eta,1)$ stands for the probability of the
alternating 4-arm event from a given site to distance 1 on $\eta\textsf{T}$.

Similarly to SLE$_6$ describing the scaling limit of a single interface of critical percolation, it is natural to expect that there are near-critical versions of SLE$_6$ which describe the near-critical scaling limits. It is well known that the tightness criterion of Aizenman and Burchard \cite{AB99} together with the near-critical RSW implies the existence of subsequential scaling limits for near-critical interfaces.  Furthermore, Nolin and Werner \cite{NW09} showed that for $\lambda\neq 0$, any subsequential scaling limits for the interfaces of near-critical percolation at $p=p_c+\lambda r(\eta)$ are singular with respect to the classical SLE$_6$ curves; in \cite{GPS18a}, the authors made a conjecture on the precise form of the driving function for such limiting curves.  These curves are the so-called ``massive" SLE$_6$ in \cite{MS10}, where the term massive refers to a system with exponential decay of correlations that is described by a massive field theory in the physics literature, while critical planar systems are described by conformal massless field theories.  Moreover, Makarov and Smirnov \cite{MS10} made a general discussion of massive SLE$_{\kappa}$'s that are conjectured to be the scaling limits of near-critical discrete models.  Note that following the Makarov-Smirnov program, convergence of near-critical interfaces has been established in the following cases: massive loop-erased random walk \cite{CW21} (massive SLE$_2$), massive harmonic explorer \cite{Pap25} (massive SLE$_4$), the interface of critical planar Ising model perturbed by an external magnetic field \cite{Pap} (massive SLE$_3$).  We remark that for $\kappa\leq 4$ the massive and usual SLE$_\kappa$ are expected to be mutually absolutely continuous.  Is it true that for $\kappa>4$ the massive and usual SLE$_\kappa$ are mutually singular? This question was posed in \cite[Question 4.6]{MS10}.  Recently, Garban and Kupiainen \cite[Section 2.5]{GK} conjectured that massive SLE$_{16/3}$ is absolutely continuous with respect to SLE$_{16/3}$.  Moreover, they proved that for near-critical FK-Ising model, any subsequential scaling limit of the collection of cluster-boundary loops (corresponding to massive CLE$_{16/3}$) is singular with respect to CLE$_{16/3}$, and conjectured that massive CLE$_3$ is also singular with respect to CLE$_3$ (see \cite[Section 2.4]{GK}).

The purpose of this paper is three-fold:
\begin{itemize}
  \item In the spirit of the CLE$_6$ scaling limit result at $p_c$, we obtain an analogous result for $\lambda$-near-critical percolation in the loop space, which extends the quad-crossing scaling limit result of \cite{GPS18a} and answers a problem proposed in \cite{CFN06}.

  \item In the spirit of the SLE$_6$ scaling limit result at $p_c$, we prove that the $\lambda$-near-critical percolation exploration paths converge in distribution to a continuous path as the mesh size tends to 0, which establishes the uniqueness of subsequential scaling limits for these paths.  It was mentioned in \cite[Section 1.1]{Pap} that ``proving uniqueness of the limiting law remains a challenge", although some progress was made in \cite{NW09,GPS18a}.

   \item For some discrete models, extracting geometric information or expressing some quantities is sometimes more straightforward from the scaling limit living in the loop space than in the quad-crossing space.  We give a few applications and basic properties of our loop-ensemble scaling limits.
\end{itemize}

\subsection{The model and near-critical coupling}\label{s3}
\begin{figure}
\centering
\includegraphics[height=0.37\textwidth]{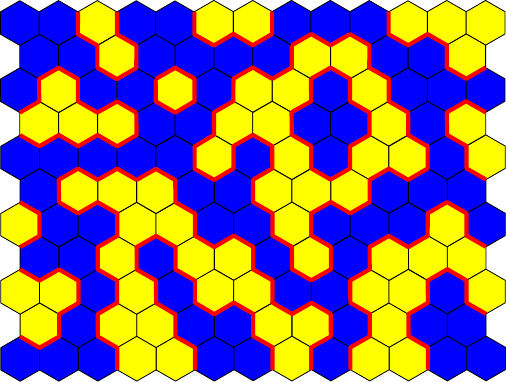}
\caption{Finite portion of a site percolation configuration on the triangular lattice $\textsf{T}$.  Each hexagon of
the hexagonal lattice $\textsf{H}$ represents a site of $\textsf{T}$, and is colored blue (open) or yellow (closed).  The red lines
represent the cluster boundaries.}\label{hexlattice}
\end{figure}

Let $\textsf{T}$ denote the triangular lattice embedded in $\mathbb{C}$ (we identify the plane
$\mathbb{R}^2$ with the set $\mathbb{C}$ of complex numbers in the usual way), with vertex (i.e., site) set
\begin{equation*}
V(\textsf{T}):=\{x+ye^{\pi i/3}\in\mathbb{C}:x,y\in\mathbb{Z}\},
\end{equation*}
and edge (i.e., bond) set $E(\textsf{T})$ obtained by connecting all pairs
$u,v\in V(\textsf{T})$ for which ${|u-v|=1}$.  Site percolation on $\textsf{T}$ is defined as follows:  For each $p\in [0,1]$, we consider the i.i.d. family $\{\omega(v):v\in V(\textsf{T})\}$ of Bernoulli random variables with parameter $p$, that is, $\omega(v)=1$ with probability $p$ and $\omega(v)=0$ with probability $1-p$.  This gives rise to a product probability measure $\mathbf{P}_p$ on the set
of configurations $\{0,1\}^{V(\textsf{T})}$.  A vertex $v$ is declared \textbf{open} in the configuration $\omega$ if $\omega(v)=1$ and \textbf{closed} if $\omega(v)=0$.  We usually represent site percolation on $\textsf{T}$ as a random two-coloring of the faces of the dual regular hexagonal lattice $\textsf{H}$, with
blue hexagons standing for open sites and yellow hexagons standing for closed sites; see Figure \ref{hexlattice}.  We will often view a site $v$ as the hexagon in $\textsf{H}$ centered at $v$.  The connected open (resp. closed) components are called open (resp. closed) \textbf{clusters}.  It is well known that site percolation on $\textsf{T}$ exhibits a phase transition at $p_c=1/2$: When $p\leq p_c$ there are almost surely no infinite open clusters, whereas when $p>p_c$ there is almost surely a unique infinite open cluster.

Recall that
\begin{equation}\label{e3}
r(\eta):=\eta^2/\alpha_4(\eta,1),
\end{equation}
where $\alpha_4(\eta,1)$ denotes the probability of the
alternating 4-arm event from a given site to distance 1 for critical site percolation on
$\eta\textsf{T}$.  See also \cite{GPS18a,GPS18b} where the same notation was used.  Recently, it has been proved in \cite{DGLZ24} that $\alpha_4(\eta,1)\asymp\eta^{5/4}$.

To study near-critical scaling limits, Garban, Pete and Schramm \cite{GPS18a} introduced the following very useful monotone coupling of near-critical percolation configurations $(\omega_{\eta}^{\lambda})_{\lambda\in\mathbb{R}}$ on $\eta\textsf{T}$, with the near-critical parametrization given in (\ref{e5}) being slightly changed.
\begin{definition}\label{def1}
The near-critical coupling $(\omega_{\eta}^{\lambda})_{\lambda\in\mathbb{R}}$ will denote the following process:
\begin{enumerate}[(i)]
\item Sample $\omega_{\eta}^{\lambda=0}$ according to the law of critical site percolation on $\eta\textsf{T}$.
\item As $\lambda$ increases, closed sites switch to open at exponential rate $r(\eta)$.
\item As $\lambda$ decreases, open sites switch to closed at exponential rate $r(\eta)$.
\end{enumerate}
\end{definition}
Then, for any $\lambda\in\mathbb{R}$, the near-critical percolation $\omega_{\eta}^{\lambda}$ corresponds exactly to a
percolation configuration on $\eta\textsf{T}$ with parameter
\begin{equation}\label{e42}
p=\left\{
\begin{aligned}
1-(1-p_c)e^{-\lambda r(\eta)}&=1-e^{-\lambda r(\eta)}/2&\mbox{ if }\lambda\geq 0,\\
p_ce^{-|\lambda| r(\eta)}&=e^{-|\lambda| r(\eta)}/2&\mbox{ if }\lambda<0.
\end{aligned}
\right.
\end{equation}
Let $\mathbf{P}_{\eta}^{\lambda}$ denote the law of $\omega_{\eta}^{\lambda}$.  Note that the process $(\omega_{\eta}^{\lambda})_{\lambda\in\mathbb{R}}$ is a time-homogeneous Markov process since $p_c=1/2$.  Moreover, the percolation parameter $p$ for $\omega_{\eta}^{\lambda}$ satisfies $p\approx p_c+\lambda r(\eta)/2$ when $\eta$ is small.

\subsection{Main results}
In this subsection, we provide a quick presentation of our main results, with some precise definitions postponed.

Given a site percolation on $\textsf{T}$ with parameter $p\in(0,1)$, we orient each edge of $\textsf{H}$ which belongs to the topological boundary of a cluster in such a way that the hexagon to the right of the oriented edge is blue (open).  It is well known that with probability one, the set of all the boundaries of the clusters forms a collection of oriented simple loops in the plane (see Figure \ref{hexlattice}), with the orientation of each loop according to the oriented edges in the loop.  For $\lambda\in\mathbb{R}$, let $\Gamma_{\eta}^{\lambda}$ denote the collection of all cluster boundary loops of the site percolation $\omega_{\eta}^{\lambda}$ on $\eta\textsf{T}$.  Let $(\mathcal {L},D_\mathcal {L})$ be the metric space of loop ensembles on $\hat{\mathbb{C}}$ defined in Section \ref{loop}.  When $\lambda=0$ (i.e., the critical case), it was proved in \cite{CN06} that $\Gamma_{\eta}^0$ converges in distribution in $(\mathcal {L},D_\mathcal {L})$ to the full-plane CLE$_6$ as $\eta\rightarrow 0$.  The following theorem extends this result to the near-critical case.
\begin{theorem}\label{t1}
Fix $\lambda\in\mathbb{R}\setminus\{0\}$.  As the mesh size $\eta$ goes to zero, $\Gamma_{\eta}^{\lambda}$ converges in law in $(\mathcal {L},D_\mathcal {L})$ to a random loop ensemble, which we denote by $\Gamma_{\infty}^{\lambda}$.
\end{theorem}

The next theorem says that for $\lambda\neq 0$, $\Gamma_{\infty}^{\lambda}$ is invariant under translations and rotations, is covariant under scalings, and has a form of orientation-reversing symmetry.  Given a transformation $f$ of $\hat{\mathbb{C}}$, we set $f(\Gamma_{\infty}^{\lambda}):=\{f(\ell):\ell\in\Gamma_{\infty}^{\lambda}\}$.  We shall view a transformation of $\mathbb{C}$ as a transformation of $\hat{\mathbb{C}}$ that maps infinity to itself.
\begin{theorem}\label{t7}
Let $\lambda\in\mathbb{R}\setminus\{0\}$.  Then $\Gamma_{\infty}^{\lambda}$ has the following properties:
\begin{enumerate}[(i)]
   \item (\textbf{Translation and rotation invariance})  For any translational/rotational transformation $f$ of the plane, $\Gamma_{\infty}^{\lambda}$ and $f(\Gamma_{\infty}^{\lambda})$ have the same distribution. \label{rotation}
   \item (\textbf{Scale covariance}) For any $\alpha>0$, let $f_{\alpha}$ denote the scaling transformation of the plane with parameter $\alpha$, that is, $f_{\alpha}(z):=\alpha z$ for $z\in \mathbb{C}$.  Then $f_{\alpha}(\Gamma_{\infty}^{\lambda})$ and $\Gamma_{\infty}^{\alpha^{-3/4}\lambda}$ have the same distribution.\label{scaling}
   \item (\textbf{Orientation-reversing symmetry}) Let $f$ be the orientation-reversing transformation of $\Gamma_{\infty}^{\lambda}$, that is, $f(\Gamma_{\infty}^{\lambda}):=\{$the reversal of $\ell:\ell\in\Gamma_{\infty}^{\lambda}\}$.  Then $f(\Gamma_{\infty}^{\lambda})$ and $\Gamma_{\infty}^{-\lambda}$ have the same distribution.\label{reversing}
\end{enumerate}
\end{theorem}

Let $\mathbf{P}_{\infty}^{\lambda}$ denote the law of $\Gamma_{\infty}^{\lambda}$.  Throughout the
paper when we declare a loop $\ell\in\Gamma_{\infty}^{\lambda}$, we always assume that $\ell$ is a non-trivial loop unless otherwise specified.  The next theorem gives more properties of $\Gamma_{\infty}^{\lambda}$.  It is well known that in the critical case, $\Gamma_{\infty}^{\lambda=0}$ has the properties (\ref{local})-(\ref{parity}) in this theorem (see \cite[Theorem 2]{CN06} and \cite[Lemma 2.13]{HS23}), but does not have property (\ref{largest}) since a fixed point of the plane is almost surely surrounded by arbitrarily large loops of $\Gamma_{\infty}^{\lambda=0}$.  This implies that $\Gamma_{\infty}^{\lambda}$'s with $\lambda\neq 0$ are singular with respect to $\Gamma_{\infty}^{\lambda=0}$.  Note that (\ref{largest}) was essentially proved in \cite{CJM09}; see Theorem 1.1(3) there.

\begin{theorem}\label{t2}
Let $\lambda\in\mathbb{R}\setminus\{0\}$.  Then $\mathbf{P}_{\infty}^{\lambda}$-almost surely, the loop ensemble $\Gamma_{\infty}^{\lambda}$ has the following properties:
\begin{enumerate}[(i)]
  \item (\textbf{Local finiteness}) For any fixed $\varepsilon,R>0$, there exist finitely many loops of $\Gamma_{\infty}^{\lambda}$ of diameter larger than $\varepsilon$ which intersect the box $[-R,R]^2$. \label{local}
  \item (\textbf{No infinite loops}) $\Gamma_{\infty}^{\lambda}$ has no infinite loops. \label{noinfinite}
  \item (\textbf{Finite chaining}) For any pair of loops $\ell,\ell'\in\Gamma_{\infty}^{\lambda}$, there is a finite set of loops
$\ell_0=\ell,\ell_1,\ldots,\ell_k=\ell'$ in $\Gamma_{\infty}^{\lambda}$ such that, for all $i\in\{1,\ldots,k\},\ell_{i-1}\cap\ell_i\neq\emptyset$. \label{chaining}
  \item (\textbf{No triple points}) No three or more loops can come together at the same point, and a single loop cannot touch the same point more than twice, nor can a loop touch a point where another loop touches itself. \label{triple}
  \item (\textbf{Parity}) Given any pair of loops $\ell,\ell'\in\Gamma_{\infty}^{\lambda}$ with $\ell\cap\ell'\neq\emptyset$, $\ell$ and $\ell'$ have opposite orientation if and only if they are nested (see Section \ref{loop} for the definition of ``nested"). \label{parity}
  \item (\textbf{Existence of outermost loops}) Given any deterministic point $z$ in the plane, all the loops in $\Gamma_{\infty}^{\lambda}$ which surround $z$ form a countable family of nested loops with diameters going to zero, and there is a largest loop in this family. \label{largest}
\end{enumerate}
\end{theorem}

We need a few definitions to state our last result, some of which are borrowed from \cite{CN06}.  Let $\mathbb{D}:=\{z\in\mathbb{C}:|z|<1\}$ be the unit disk centered at 0, and let $\overline{\mathbb{D}}$ denote its closure.  Denote by $\mathbb{D}^{\eta}$ the union of all the hexagons of $\eta\textsf{H}$ contained in $\mathbb{D}$, and let $\partial\mathbb{D}^{\eta}$ be its topological boundary.  Let $\Delta\mathbb{D}^{\eta}$ denote the \textbf{site boundary} of $\mathbb{D}^{\eta}$, that is, the set of $\eta$-hexagons that do not belong to $\mathbb{D}^{\eta}$ but are adjacent to $\eta$-hexagons in $\mathbb{D}^{\eta}$.  A vertex $v$ of $\eta\textsf{H}$ that belongs to $\partial\mathbb{D}^{\eta}$ can be either of two types, according to whether the edge incident on $v$ that is not in $\partial\mathbb{D}^{\eta}$ belongs to a hexagon in $\mathbb{D}^{\eta}$ or not.  We call a vertex of the second type an exposed vertex.  For two distinct points $a,b\in\partial\mathbb{D}$, choose the pair $(a^{\eta},b^{\eta})$ of exposed vertices in $\partial D^{\eta}$ closest to, respectively, $a$ and $b$.  We denote by $\partial_{a,b}\mathbb{D}^{\eta}$ the portion of $\partial\mathbb{D}^{\eta}$ traversed counterclockwise from $a$ to $b$; the remaining part of the boundary is denoted by $\partial_{b,a}\mathbb{D}^{\eta}$.  Similarly, we let $\Delta_{a,b}\mathbb{D}^{\eta}$ denote the portion of $\Delta\mathbb{D}^{\eta}$ whose hexagons are adjacent to $\partial_{a,b}\mathbb{D}^{\eta}$; the remaining part of $\Delta\mathbb{D}^{\eta}$ is denoted by $\Delta_{b,a}\mathbb{D}^{\eta}$.

For a percolation configuration on $\textsf{T}$, a \textbf{boundary path} is a sequence $(e_1,\ldots,e_n)$ of distinct edges of $\textsf{H}$ belonging to the topological boundary of a cluster and such that $e_i$ and $e_{i+1}$ meet at a vertex of $\textsf{H}$ for all $i=1,\ldots,n-1$.  A boundary path is often viewed as a curve in the plane, with its orientation according to the oriented edges in the path.  Recall that the hexagon to the right of an oriented edge is blue.

\begin{figure}
\centering
\includegraphics[height=0.42\textwidth]{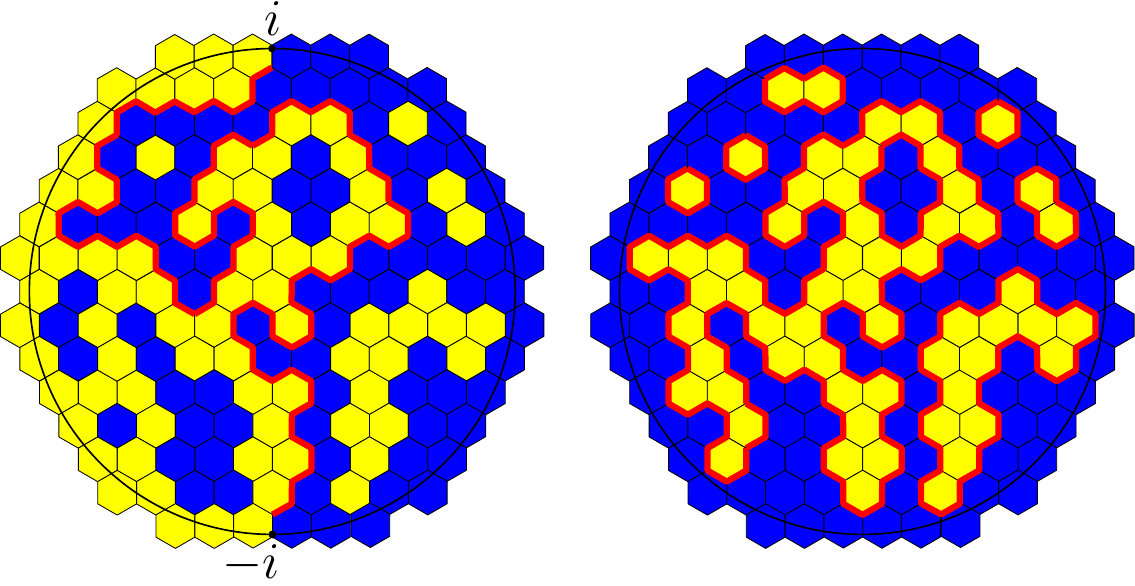}
\caption{\emph{Left}: The percolation exploration path $\gamma_{\eta}^{\lambda}(\mathbb{D},-i,i)$ is indicated by the red line.  \emph{Right}: The collection of cluster boundary loops in $\mathbb{D}$ with monochromatic blue boundary condition, $\Gamma_{\eta,\mathbb{D}}^{\lambda}$, is indicated by the collection of red loops.}\label{interface}
\end{figure}

Given a percolation configuration $\omega_{\eta}^{\lambda}$ on $\mathbb{D}^{\eta}$ and two distinct points $a,b\in\partial\mathbb{D}$, two different choices of boundary conditions will be of interest here:
\begin{itemize}
\item Dobrushin boundary conditions: coloring blue all the hexagons in $\Delta_{a,b}\mathbb{D}^{\eta}$ and yellow all
those in $\Delta_{b,a}\mathbb{D}^{\eta}$.  Then there is a unique boundary path from $a^{\eta}$ to $b^{\eta}$ in $\mathbb{D}^{\eta}$ which separates the blue cluster adjacent to $\Delta_{a,b}\mathbb{D}^{\eta}$ from the yellow cluster adjacent to $\Delta_{b,a}\mathbb{D}^{\eta}$.  We let $\gamma_{\eta}^{\lambda}(\mathbb{D},a,b)$ denote this path, and call it the \textbf{percolation exploration path} from $a^{\eta}$ to $b^{\eta}$ in $\mathbb{D}^{\eta}$.  See the left-hand side of Figure \ref{interface}.
\item Monochromatic blue boundary conditions: coloring blue all the hexagons in $\Delta\mathbb{D}^{\eta}$.  Then we denote by $\Gamma_{\eta,\mathbb{D}}^{\lambda}$ the collection of all the cluster boundary loops in $\mathbb{D}$.  See the right-hand side of Figure \ref{interface}.
\end{itemize}

Relying on Theorem \ref{t1}, we obtain the following near-critical analogs of the SLE$_6$ scaling limit \cite{Smi01,CN07} and the CLE$_6$ scaling limit in a disk \cite{CN06}, though we only establish the convergence here, not deriving the SDE for the limiting curve.  As mentioned in Section \ref{bm}, it was proved in \cite{NW09} that near-critical limiting curves are singular with respect to SLE$_6$ curves.  Moreover, the following theorem says that the near-critical limiting curve can be constructed from $\Gamma_{\infty}^{\lambda}$; analogous result for the critical case was proved right after Theorem 7 in \cite{CN06} (see also \cite[Proposition 2.8]{GPS13}).
\begin{theorem}\label{t3}
Let $\lambda\in\mathbb{R}\setminus\{0\}$.  Consider a coupling of the measures $\{\mathbf{P}_{\eta}^{\lambda}\}_{\eta}$ and $\mathbf{P}_{\infty}^{\lambda}$ on $(\mathcal{L},D_{\mathcal{L}})$ in which $\Gamma_{\eta}^{\lambda}$ converges to $\Gamma_{\infty}^{\lambda}$ almost surely as $\eta\rightarrow 0$.  This can be done by using Theorem \ref{t1} and Skorokhod representation theorem.  Then the following statements hold:
\begin{enumerate}[(i)]
  \item The percolation exploration path $\gamma_{\eta}^{\lambda}(\mathbb{D},a,b)$ converges in probability in the metric (\ref{e28}) to a random curve $\gamma_{\infty}^{\lambda}(\mathbb{D},a,b)$ from $a$ to $b$ in $\overline{\mathbb{D}}$ which can be constructed as a measurable function of $\Gamma_{\infty}^{\lambda}$. \label{sle}
  \item $\Gamma_{\eta,\mathbb{D}}^{\lambda}$ converges in probability in the metric (\ref{e27}) to a loop ensemble $\Gamma_{\infty,\mathbb{D}}^{\lambda}$ on $\overline{\mathbb{D}}$ which can be constructed as a measurable function of $\Gamma_{\infty}^{\lambda}$. \label{cle}
\end{enumerate}
\end{theorem}
\begin{remark}
For simplicity, Theorem \ref{t3} is stated and proved for the unit disk.  With few modifications, our proof also works for any Jordan domain with a piecewise smooth boundary.  Some additional technicalities may be needed if one want to extend the theorem to general domains with rougher boundary.
\end{remark}
\begin{remark}
Combining Theorem \ref{t3} and the powerful topological tool developed by Kemppainen and Smirnov \cite{KS17} (see Theorem 1.5 and Corollary 1.7 in \cite{KS17}, which generalize Aizenman-Burchard tightness criterion \cite{AB99}), it is easy to show that $\gamma_{\infty}^{\lambda}(\mathbb{D},a,b)$ can almost surely be reparametrized by capacity as the trace of a Loewner chain, with its driving process being the limit of the driving process of $\gamma_{\eta}^{\lambda}(\mathbb{D},a,b)$ as $\eta\rightarrow 0$.  Notably, in Section 13.3 of \cite{GPS18a} the authors proposed a conjectural SDE for the limiting driving process.
\end{remark}

\subsection{Strategy of the proof of Theorem \ref{t1}}\label{strategy}
Our proof is inspired by the one in \cite{GPS18b}, where the near-critical coupling $(\omega_{\eta}^{\lambda})_{\lambda\in\mathbb{R}}$ and the enhanced networks were used to prove scaling limit result for the minimal spanning tree on $\textsf{T}$.  \emph{To be consistent with the notation in \cite{GPS18b}, we let $\bar{\lambda}=(\lambda,\lambda')$ be a pair of near-critical parameters with $\lambda<\lambda'$, and we will assume that $\lambda'$ is always fixed and represents $\lambda$ in Theorems \ref{t1}-\ref{t3}, while $\lambda$ in the pair $\bar{\lambda}$ is an auxiliary parameter that will tend to $-\infty$.}

Let us first briefly introduce the enhanced networks since they are powerful tools to help us utilize the quad-crossing scaling limit to study the loop scaling limit.  The enhanced networks, created in \cite{GPS18b} by adding a bit more structure to the networks introduced in \cite{GPS18a}, are planar ``skeletons" of coupled percolation configurations.  Roughly speaking, a network is a graph with vertex set $X$ given by the $\varepsilon$-important points (see Section \ref{pivotal}) in the configuration $\omega^{\lambda}$ which switch between parameters $\lambda$ and $\lambda'$, and edges given by the open connections in $\omega^{\lambda}\setminus X$.  Compared with the network, an enhanced network adds the information about the colors of the points of $X$ in $\omega^{\lambda}$ and which of these points are connected together by an open cluster of $\omega^{\lambda}\setminus X$; it has two types of vertices, representing those open clusters and colored $\varepsilon$-important points, respectively, and its edges represent open connections in $\omega^{\lambda}\setminus X$ between the two types of vertices.  When $\lambda$ is very negative and much smaller than $\lambda'$, then all the $\lambda$-clusters in a bounded domain are of small macroscopic size, and from the enhanced network we know how the $\lambda$-clusters merge to form the $\lambda'$-clusters when we raise the near-critical parameter, which allows us to describe the macroscopic structure of $\omega^{\lambda'}$.

Let us now outline the proof of Theorem \ref{t1}.  Our method involves three steps:

\textbf{Step 1. Convergence of good embeddings of enhanced networks.}  In \cite{GPS18b}, the stabilization and $\eta$-convergence of enhanced networks were established, and the cut-off trees are extracted from enhanced networks in which the embedded edges are straight line segments.  Since we are studying scaling limits in loop space, we need to construct some ``good" embeddings of enhanced networks, and prove that such embeddings exist both in the discrete and continuum cases with high probability, with these embeddings also satisfying the stabilization and $\eta$-convergence.

\textbf{Step 2. Construction of the cut-off loop ensembles.} From the good embeddings of enhanced networks in Step 1, we build the router graphs (see Figure \ref{routerg}) and extract the cut-off loop ensembles $\Gamma^{\bar{\lambda},\varepsilon}$ (see Figure \ref{cutoffloop}) as the collections of boundary loops of the router graphs.  Under the coupling, the discrete and continuum cut-off loop ensembles $\Gamma_{\eta}^{\bar{\lambda},\varepsilon}$ and $\Gamma_{\infty}^{\bar{\lambda},\varepsilon}$ are equal with high probability.  Moreover, for discrete percolation we build the good router graphs (see Figure \ref{routerg2}), which can be viewed as skeletons embedded into $\lambda'$-open clusters; we extract the good cut-off loop ensembles $\tilde{\Gamma}_{\eta}^{\bar{\lambda},\varepsilon}$ (see Figure \ref{cutoffloop2}) from such graphs, and show that they approximate $\Gamma_{\infty}^{\bar{\lambda},\varepsilon}$ well for all $\eta>0$ small enough.

\textbf{Step 3. Approximation of $\Gamma_{\eta}^{\lambda'}$ by the good cut-off loop ensembles $\tilde{\Gamma}_{\eta}^{\bar{\lambda},\varepsilon}$.}
Using chains of $\lambda$-open clusters with the ``linking" ($\lambda$-closed) sites being $\varepsilon$-important points, we prove that for all $\eta>0$ small enough, $\Gamma_{\eta}^{\lambda'}$ is well-approximated by $\tilde{\Gamma}_{\eta}^{\bar{\lambda},\varepsilon}$ if $\lambda$ is very negative and $\varepsilon$ is small.  Combined with Step 2, we deduce that $\Gamma_{\eta}^{\lambda'}$ is close to $\Gamma_{\infty}^{\bar{\lambda},\varepsilon}$, which implies the existence and uniqueness of the scaling limit $\Gamma_{\infty}^{\lambda'}$.  The following is a pictorial sketch of this strategy (see also the analogous pictorial summary of the strategy in \cite[Section 1.4]{GPS18b}):
\begin{figure}[H]
\centering
\includegraphics[height=0.11\textwidth]{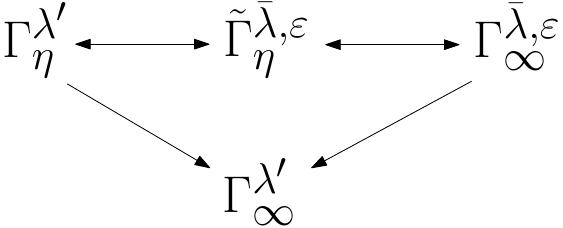}
\label{pictorial}
\end{figure}

In fact, to avoid boundary issues, similarly to \cite{GPS18b} we shall first work on a torus, and then extend the results to the full plane.

Note that the scaling limit of the minimal spanning trees was constructed in \cite{GPS18b} by taking $\lambda\rightarrow-\infty,\lambda'\rightarrow\infty$ and then letting $\varepsilon\rightarrow 0$, while our near-critical loop scaling limit is constructed by fixing $\lambda'$ and taking $\lambda\rightarrow-\infty$ and then letting $\varepsilon\rightarrow 0$.  Although our strategy is analogous to that of \cite{GPS18b}, there are new difficulties and technicalities: First, we need to construct good embeddings of enhanced networks because of the loop topology.  Second, to approximate cluster boundaries we need to build good cut-off loops and handle many topological issues by using chains of clusters.  Third, as we mentioned above, in \cite{GPS18b} the parameter $\lambda'$ is going to $\infty$ and the system will be very supercritical, but here $\lambda'$ is fixed.

Recall that based on Cardy's formula proved in \cite{Smi01}, the convergence of critical percolation exploration paths was established in \cite{Smi01,CN07}, which was employed to prove the CLE$_6$ scaling limit result in \cite{CN06}.  Since a near-critical Cardy theorem was established in \cite[(9.2)]{GPS18a}, one may wonder if the method of \cite{CN06} can be used directly to show Theorem \ref{t1}.  However, their method relies on the conformal invariance of the critical scaling limit, while in the near-critical case such property does not hold.

Instead of the strategy of the present paper, one could try to prove Theorem \ref{t1} without using the near-critical coupling.  Our method has the following advantages: First, it keeps track of how the $\lambda$-open clusters merge to form the $\lambda'$-open clusters, from which one can extract the ``merging points" together with their parameters (i.e., merging times).  Thus, it can be applied to study the scaling limit of invasion percolation in loop space (see Section \ref{invasion}).  Second, unlike \cite{CN06}, it proves the full scaling limit result directly, not relying on the single-interface scaling limit result for general simply connected domains that may have fractal non-simple boundaries, which allows us to avoid the technicalities related to boundary issues.

\subsection{Applications and related models}\label{app}
In this subsection, we wish to list a few applications of the results of the present paper, and discuss some related models and open problems.

\subsubsection{Scaling limit of the near-critical coupling}
It was proved in \cite[Theorem 1.5]{GPS18a} that as $\eta\rightarrow 0$, the process $\lambda\mapsto\omega_{\eta}^{\lambda}$ converges in law, in a
Skorokhod space of c\`{a}dl\`{a}g paths with values in the quad-crossing space, to a limiting random process $\lambda\mapsto\omega_{\infty}^{\lambda}$.  Similarly, it is expected that as $\eta\rightarrow 0$, the process $\lambda\mapsto\Gamma_{\eta}^{\lambda}$ converges in law, in a Skorokhod space of c\`{a}dl\`{a}g paths with values in $(\mathcal {L},D_\mathcal {L})$, to a limiting random process $\lambda\mapsto\Gamma_{\infty}^{\lambda}$.  It is not hard to show this by using Theorem \ref{t1} and its proof; we plan to prove this together with some properties of the limiting process in a future work.

\subsubsection{Applications to first-passage percolation}
First-passage percolation (FPP) is a model for the spread of a fluid in a random medium (see \cite{ADH17} for a recent survey).  \textbf{Bernoulli (site) FPP} is equivalent to Bernoulli site percolation, but studied from the FPP point of view: For every site $v$ we let $t(v)=0$ if $v$ is open and $t(v)=1$ if $v$ is closed, where $t(v)$ is the passage time of $v$.  One of the main subjects of FPP is studying the asymptotics of first-passage times between two sets that are far from each other.  In \cite{JY19,Yao14,Yao18,Yao19,Yao24}, scaling limits of critical and near-critical percolation are used to give exact asymptotics for first-passage times of Bernoulli FPP.

\begin{figure}
\centering
\includegraphics[height=0.22\textwidth]{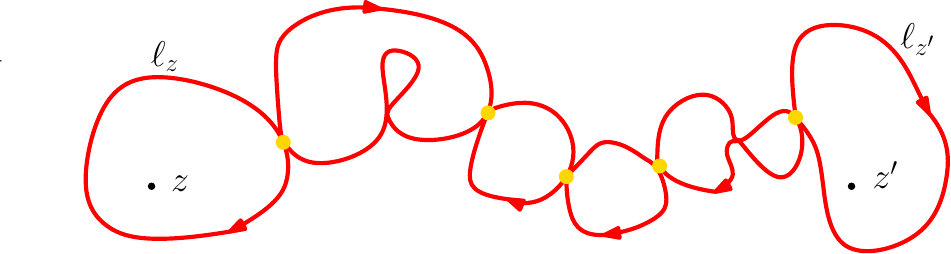}
\caption{Illustration of a chain $\mathscr{L}$ of loops from $\ell_z$ to $\ell_{z'}$ in $\Gamma_{\infty}^{-1}$ such that $T_{z,z'}=|\mathscr{L}|-1=5$.  Yellow dots represent touching points of consecutive loops in $\mathscr{L}$.  Note that the loops in $\mathscr{L}$ are outermost loops in $\Gamma_{\infty}^{-1}$, which are oriented clockwise.}\label{fpploopchain}
\end{figure}

In \cite{Yao24}, we use the near-critical scaling limit of the collection of all open clusters (viewed as random compact
sets) to study Bernoulli FPP on $\textsf{T}$ when $p\uparrow p_c$.  See Figure 4 in \cite{Yao24} for the ``point-to-point" passage time defined by chains of ``continuum clusters".   By Theorem \ref{t2} (\ref{chaining}) and (\ref{largest}), we have a very natural way to define the point-to-point passage time between two points $z$ and $z'$ in the plane for $\Gamma_{\infty}^{-1}$, thus constructing a ``continuum" FPP:
\begin{equation*}
T_{z,z'}:=\inf\{|\mathscr{L}|-1:\mathscr{L}\mbox{ is a chain of loops from $\ell_z$ to $\ell_{z'}$ in $\Gamma_{\infty}^{-1}$}\},
\end{equation*}
where $\ell_z$ denotes the outermost loop of $\Gamma_{\infty}^{-1}$ surrounding $z$.  See Figure \ref{fpploopchain} for an example.

\subsubsection{Applications to invasion percolation}\label{invasion}
Invasion percolation is a paradigmatic model of self-organized criticality \cite{BTW87}, and is closely related to critical and near-critical Bernoulli percolation.  It has been extensively studied in the probability and statistical physics communities; see, for example, \cite{CCN85,DS11,DS12,DSV09,GPS18b,Jar03} for invasion percolation on lattices.

Garban, Pete and Schramm \cite{GPS18b} have constructed the scaling limit of \textbf{invasion percolation tree} (with the constraint that invaded edges are forbidden to form circuits) on $\textsf{T}$ in the space of essential spanning forests introduced in \cite{ABNW99}.  Using the method developed in the present paper, one may prove that the collection of all the boundary loops of the so-called \textbf{invasion percolation cluster} (IPC) converges to a loop-ensemble scaling limit, as the mesh size tends to zero.  This scaling limit may be useful in improving some results of \cite{DS11,DS12}.  Moreover, based on the singularity result in \cite{DSV09}, it is not hard to show that such scaling limit of IPC and the scaling limit of incipient infinite cluster constructed in \cite{Yao18b} are mutually singular.

\subsubsection{Applications to chemical distance}
Very recently, \DJ ankovi\'{c}, Markering, Miller and Yuan \cite{DMMY} proved that the intrinsic (i.e., chemical distance) metric on the critical open clusters on $\textsf{T}$ converges in the scaling limit to the geodesic CLE$_6$ metric, and established the existence of the chemical distance exponent.  It is expected that analogous results also hold in the near-critical case, with the limiting geodesic metric associated to $\Gamma_{\infty}^{\lambda}$ or
$\Gamma_{\infty,\mathbb{D}}^{\lambda}$.  Such results seem to have applications to the following conjecture made by Benjamini \cite{Ben10}:  Suppose $p>p_c$.  Condition that the origin $0$ is in the infinite cluster.  Consider the balls centered at 0 with radius $n$ for the chemical distance on the infinite cluster.  Then the limiting shape of these balls becomes round as $p\downarrow p_c$.

\subsubsection{A discussion on the equivalence of near-critical scaling limits in quad-crossing space and loop space}\label{equivalence}
There are several natural notions of the ``full" scaling limit for percolation in the plane that encode all macroscopic connectivity information, and the general belief is that these different notions are all equivalent in the sense that the associated $\sigma$-algebras are the same; see Section 1.1 of \cite{SS11} for a nice discussion.  In the critical case, Camia and Newman \cite{CN06} essentially proved that CLE$_6$ determines the quad-crossing scaling limit, as pointed out in \cite{GPS13}; recently, Holden and Sun \cite{HS23} proved the measurability in the reverse direction, thus establishing the equivalence of these two scaling limits.  When $\lambda\neq 0$, it is expected that $\Gamma_{\infty}^{\lambda}$ and $\omega_{\infty}^{\lambda}$ (i.e., the near-critical quad-crossing scaling limit) are also equivalent.  Similarly to the critical case, by using the arguments in Section 2.3 of \cite{GPS13} and an extension of Theorem \ref{t3} to piecewise smooth quads, one can show that $\Gamma_{\infty}^{\lambda}$ almost surely determines $\omega_{\infty}^{\lambda}$.  However, the proof of \cite[Theorem 6.10]{HS23} is not directly applicable to show the measurability in the reverse direction, since it relies on the conformal invariance of critical scaling limit (which allows the authors to extend measurability result for a disk to SLE$_6$ cut-out domains), while near-critical scaling limits are not conformally invariant.

\subsubsection{A discussion on the scaling limits of near-critical planar FK-Ising model}\label{FK}
In the monotone coupling of near-critical percolation in Definition \ref{def1}, sites have independent exponential clocks and turn open in a Poissonian way, and near-critical behavior is governed by pivotal sites at criticality, due to near-critical stability.  However, as explained in \cite{DCGP14}, in any monotone coupling of planar FK-Ising configurations with parameters $p$, as one raises $p$ near $p_c$, the new edges appear simultaneously in ``clouds", in a self-organized way, so that the correlation length defined via crossing probabilities is not governed anymore by the number of pivotal edges at criticality.  Driven by this phenomenology, Duminil-Copin and Manolescu \cite{DCM22} introduced a crucial quantity called ``mixing rate" that encodes how much an edge is influenced by boundary conditions at a distance $R$, which plays a role in FK-Ising model similar to the role played by the probability of being pivotal in percolation.  Note that it was proved in \cite{GK} that any subsequential scaling limit of the collection of FK-Ising cluster boundaries is singular with respect to CLE$_{16/3}$, as mentioned in Section \ref{bm}.

\subsection{Organization of the paper}
In Section \ref{pre}, we introduce the basic notation and background on the quad-crossing and loop-ensemble spaces, arm events, Poisson point process $\textsf{PPP}_{\bar{\lambda}}^{\varepsilon}$ of pivotal points, enhanced network $\textnormal{\textsf{EN}}_M^{r\text{-}\textnormal{\textsf{meso}}}(\omega^{\lambda},\textnormal{\textsf{PPP}}_{\bar{\lambda}}^{\varepsilon})$ and correlation length, collect some classical results about critical and near-critical percolation, and prove a property of large clusters related to their pivotals in $\textsf{PPP}_{\bar{\lambda}}^{\varepsilon}$.  In Section \ref{clusterchain}, we give some geometric properties for chains of clusters.  In Section \ref{routergraph}, we use some ``good embedding" of enhanced network, called the router graph, to define the cut-off loop ensemble $\Gamma_{\infty}^{\bar{\lambda},\varepsilon}$ for $\omega_{\infty}^{\lambda}$; we also define the good cut-off loop ensemble $\tilde{\Gamma}_{\eta}^{\bar{\lambda},\varepsilon}$ for $\omega_{\eta}^{\lambda}$, and prove that $\tilde{\Gamma}_{\eta}^{\bar{\lambda},\varepsilon}$ is well-approximated by $\Gamma_{\infty}^{\bar{\lambda},\varepsilon}$.  In Section \ref{mainproof}, by using cluster chains we show that $\Gamma_{\eta}^{\lambda'}$ is well-approximated by $\tilde{\Gamma}_{\eta}^{\bar{\lambda},\varepsilon}$, and then prove our main result, Theorem \ref{t1}.  Based on Theorem \ref{t1}, we prove Theorems \ref{t7} and \ref{t2} in Section \ref{property}, and prove Theorem \ref{t3} in Section {\ref{disk}}.

\section{Notation and preliminaries}\label{pre}

Let $\mathbb{N}=\{1,2,\ldots\}$ denotes the set of natural integers.  We identify the plane
$\mathbb{R}^2$ with the set $\mathbb{C}$ of complex numbers in the usual way.

For $M>0$, let $\Lambda_M=[-M,M]^2$ denote the square centered at the origin, and let $\mathbb{T}_M^2$ denote the \textbf{torus} of side-length $2M$, which is realized as the subdomain $[-M,M)^2$ of $\mathbb{C}$ with a periodic boundary condition.  For
$\varepsilon>0$ and $z\in\mathbb{C}$, we usually let $B_{\varepsilon}(z)$ denote the $\varepsilon$-square $[i,i+1)\varepsilon\times[j,j+1)\varepsilon$ in the lattice $\varepsilon\mathbb{Z}^2$ that contains $z$.  We often call the squares $\Lambda_M$ and $B_{\varepsilon}(z)$ boxes, while we tend to use $B_{\varepsilon}(z)$ for boxes much smaller than $\Lambda_M$, and often consider a collection of such boxes with their union covering $\Lambda_M$.

For $0<r<R$, define the annulus $A(r,R):=\Lambda_R\setminus\Lambda_r$.  For $z\in\mathbb{C}$, set $A(z;r,R):=z+A(r,R)$.

Let $\|\cdot\|_{\infty}$ denote the $L^{\infty}$-norm on $\mathbb{C}$.  For a set $A\subset\mathbb{C}$, let $\diam(A)$ denote the $L^{\infty}$-diameter of $A$.  For convenience, we will just say diameter to mean $L^{\infty}$-diameter in the rest of this paper.  Similarly, for two sets $A,B\subset\mathbb{C}$, when we say $A$ is at distance $r$ from $B$, we mean that their $L^{\infty}$-distance is $r$.

For a bounded set $A\subset\mathbb{C}$, we write $\overline{A}$ for the closure of $A$, and denote by $F(A)$ the \textbf{filling} of $A$, that is, the complement of the unbounded connected component of $\mathbb{C}\setminus A$.

A \textbf{path} on $\textsf{T}$ is a finite sequence $(v_1,\ldots,v_n)$ of sites such that $v_i$
and $v_{i+1}$ are neighbors for all $i=1,\ldots,n-1$.  A path will be called open or closed if all its sites are open or closed, respectively.  A \textbf{circuit} is a path on $\textsf{T}$ that starts and ends at the same site.  A path (resp. circuit) is sometimes identified with a curve (resp. loop) formed by its edges in the plane.

An open cluster of $\omega_{\eta}^{\lambda}$ is often called a $\lambda$-cluster.

Throughout this paper, $\varepsilon,\delta,C,K,C_i$ and $K_i$ stand for (strictly) positive
constants that may change from line to line according to the context.

\subsection{Loop-ensemble space}\label{loop}
In dealing with the scaling limits, topological spaces of curves and loops were introduced in \cite{AB99,CN06}.  We borrow the notation and definitions from \cite{CN06} with small modifications.

When taking the scaling limit of percolation on the whole plane, we compactify $\mathbb{C}$ into $\hat{\mathbb{C}}:=
\mathbb{C}\cup\{\infty\}$ as follows.  First, we replace the Euclidean metric with a distance function
$\Delta(\cdot,\cdot)$ defined on $\mathbb{C}\times\mathbb{C}$ by
\begin{equation}\label{e35}
\Delta(x,y):=\inf_{\varphi}\int\frac{2}{1+|\varphi(s)|^2}ds,
\end{equation}
where the infimum is over all smooth curves $\varphi$ joining $x$
with $y$, parameterized by arc length, and $|\cdot|$ denotes
the Euclidean norm.  Then, we add a single point at
infinity to get the compact space $\hat{\mathbb{C}}$ which is
isometric, via stereographic projection, to the unit sphere.  The metric $\Delta$ is equivalent to the Euclidean metric
in bounded regions, while the $\Delta$-distance between any two points outside $[-M,M]^2$ in $\hat{\mathbb{C}}$ is at
most $O(1/M)$.

A \textbf{loop} is a continuous map from the circle $\mathbb{R}/\mathbb{Z}$ to $\hat{\mathbb{C}}$ or $\mathbb{T}_M^2$ identified
up to reparametrizations by orientation-preserving homeomorphisms of the circle.  The pseudo-distance $D$ between two loops $\ell,\ell'$ on $\hat{\mathbb{C}}$, induced by $\Delta$, is defined by
\begin{equation*}
D(\ell,\ell'):=\inf_{\phi}\sup_{t\in\mathbb{R}/\mathbb{Z}}\Delta(\ell(t),\ell'(\phi(t))),
\end{equation*}
where the infimum is taken over all orientation-preserving homeomorphisms ${\phi:\mathbb{R}/\mathbb{Z}\rightarrow\mathbb{R}/\mathbb{Z}}$.

A closed set of loops on $\hat{\mathbb{C}}$ with respect to the $D$-metric is called
a \textbf{loop ensemble on $\hat{\mathbb{C}}$}.  Let $\mathcal {L}$ be the space of loop ensembles on $\hat{\mathbb{C}}$.  The distance
between two loop ensembles $\mathcal {F},\mathcal {F}'$ in $\mathcal {L}$ is defined by the Hausdorff metric induced by $D$:
\begin{equation}\label{e26}
D_{\mathcal {L}}(\mathcal {F},\mathcal {F}'):=\inf\{\varepsilon>0:\forall
\ell\in \mathcal {F},\exists \ell'\in\mathcal {F}'\mbox{ such
that }D(\ell,\ell')\leq\varepsilon\mbox{ and vice
versa}\}.
\end{equation}
Note that the space $\mathcal {L}$ equipped with the metric $D_{\mathcal {L}}$ is a complete separable metric space.

A loop whose range is a single point in $\hat{\mathbb{C}}$ is called a \textbf{trivial loop}.  In Theorem \ref{t1}, to make the collection $\Gamma_{\eta}^{\lambda}$ of all cluster boundary loops be a closed set so that it is an element of $\mathcal {L}$, we need to add to $\Gamma_{\eta}^{\lambda}$ the trivial loop at the point $\{\infty\}$.  Throughout the paper when we declare a loop $\ell$ to be in $\Gamma_{\eta}^{\lambda}$ or $\Gamma_{\infty}^{\lambda}$, we always assume that $\ell$ is a non-trivial loop unless otherwise stated.

A loop on $\mathbb{C}$ is called a \textbf{quasi-simple loop} if there exists a sequence of simple loops on $\mathbb{C}$ converging to it in the metric $D$.  Let $\ell$ be a non-trivial quasi-simple loop.  It is clear that the winding number of $\ell$ around a point $z\in\mathbb{C}\setminus\ell$ is equal to 0 or $\pm 1$.  In the latter case, we say that $z$ is \textbf{surrounded} by $\ell$.  Let $\textrm{reg}(\ell)$ denote the region enclosed by $\ell$, that is, the closure of the union of all the points surrounded by $\ell$.  Given two distinct non-trivial quasi-simple loops $\ell$ and $\ell'$ on $\mathbb{C}$, we say that $\ell$ and $\ell'$ are \textbf{nested} if and
only if $\ell\subset\textrm{reg}(\ell')$ or $\ell'\subset\textrm{reg}(\ell)$.

Instead of the metric $D$, when we work in a bounded region of the plane, it is more convenient to use the metric $d$ between two loops $\ell,\ell'$ on $\mathbb{C}$ induced by $\|\cdot\|_{\infty}$:
\begin{equation}\label{e12}
d(\ell,\ell'):=\inf_{\phi}\sup_{t\in\mathbb{R}/\mathbb{Z}}\|\ell(t)-\ell'(\phi(t))\|_{\infty},
\end{equation}
where the infimum is taken over all orientation-preserving homeomorphisms ${\phi:\mathbb{R}/\mathbb{Z}\rightarrow\mathbb{R}/\mathbb{Z}}$.  For two sets of loops $\mathcal {F},\mathcal {F}'$ on $\mathbb{C}$, we define the following distance by the Hausdorff metric induced by $d$:
\begin{equation}\label{e27}
d_{\mathcal {L}}(\mathcal {F},\mathcal {F}'):=\inf\{\varepsilon>0:\forall
\ell\in \mathcal {F},\exists \ell'\in\mathcal {F}'\mbox{ such
that }d(\ell,\ell')\leq\varepsilon\mbox{ and vice
versa}\}.
\end{equation}
A closed set of loops on a bounded closed set $S$ of the plane with respect to the $d$-metric is called
a \textbf{loop ensemble on $S$}.  Let $\mathcal {L}_S$ denote the space of loop ensembles on $S$.  The space $\mathcal {L}_S$ equipped with the metric $d_{\mathcal {L}}$ is also a complete separable metric space, and the metrics $d_{\mathcal {L}}$ and $D_{\mathcal {L}}$ are equivalent on $\mathcal {L}_S$.

A \textbf{curve} is a continuous map from the unit interval $[0,1]$ to the plane $\mathbb{C}$ identified
up to reparametrizations by increasing homeomorphisms of $[0,1]$.  So, a loop on $\mathbb{C}$ can be viewed as a curve modulo increasing reparametrization and rerooting.  With a slight abuse of notation, the metric $d$ (the same notation is used in (\ref{e12})) between two curves $\gamma$ and $\gamma'$, induced by $\|\cdot\|_{\infty}$, is defined by
\begin{equation}\label{e28}
d(\gamma,\gamma'):=\inf_{\phi}\sup_{t\in[0,1]}\|\gamma(t)-\gamma'(\phi(t))\|_{\infty},
\end{equation}
where the infimum is taken over all increasing bijections $\phi$ from $[0,1]$ onto itself.  In this paper, we will always assume that the meaning of $d(\gamma,\gamma')$ is as in (\ref{e12}) if $\gamma$ and $\gamma'$ are loops, and is as in (\ref{e28}) if $\gamma$ and $\gamma'$ are curves.

\subsection{Quad-crossing space}\label{quad}
To describe the scaling limits of planar percolation, Schramm and
Smirnov \cite{SS11} introduced the quad-crossing space, which is developed further in \cite{GPS13,GPS18a,GPS18b}.  Roughly speaking, a quad-crossing scaling limit encodes all macroscopic quad-crossing events; the quads are ``building blocks" for more complicated geometric structures.  In this section we briefly recall the notation and definitions for the quad-crossing topology, borrowed from \cite{GPS18b}, together with some basic results.

Let $D\subset\hat{\mathbb{C}}$ be open, or be equal to the torus $\mathbb{T}_M^2$.  A \textbf{quad} in $D$ is a homeomorphism ${Q:[0,1]^2\rightarrow D}$.  The space of all quads in $D$, denoted by $\mathcal
{Q}_D$, can be equipped with the following metric:
\begin{equation*}
d_{\mathcal
{Q}}(Q_1,Q_2):=\inf_{\phi}\sup_{z\in\partial[0,1]^2}|Q_1(z)-Q_2(\phi(z))|,
\end{equation*}
where the infimum is over all homeomorphisms $\phi:
[0,1]^2\rightarrow[0,1]^2$ which preserve the four corners of the
square.  A \textbf{crossing} of a quad $Q$ is a connected closed
subset of $[Q]:=Q([0,1]^2)$ that intersects both $\partial_1Q=Q(\{0\}\times[0,1])$
and $\partial_3Q=Q(\{1\}\times[0,1])$.  We say that $Q$ has a \textbf{dual crossing} between $\partial_1Q$ and $\partial_3Q$ by some closed subset
$S\subseteq [Q]$ if there is no crossing in $S$ between $\partial_2 Q=Q([0,1]\times\{0\})$ and $\partial_4 Q=Q([0,1]\times\{1\})$.

We will sometimes view a quad $Q$ as its range $[Q]$ together with four ordered marked points $Q(0,0),Q(1,0),Q(1,1),Q(0,1)$.

There is a natural partial order on $\mathcal {Q}_D$ induced by the crossings.  We write
$Q_1\leq Q_2$ if any crossing of $Q_2$ contains a crossing of $Q_1$.
Furthermore, we write $Q_1<Q_2$ if there are open neighborhoods
$\mathcal {N}_i$ of $Q_i$ (in the uniform metric), such that $N_1\leq N_2$
holds for any $N_i\in\mathcal {N}_i$.  We say that a subset
$\mathcal {S}\subset\mathcal {Q}_D$ is \textbf{hereditary} if,
whenever $Q\in \mathcal {S}$ and $Q'\in\mathcal {Q}_D$ satisfies
$Q'<Q$, we have $Q'\in \mathcal {S}$.  The collection of all closed
hereditary subsets of $\mathcal {Q}_D$ will be denoted by
$\mathscr{H}_D$.

Note that any near-critical percolation configuration $\omega_{\eta}^{\lambda}$, considered as a
union of the topologically closed percolation-wise open hexagons in $D$, naturally defines an element of
$\mathscr{H}_D$: the set of all quads for which $\omega_{\eta}^{\lambda}$ contains a crossing.  With a slight abuse of notation, we will still
denote by $\omega_{\eta}^{\lambda}$ the element in $\mathscr{H}_D$ induced by the percolation configuration $\omega_{\eta}^{\lambda}$, and by $\mathbf{P}_{\eta}^{\lambda}$ the induced probability measure on $\mathscr{H}_D$.

By introducing a natural topology, $\mathscr{H}_D$ can be made into
a compact metric space.  Indeed, let
\begin{align*}
&\boxminus_{Q}:=\{\mathcal {S}\in\mathscr{H}_D: Q\in\mathcal {S}\}\qquad\mbox{for any }Q\in\mathcal {Q}_D,\\
&\boxdot_{\mathcal {U}}:=\{\mathcal {S}\in\mathscr{H}_D: \mathcal
{U}\cap\mathcal {S}=\emptyset\}\qquad\mbox{for any open subset
}\mathcal {U}\subset\mathcal {Q}_D.
\end{align*}
Then we define $\mathcal{T}_D$ to be the minimal topology on $\mathscr{H}_D$ that contains every $\boxminus_{Q}^c$
and $\boxdot_{\mathcal {U}}^c$ as open sets.  It is proved in
\cite[Theorem 1.13]{SS11} that for any nonempty open $D$, the topological space
$(\mathscr{H}_D,\mathcal{T}_D)$ is compact, Hausdorff and metrizable.  Furthermore, for any dense $\mathcal {Q}_0\subset\mathcal
{Q}_D$, the events $\{\boxminus_{Q}:Q\in\mathcal {Q}_0\}$ generate
the Borel $\sigma$-field of $\mathscr{H}_D$.  An arbitrary metric
generating the topology $\mathcal{T}_D$ will be denoted by
$d_{\mathscr{H}}$.  Moreover, the space of Borel probability
measures of $(\mathscr{H}_D,\mathcal{T}_D)$, equipped with the weak*-topology is a compact metrizable Hausdorff space; see \cite[Corollary 1.15]{SS11}.  This implies that $\mathbf{P}_{\eta}^{\lambda}$ has subsequential scaling limits on $\mathscr{H}_D$ as $\eta\rightarrow 0$.  The following theorem shows the uniqueness of the near-critical scaling limit.

\begin{theorem}[Theorem 1.4 and Corollary 10.5 in \cite{GPS18a}]\label{t4}
For any $\lambda\in\mathbb{R}$, there is a unique probability measure $\mathbf{P}_{\infty}^{\lambda}$ for percolation configurations $\omega_{\infty}^{\lambda}$ in
$(\mathscr{H}_D,\mathcal{T}_D)$, such that $\omega_{\eta}^{\lambda}$ converges in law to $\omega_{\infty}^{\lambda}$ as $\eta\rightarrow 0$.  Moreover, $\mathbf{P}_{\infty}^{\lambda}$ has a scaling covariance property: For any $\alpha>0$, let $f_{\alpha}(z):=\alpha z$ for $z\in \mathbb{C}$.  Then $f_{\alpha}(\omega_{\infty}^{\lambda})$ and $\omega_{\infty}^{\alpha^{-3/4}\lambda}$ have the same law.
\end{theorem}

In order to state a result for quad-crossing events, we need to introduce a family of polygonal quads in $\mathbb{T}_M^2$.
\begin{definition}[A dyadic family of quads in $\mathbb{T}_M^2$]\label{def9}
For any $M,k\geq 1$, let $\mathcal{Q}_M^k=\{Q_n^k\}_{1\leq n\leq N_{M,k}}$ be the family
of all polygonal quads in $\mathbb{T}_M^2$ such that their boundary $\partial Q_n^k$ consists of edges of
$\mathbb{T}_M^2\cap 2^{-k}\mathbb{Z}^2$ and the four marked points are vertices of $\mathbb{T}_M^2\cap 2^{-k}\mathbb{Z}^2$.
Notice that $\mathcal{Q}_M^k\subset\mathcal{Q}_M^{k+1}$.
\end{definition}

The following lemma shows the convergence of the $\mathbf{P}_{\eta}^{\lambda}$-probability of the intersection of any finite set of quad-crossing and non-crossing events.  Note that (\ref{e13}) is
\cite[Lemma 5.1]{SS11} in case $\lambda=0$, and the same proof works for general $\lambda$; see \cite[Lemma 2.3]{GPS18b}. (\ref{e14}) is a direct consequence of (\ref{e13}) and Theorem \ref{t4}; see also \cite[Corollary 5.2]{SS11} and \cite[Lemma 2.5]{GHSS21} for similar results in case $\lambda=0$.

\begin{lemma}\label{l1}
Fix $\lambda\in\mathbb{R}$.  For any fixed $Q\in\mathcal {Q}_{D}$, we have
\begin{equation}\label{e13}
\mathbf{P}_{\infty}^{\lambda}[\partial\boxminus_{Q}]=0.
\end{equation}
Therefore, by the weak convergence of $\mathbf{P}_{\eta}^{\lambda}$ to $\mathbf{P}_{\infty}^{\lambda}$,
\begin{equation}\label{e59}
\mathbf{P}_{\eta}^{\lambda}[\boxminus_{Q}]\rightarrow \mathbf{P}_{\infty}^{\lambda}[\boxminus_{Q}]\quad\mbox{as $\eta\rightarrow 0$}.
\end{equation}
Let $\mathbf{P}$ be a coupling of the measures $\{\mathbf{P}_{\eta}^{\lambda}\}$ and $\mathbf{P}_{\infty}^{\lambda}$ on $(\mathscr{H}_D,\mathcal{T}_D)$ in which $\omega_{\eta}^{\lambda}\rightarrow\omega_{\infty}^{\lambda}$ a.s. as $\eta\rightarrow 0$.  For any event $\mathcal{E}$ that can be expressed using finitely many of the events $\boxminus_{Q},Q\in\mathcal {Q}_{D}$ and the operations of union, intersection and taking complements, we have
\begin{equation}\label{e14}
\mathbf{P}[\{\omega_{\eta}^{\lambda}\in\mathcal{E}\}\Delta\{\omega_{\infty}^{\lambda}\in\mathcal{E}\}]\rightarrow 0\quad\mbox{as }\eta\rightarrow 0.
\end{equation}
In particular, given any $M,k,n\geq 1$, for $D=\mathbb{T}_M^2$ and any $Q_1,\ldots,Q_n\in\mathcal{Q}_M^k$, we have
\begin{equation}\label{e4}
\mathbf{P}[\{\omega_{\eta}^{\lambda}\in\cap_{i=1}^n\boxminus_{Q_i}\}\Delta\{\omega_{\infty}^{\lambda}\in\cap_{i=1}^n\boxminus_{Q_i}\}]\rightarrow 0\quad\mbox{as }\eta\rightarrow 0.
\end{equation}
\end{lemma}

\begin{figure}
\centering
\includegraphics[height=0.33\textwidth]{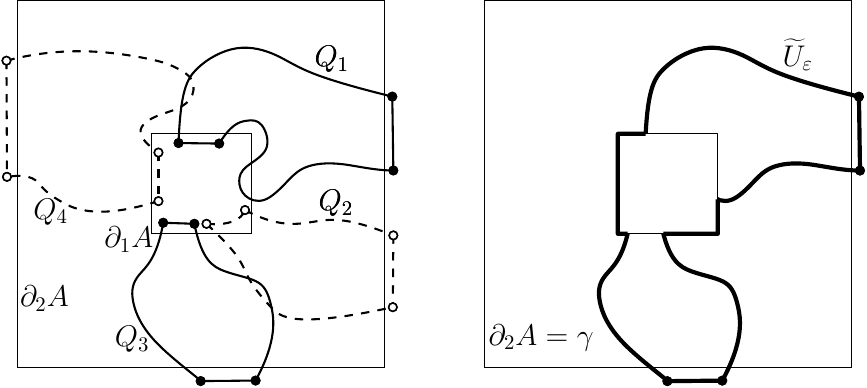}
\caption{Defining the alternating 4-arm event in Section \ref{quad} and the color of a pivotal point in
Section \ref{pivotal} using quad-crossings. Quads with a primal crossing are indicated by heavy solid lines; quads with a dual
crossing are indicated by dashed lines.
}\label{colored4arm}
\end{figure}

Now let us define the alternating 4-arm events in a way that makes them measurable w.r.t. the quad-crossing topology; general arm events for quad-crossing space were introduced and studied in \cite{GPS13}.  For any topological annulus $A\subset D$ with piecewise smooth inner and outer boundary pieces $\partial_1 A$ and $\partial_2 A$ (and for the case of $D=\mathbb{T}_M^2$, we also require $A$ to be null-homotopic), we define the
\textbf{alternating 4-arm event} in $A$ as $\mathcal{A}_4=\mathcal{A}_4(A)=\bigcup_{\delta>0}\mathcal{A}_4^{\delta}$, where $\mathcal{A}_4^{\delta}$ is the existence of quads $Q_i\subset D,i=1,2,3,4$, with the following properties (see the left-hand side of Figure \ref{colored4arm}):
\begin{enumerate}[(i)]
 \item $Q_1$ and $Q_3$ are disjoint and are at distance at least $\delta$ from each other; the
same for $Q_2$ and $Q_4$.
\item For $i\in\{1,3\}$, the sides $\partial_1 Q$ lie inside $\partial_1 A$ and the
sides $\partial_3 Q$ lie outside $\partial_2 A$; for $i\in\{2,4\}$, the sides $\partial_2 Q$
lie inside $\partial_1 A$ and the sides $\partial Q_4$ lie outside $\partial_2 A$; all these sides are at distance at least $\delta$ from the annulus $A$ and from the other $Q_j$'s.
\item The four quads are ordered cyclically around $A$ according to their indices.
\item For $i\in\{1,3\}$, we have $\omega\in\boxminus_{Q_i}$, while for $i\in\{2,4\}$, we have $\omega\in\boxminus_{Q_i}^c$.
\end{enumerate}

\subsection{Pivotal measures and Poisson point processes on the set of pivotal points}\label{pivotal}
Roughly speaking, a site is called pivotal for some macroscopic connectivity event in a discrete percolation configuration if flipping its state changes the outcome of the event.  The pivotal points for discrete percolation are easy to define, but the situation is more complicated for the quad-crossing scaling limit, and this case is treated carefully in Section 2.3 of \cite{GPS18b}.  In this section, we shall recall the definitions and some basic results of the set of pivotal points together with pivotal measures and Poisson point processes on this set for both the discrete and continuum percolation, studied in \cite{GPS13,GPS18a,GPS18b}.

Let $x$ be a point surrounded by a piecewise smooth Jordan curve $\gamma\subset D$, where ``surrounded" means that $D\setminus\gamma$ has two connected
components, with the one containing $x$ being homeomorphic to a disk.  For any
$\varepsilon>0$, fix a lattice $\varepsilon\mathbb{Z}^2$ in $D$, and let $B_{\varepsilon}(x)$ be the $\varepsilon$-box $[i,i+1)\varepsilon\times[j,j+1)\varepsilon$ in the lattice that contains $x$.  We say that $x$ is \textbf{pivotal for $\gamma$} in $\omega_{\infty}^{\lambda}$ if, for any $\varepsilon>0$ such
that $B_{2\varepsilon}(x)$ is surrounded by $\gamma$, the alternating 4-arm event occurs in the annulus with boundary pieces $\partial B_{2\varepsilon}(x)$ and $\gamma$.  Let $\mathcal{P}^{\gamma}$ denote the set of pivotal points for $\gamma$ in $D$.  Furthermore, we can identify the \textbf{color of a
pivotal point} $x\in\mathcal{P}^{\gamma}$ as \textbf{open} (blue) versus \textbf{closed} (yellow), as follows.  We let $\mathcal{P}_{\rm open}^{\gamma,\varepsilon}$ denote the set of points $x$ for which $\gamma$ surrounds $x$ without intersecting
or touching $B_{2\varepsilon}(x)$, and there exist quads $Q_{\varepsilon,i},i=1,2,3,4$, exhibiting the alternating 4-arm event from $\partial B_{2\varepsilon}(x)$ to $\gamma$ such that the quad $\tilde{U}_{\varepsilon}$, given by taking the union of
$U_{\varepsilon}:=Q_{\varepsilon,1}\cup Q_{\varepsilon,3}\cup B_{2\varepsilon}(x)$ and the bounded components of $\mathbb{C}\setminus U_{\varepsilon}$, is crossed between the boundary pieces $Q_{\varepsilon,1}(\{1\}\times[0,1])$ and $Q_{\varepsilon,3}(\{1\}\times[0,1])$ in $\omega_{\infty}^{\lambda}$; see the right-hand side of Figure \ref{colored4arm}.  Then the set of \textbf{open pivotals} for $\gamma$ is defined by
\begin{equation*}
\mathcal{P}_{\rm open}^{\gamma}(\omega_{\infty}^{\lambda}):=\{x\in D: x\in\mathcal{P}_{\rm open}^{\gamma,\varepsilon}\mbox{ for all $\varepsilon>0$ s.t. $B_{2\varepsilon}(x)$ is surrounded by $\gamma$}\}.
\end{equation*}
We set $x\in \mathcal{P}_{\rm closed}^{\gamma}(\omega_{\infty}^{\lambda})$ if the analogous dual crossing holds in the quad $\tilde{U}_{\varepsilon}$, with $U_{\varepsilon}$ replaced by $Q_{\varepsilon,2}\cup Q_{\varepsilon,4}\cup B_{2\varepsilon}(x)$, for each small enough $\varepsilon>0$.

For a discrete percolation configuration $\omega_{\eta}^{\lambda}$, we modify the above definitions and define the set of open pivotals for $\gamma$ as
\begin{equation*}
\mathcal{P}_{\rm open}^{\gamma}(\omega_{\eta}^{\lambda}):=\left\{x\in \eta\textsf{T}\cap D:
\begin{aligned}
&\mbox{the $\eta$-hexagon of $x$ is open and is surrounded by $\gamma$, and}\\
&\mbox{there exist four alternating arms from neighbors of $x$ to $\gamma$}
\end{aligned}
\right\}.
\end{equation*}
The set of closed pivotals for $\gamma$, $\mathcal{P}_{\rm closed}^{\gamma}(\omega_{\eta}^{\lambda})$, is defined similarly.  Let $\mathcal{P}^{\gamma}(\omega_{\eta}^{\lambda}):=\mathcal{P}_{\rm open}^{\gamma}(\omega_{\eta}^{\lambda})\cup\mathcal{P}_{\rm closed}^{\gamma}(\omega_{\eta}^{\lambda})$ denote the set of pivotals for $\gamma$.

Now let us introduce the normalized counting measures on the set of pivotals.  Let $D$ be a bounded domain in $\mathbb{C}$, or be equal to the torus $\mathbb{T}_M^2$. For any fixed $\varepsilon>0$, define the set of \textbf{$\varepsilon$-important points} $\mathscr{P}^{\varepsilon}(\omega_{\eta}^{\lambda})$ of any discrete percolation configuration in $D$, relative to the $(\varepsilon,3\varepsilon)$-annuli given by a fixed lattice $\varepsilon\mathbb{Z}^2$: For $x\in D$, we let $B_{\varepsilon}(x)$ be the lattice $\varepsilon$-box as before, let $\tilde{B}_{\varepsilon}(x)$ be the $3\varepsilon$-box
centered at $B_{\varepsilon}(x)$ and let $x\in\mathscr{P}^{\varepsilon}(\omega_{\eta}^{\lambda})$ if and only if $x\in\mathcal{P}^{\partial\tilde{B}_{\varepsilon}(x)}$.
The normalized counting measure $\mu^{\varepsilon}(\omega_{\eta}^{\lambda})$ on the set $\mathscr{P}^{\varepsilon}(\omega_{\eta}^{\lambda})$ is defined by
\begin{equation*}
\mu^{\varepsilon}(\omega_{\eta}^{\lambda}):=\sum_{x\in\mathscr{P}^{\varepsilon}(\omega_{\eta}^{\lambda})}\delta_x\eta^2\alpha_4(\eta,1)^{-1},
\end{equation*}
where $\delta_x$ denotes the Dirac measure concentrated at $x$.  Then we have
\begin{theorem}[Theorem 2.7 in \cite{GPS18b}]\label{t6}
For any $\varepsilon>0$, there exists a measurable map $\mu^{\varepsilon}:\mathscr{H}_D\rightarrow \mathfrak{M}_D$, into the space of finite Borel measures on $D$, such that, for $\lambda\in\mathbb{R}$, as $\eta\rightarrow 0$,
\begin{equation*}
(\omega_{\eta}^{\lambda},\mu^{\varepsilon}(\omega_{\eta}^{\lambda}))\overset{d}{\longrightarrow}(\omega_{\infty}^{\lambda},\mu^{\varepsilon}(\omega_{\infty}^{\lambda}))
\end{equation*}
in the quad-crossing topology $(\mathscr{H}_D,\mathcal{T}_D)$ in the first coordinate and in the Prokhorov distance of measures in the second one.  Furthermore,
\begin{equation*}
(\mathscr{P}_{\rm open}^{\varepsilon}(\omega_{\eta}^{\lambda}),\mathscr{P}_{\rm closed}^{\varepsilon}(\omega_{\eta}^{\lambda}))\overset{d}{\longrightarrow}
(\mathscr{P}_{\rm open}^{\varepsilon}(\omega_{\infty}^{\lambda}),\mathscr{P}_{\rm closed}^{\varepsilon}(\omega_{\infty}^{\lambda}))
\end{equation*}
in the Hausdorff metric of closed sets.
\end{theorem}

Let us use the pivotal measures in Theorem \ref{t6} to define the Poisson point processes of pivotal points, which will be used to build enhanced networks in the next subsection.
\begin{definition}[Definition 3.1 in \cite{GPS18b}]\label{def8}
Let $\bar{\lambda}=(\lambda,\lambda')\in\mathbb{R}^2$ be any pair of near-critical parameters with $\lambda<\lambda'$, and let $\varepsilon>0$ be fixed.  Let $\omega^{\lambda}$ be a near-critical configuration $\omega_{\eta}^{\lambda}$ or $\omega_{\infty}^{\lambda}$ in $\mathbb{T}_M^2$.  We will denote by $\textsf{PPP}_{\bar{\lambda}}^{\varepsilon}=\textsf{PPP}_{\bar{\lambda}}^{\varepsilon}(\mu^{\varepsilon}(\omega^{\lambda}))$ the Poisson point process
\begin{equation*}
\textsf{PPP}_{\bar{\lambda}}^{\varepsilon}:=\{(x_i,t_i),1\leq i\leq q\}\subset\mathscr{P}^{\varepsilon}(\omega^{\lambda})\times[\lambda,\lambda']
\end{equation*}
of intensity measure $\mu^{\varepsilon}(\omega^{\lambda})dx\times \mathbf{1}_{[\lambda,\lambda']}(t)dt$.  The set $\{x_1,\ldots,x_q\}$ of pivotals will usually be denoted by $X$.  The process $\textsf{PPP}_{\bar{\lambda}}^{\varepsilon}(\mu^{\varepsilon}(\omega_{\eta}^{\lambda}))$ can be constructed measurably from $\omega_{\eta}^{[\lambda,\lambda']}$, and we will always work in this natural coupling.
\end{definition}

\subsection{Enhanced networks}\label{s2}
To prove a scaling limit result, one often uses auxiliary graphs (or networks) to encode macroscopic connectivity information from percolation configurations, and these graphs are measurable with respect to the quad-crossing space (or some other topological space) and defined both for the discrete model and the scaling limit.  In this section we recall the definitions and basic properties of networks and enhanced networks from \cite{GPS18a,GPS18b}; see the second paragraph of Section \ref{strategy} for an informal description of them.  Before giving the definitions, we want to mention two related works which also use auxiliary graphs to study scaling limits.
\begin{itemize}
\item In order to prove a mesh-independent gluing theorem for a quad $Q$ with a curve $\gamma$ cutting it (see \cite[Proposition 4.1]{SS11} for the precise statement), Schramm and Smirnov \cite{SS11} constructed auxiliary graphs for a macroscopic discretization of $Q$ and the corresponding coarse-graining finite $\sigma$-field, which were used to describe the connectivity structures away from $\gamma$ and near $\gamma$.
\item In \cite{CCK19}, Camia, Conijn and Kiss proved that as the mesh $\eta\rightarrow 0$, the collection of critical percolation clusters, seen as a collection of random compact subsets of the plane, converges to a collection of continuum clusters.  In their proof, ``good" graphs (see \cite[Definition 5]{CCK19}) were used as approximations of large open clusters, with the graphs being measurable with respect to the quad-crossing percolation configurations.
\end{itemize}

\emph{In the rest of this subsection, we will assume that $M\geq 1$ and $r\in 2^{-\mathbb{N}}$.}  To define networks, we start by introducing families of dyadic squares.
\begin{definition}[A nested family of dyadic coverings]\label{coverings}
For any $b>0$ in $2^{-\mathbb{N}}={\{2^{-k}:k\in\mathbb{N}\}}$, let $S^b$ be the disjoint covering of $\mathbb{T}_M^2$ using the $b$-squares ${\{[j,j+1)b\times[k,k+1)b:(j,k)\in \mathbb{Z}^2\}}$.  Now, for any $r\in 2^{-\mathbb{N}}$ and any finite subset $X=\{x_1,\ldots,x_q\}\subset\mathbb{T}_M^2$, we associate uniquely $r$-squares $B_{x_1}^r,\ldots,B_{x_q}^r$ in the following manner: For all $1\leq i\leq q$, there is a unique square $\tilde{B}_{x_i}\in S^{r/2}$ which contains $x_i$ and we define $B_{x_i}^r$ to be the square of side length $r$ centered around the $r/2$-square $\tilde{B}_{x_i}$.  We will denote by $B^r(X)$ this family of $r$-squares.
\end{definition}

Let $X=\{x_1,\ldots,x_q\}$ be the set of pivotals of the Poisson point process $\textsf{PPP}_{\bar{\lambda}}^{\varepsilon}$ as in Definition \ref{def8}.  Although it is possible that there exist two distinct $i,j\in\{1,\ldots,q\}$ such that $x_i=x_j$, it was noted right after \cite[Definition 3.2]{GPS18b} that in the continuum case, it was essentially proved in \cite[Proposition 5.2]{GPS18a} that the points in $X(\mu^{\varepsilon}(\omega_{\infty}^{\lambda}))$ are at positive distance from each other almost surely; the proof also works for the discrete case, and one gets that with high probability, for all small enough $\eta>0$, the points in $X(\mu^{\varepsilon}(\omega_{\eta}^{\lambda}))$ are at positive distance from each other.
\begin{definition}\label{r^*}
Let $r^*(X)$ denote the smaller of $1/10$ and one-tenth of the smallest distance between any pair $x_i\neq x_j\in X$.
\end{definition}

Now we are ready to introduce the mesoscopic network.
\begin{definition}[Mesoscopic network]\label{network}
For $0<r<r^*(X)$, the \textbf{$r$-mesoscopic $\varepsilon$-network} $\textsf{N}_M^{r\text{-} \textsf{meso}}(\omega^{\lambda},\textsf{PPP}_{\bar{\lambda}}^{\varepsilon})$ associated to $\omega^{\lambda}$ in the torus $\mathbb{T}_M^2$ and $\textsf{PPP}_{\bar{\lambda}}^{\varepsilon}$ is the graph with vertex set $X$ and an edge set, with an edge connecting two distinct points $x_i,x_j\in X$ if there exists a quad $Q$ such that $\partial_1 Q$ and $\partial_3 Q$ remain strictly inside $B_{x_i}^r$ and $B_{x_j}^r$, respectively, and $Q$ remains strictly away from the squares $B_{x_k}^r, k\notin\{i,j\}$, and for which $\omega^{\lambda}\in\boxminus_Q$.
\end{definition}

Two $r$-mesoscopic networks are seen as the same if the $r$-squares for the vertices and the graph structures coincide.  For $r_1<r_2$, we can compare an $r_1$-mesoscopic network with an $r_2$-mesoscopic network by considering the unique $r_2$-squares containing the $r_1$-squares of the first network.  Now let us take $r\rightarrow 0$ to get a network for $\omega_{\infty}^{\lambda}$:
\begin{proposition}[$r$-stabilization of networks, Proposition 3.4(i) in \cite{GPS18b}]\label{p9}
There is a measurable $r_M=r_M(\omega_{\infty}^{\lambda},\textnormal{\textsf{PPP}}_{\bar{\lambda}}^{\varepsilon})\in(0,r^*(X))$, such that for all $r\in(0,r_M)$ we get the same $r$-mesoscopic $\varepsilon$-network $\textnormal{\textsf{N}}_M^{r\text{-} \textnormal{\textsf{meso}}}(\omega_{\infty}^{\lambda},\textnormal{\textsf{PPP}}_{\bar{\lambda}}^{\varepsilon})$.  This stabilized network will be called the \textbf{$\varepsilon$-network} $\textnormal{\textsf{N}}_{\infty}^{\bar{\lambda},\varepsilon}=\textnormal{\textsf{N}}_M(\omega_{\infty}^{\lambda},\textnormal{\textsf{PPP}}_{\bar{\lambda}}^{\varepsilon})$.
\end{proposition}
For discrete percolation, the \textbf{$\varepsilon$-network} $\textnormal{\textsf{N}}_{\eta}^{\bar{\lambda},\varepsilon}=\textnormal{\textsf{N}}_M(\omega_{\eta}^{\lambda},\textnormal{\textsf{PPP}}_{\bar{\lambda}}^{\varepsilon})$ is defined in the following natural way: It is the graph with vertex set $X$ and an edge set, with an edge connecting two distinct points $x_i,x_j\in X$ if there is an open path of $\omega_{\eta}^{\lambda}\setminus X$ adjacent to both $x_i$ and $x_j$.

In order to describe which of the pivotals in $X$ are connected together by an open cluster of $\omega^{\lambda}\setminus X$, let us introduce the concept of sub-routers:
\begin{definition}[Mesoscopic sub-routers]\label{subrouter}
Fix any $r<\rho$ in $2^{-\mathbb{N}}$.  As in Definition \ref{coverings}, let $B^r(\mathbb{T}_M^2)$ be the finite covering of $\mathbb{T}_M^2$ by overlapping $r$-squares.  Given a subset $Y$ of the set $X=\{x_1,\ldots,x_q\}$ of the pivotals in $\textsf{PPP}_{\bar{\lambda}}^{\varepsilon}$ with $|Y|\geq 2$, an \textbf{$(r,\rho)$-mesoscopic sub-router for $Y$} is an $r$-square $B\in B^r(\mathbb{T}_M^2)$ with the following properties:
\begin{itemize}
\item It is at distance at least $2\rho$ from each $x_i\in X$.
\item There exist open circuits (i.e., no dual arm) surrounding $B$ in the square annulus with inner face $B$ and outer radius $\rho$; the largest $r'$-square with
$r'\in 2^{-\mathbb{N}}$ that is concentric with $B$, contains it, and is surrounded by all such open circuits will be denoted by $\hat{B}$.
\item For each $x_i\in Y$, there exists a quad $Q$ with $\partial_1 Q$ strictly inside $\hat{B}$, $\partial_3 Q$ strictly inside $B_{x_i}^r$, remaining strictly away from all the squares $B_{x_k}^r$ with $x_k\in X\setminus\{x_i\}$, and for which $\omega^{\lambda}\in\boxminus_Q$.
\end{itemize}
\end{definition}
The following lemma says that the $(r,\rho)$-mesoscopic sub-routers exist.
\begin{lemma}[Lemma 3.7 in \cite{GPS18b}]\label{l13}
With probability tending to 1 as $\eta\rightarrow 0$ and then $r\rightarrow 0$, for any $0<\rho<r^*(X)$, for all $Y\subseteq X$ with $|Y|\geq 2$ whose sites are connected together by a single open cluster of $\omega_{\eta}^{\lambda}\setminus X$ (i.e., there is an open cluster of $\omega_{\eta}^{\lambda}\setminus X$ that neighbors each site in $Y$), the set of $(r,\rho)$-mesoscopic sub-routers for $Y$ is nonempty.
\end{lemma}

For any $Y\subseteq X$ with $|Y|\geq 2$, if there exists an open cluster of $\omega_{\eta}^{\lambda}\setminus X$ neighboring each site in $Y$ but not neighboring any site in $X\setminus Y$, then we call $Y$ a \textbf{cluster of pivotals} for $\omega_{\eta}^{\lambda}$.

Let $B_1,B_2$ be $(r,\rho)$ sub-routers for $Y_1,Y_2\subseteq X$, respectively.  We will call $B_1$ and $B_2$ connected if there exists a quad $Q$ with $\partial_1 Q$ contained in $\hat{B}_1$, $\partial_3 Q$ contained in $\hat{B}_2$, remaining strictly away from all the squares in $B^r(X)$, and for which $\omega^{\lambda}\in\boxminus_Q$.  The relation of being connected is an equivalence relation.  If $B_1$ and $B_2$ are connected, then it is easy to see that both of them are $(r,\rho)$ sub-routers for $Y_1\cup Y_2$.  Thus, for each equivalence class of $(r,\rho)$ sub-routers there exists a maximal subset $Y\subseteq X$ for which all elements of the equivalence class are sub-routers.  Such a maximal subset $Y$ will be called a \textbf{cluster of pivotals} for $(r,\rho)$ sub-routers.

In each equivalence class of sub-routers, we will single out one of them in a way that is typically continuous w.r.t. $\omega^{\lambda}$.  For this, we restrict ourselves to the case when every open cluster of $\omega^{\lambda}$ in $\mathbb{T}_M^2$ has diameter less than $M/3$ (for $\omega_{\infty}^{\lambda}$, the diameter is the $\limsup_{r\rightarrow 0}$ of distances between $r$-boxes that are connected in the usual sense that there is a crossed quad with its opposite sides contained in the $r$-boxes); this case occurs with high probability when $\lambda$ is very negative.  Then, for $r<M/10$, the set of $(r,\rho)$ sub-routers for any
$Y$ has an isometric embedding into $[-M/2,M/2)^2\subset\mathbb{C}$.  The leftmost sub-router of the lowermost ones in such an embedding will be the same in any of these embeddings; moreover, its location in $\mathbb{T}_M^2$ can change only a little if we move each sub-router a little.  The set of these ``leftmost of lowermost" sub-routers is called the set of \textbf{$(r,\rho)$-mesoscopic routers of $X$}, or the set of \textbf{$r$-mesoscopic routers} when
\begin{equation}\label{rho}
\rho=\rho_M:=\max\{\rho'\in 2^{-\mathbb{N}}:\rho'\leq r_M/2\wedge\varepsilon\},
\end{equation}
where $r_M$ is as in Proposition \ref{p9}.  For an $r$-mesoscopic router $B$, we call the box $\hat{B}$ an \textbf{$r$-mesoscopic expanded router}.  Note that clusters containing only one pivotal of $X$ do not have routers.

Now we are ready to define what a mesoscopic enhanced network is.
\begin{definition}[Mesoscopic enhanced network]\label{def10}
Assume that the diameter of every open cluster of $\omega^{\lambda}$ in $\mathbb{T}_M^2$ is at most $M/3$, and $r$ is smaller than $\rho_M$ defined in (\ref{rho}).  The \textbf{$r$-mesoscopic enhanced $\varepsilon$-network} $\textsf{EN}_M^{r\text{-}\textsf{meso}}(\omega^{\lambda},\textsf{PPP}_{\bar{\lambda}}^{\varepsilon})$ is the following vertex-labeled bipartite graph.  One part of the vertex set is the set $X$ of the pivotals of $\textsf{PPP}_{\bar{\lambda}}^{\varepsilon}$, and the other part is the set of $r$-mesoscopic routers of $X$.  The vertices in $X$ are colored open or closed, according to the definitions in Section \ref{pivotal}; the routers are colored open.  The edge set consists of the connections between the routers and the elements of their maximal $Y\subseteq X$.  The edges are drawn on $\mathbb{T}_M^2$ so that they are homotopic (with fixed endpoints) to the connections they represent; one can also achieve that they do not intersect each other; see \cite[Figure 10]{GPS18b} (where dual routers and dual edges were also shown, but we only consider primal routers and primal edges here).

If the assumptions on the diameters and $r$ are not satisfied, let $\textsf{EN}_M^{r\text{-}\textsf{meso}}(\omega^{\lambda},\textsf{PPP}_{\bar{\lambda}}^{\varepsilon})$ be empty.
\end{definition}

Note that if $\textsf{EN}_M^{r\text{-}\textsf{meso}}(\omega^{\lambda},\textsf{PPP}_{\bar{\lambda}}^{\varepsilon})$ is nonempty, then it induces naturally the $r$-mesoscopic $\varepsilon$-network $\textsf{N}_M^{r\text{-} \textsf{meso}}(\omega^{\lambda},\textsf{PPP}_{\bar{\lambda}}^{\varepsilon})$ of Definition \ref{network}: There is an edge between $x_i,x_j\in X$ in $\textsf{N}_M^{r\text{-} \textsf{meso}}(\omega^{\lambda},\textsf{PPP}_{\bar{\lambda}}^{\varepsilon})$ if and only if there is a router with edges to $x_i,x_j$ in $\textsf{EN}_M^{r\text{-}\textsf{meso}}(\omega^{\lambda},\textsf{PPP}_{\bar{\lambda}}^{\varepsilon})$.

For the continuum percolation $\omega_{\infty}^{\lambda}$, we will now take $r\rightarrow 0$ to get an enhanced network $\textnormal{\textsf{EN}}_{\infty}^{\bar{\lambda},\varepsilon}$, and then compare the enhanced networks for $\omega_{\eta}^{\lambda}$ and $\omega_{\infty}^{\lambda}$.  The following proposition is essentially \cite[Proposition 3.9]{GPS18b} but with a slight modification, and is stated in a form which is more suitable for us; the proof in \cite{GPS18b} works for our modification.  The statements for the degrees of pivotals need more words of explanation: If there is a vertex $x\in X$ having at least three router-neighbors in $\textnormal{\textsf{EN}}_M^{r\text{-}\textnormal{\textsf{meso}}}(\omega_{\eta}^{\lambda},\textnormal{\textsf{PPP}}_{\bar{\lambda}}^{\varepsilon})$, then we have a polychromatic 6-arm event from $B_x^r$ to distance $\rho_M/2$ in $\omega_{\eta}^{\lambda}$.  The 6-arm exponent being larger than 2 (see (\ref{e67})) implies that when $r$ is small, this event does not occur with high probability, which gives the first part of (\ref{num4}) below; the second part of (\ref{num4}) follows similarly.  Combining (\ref{num3}) and the first part of (\ref{num4}) we get the last sentence in part (i).  Moreover, notice that the first sentence in (\ref{num1}) follows from the argument right above Definition \ref{r^*}, and the second sentence in (\ref{num1}) has an analog in \cite[Proposition 3.4(ii)]{GPS18b}.

\begin{proposition}[$r$-stabilization and $\eta$-convergence of enhanced networks, Proposition 3.9 in \cite{GPS18b}]\label{p6}
(i) Let $\rho_M$ be as in (\ref{rho}).  There is a measurable scale $\tilde{r}_M=\tilde{r}_M(\omega_{\infty}^{\lambda}, \textnormal{\textsf{PPP}}_{\bar{\lambda}}^{\varepsilon})\in(0,\rho_M)$ such that for all $r\in(0,\tilde{r}_M)$ we get the same $r$-mesoscopic enhanced $\varepsilon$-network $\textnormal{\textsf{EN}}_M^{r\text{-}\textnormal{\textsf{meso}}}(\omega_{\infty}^{\lambda},\textnormal{\textsf{PPP}}_{\bar{\lambda}}^{\varepsilon})$ in the sense that the networks are the same, plus the colors in $X$ and the collections of clusters of pivotals are also the same.  (The corresponding routers do not exactly stabilize, since for a smaller $r$ new $(r,\rho_M)$ routers can appear, with each router converging to a point in $\mathbb{T}_M^2$ as $r\rightarrow 0$.)  This stabilized network will be called the \textbf{enhanced $\varepsilon$-network} $\textnormal{\textsf{EN}}_{\infty}^{\bar{\lambda},\varepsilon}=\textnormal{\textsf{EN}}_M(\omega_{\infty}^{\lambda},\textnormal{\textsf{PPP}}_{\bar{\lambda}}^{\varepsilon})$.
Recall that all elements of $X$ are distinct almost surely.  Furthermore, every vertex in $X$ has at most two router-neighbors in $\textnormal{\textsf{EN}}_{\infty}^{\bar{\lambda},\varepsilon}$ almost surely.

(ii) Write $X_{\eta}=X(\omega_{\eta}^{\lambda})$ and $X_{\infty}=X(\omega_{\infty}^{\lambda})$.  For any $\alpha>0$, there is a scale $r_{\alpha}=r_{\alpha}(M,\bar{\lambda},\varepsilon)\in 2^{-\mathbb{N}}$, such that in any coupling with $\omega_{\eta}^{\lambda}\xrightarrow{a.s.}\omega_{\infty}^{\lambda}$ in $\mathbb{T}_M^2$, for any fixed $r\in(0,r_{\alpha}]$ and all sufficiently small $\eta>0$ there is a coupling of $\textnormal{\textsf{PPP}}_{\bar{\lambda}}^{\varepsilon}(\omega_{\eta}^{\lambda})$ and $\textnormal{\textsf{PPP}}_{\bar{\lambda}}^{\varepsilon}(\omega_{\infty}^{\lambda})$, such that
with probability at least $1-\alpha/2$ the following hold:
\begin{enumerate}
  \item All elements of $X_{\eta}$ are distinct.  Moreover, $r_{\alpha}$ is less than both $\tilde{r}_M$ and $r^*(X_{\eta})$. \label{num1}
  \item The networks $\textnormal{\textsf{N}}_{\eta}^{\bar{\lambda},\varepsilon}$ and $\textnormal{\textsf{N}}_M^{r\text{-} \textnormal{\textsf{meso}}}(\omega_{\eta}^{\lambda},\textnormal{\textsf{PPP}}_{\bar{\lambda}}^{\varepsilon})$ are the same, and the collection of clusters of pivotals for $\omega_{\eta}^{\lambda}$ is the same as that for $(r,\rho_{M})$ sub-routers.  Furthermore, if all the open clusters of $\omega_{\eta}^{\lambda}$ on $\mathbb{T}_M^2$ have diameter less than $M/3$, then for every $r$-mesoscopic router $B$ with the corresponding cluster of pivotals $Y\subseteq X_{\eta}$, the open cluster of $\omega_{\eta}^{\lambda}\setminus X_{\eta}$ which contains the open circuits surrounding $\hat{B}$ in the square annulus with inner face $\hat{B}$ and outer radius $\rho_M$ neighbors every pivotal in $Y$.
      \label{num2}
  \item $\textnormal{\textsf{EN}}_M^{r\text{-}\textnormal{\textsf{meso}}}(\omega_{\eta}^{\lambda},\textnormal{\textsf{PPP}}_{\bar{\lambda}}^{\varepsilon})$ is the same as $\textnormal{\textsf{EN}}_M^{r\text{-}\textnormal{\textsf{meso}}}(\omega_{\infty}^{\lambda},\textnormal{\textsf{PPP}}_{\bar{\lambda}}^{\varepsilon})$  in the following sense: The collections $B^r(X_{\eta})$ and $B^r(X_{\infty})$ are the same; the collections of $r$-mesoscopic routers (resp. expanded routers) are the same;  for any $x_{\eta}\in X_{\eta},x_{\infty}\in X_{\infty}$ contained in the same box in $B^r(X)$ and for any $r$-mesoscopic router $B$, we have that the color of $x_{\eta}$ is the same as that of $x_{\infty}$, and there is an edge between $x_{\eta}$ and $B$ in $\textnormal{\textsf{EN}}_M^{r\text{-}\textnormal{\textsf{meso}}}(\omega_{\eta}^{\lambda},\textnormal{\textsf{PPP}}_{\bar{\lambda}}^{\varepsilon})$ if and only if there is an edge between $x_{\infty}$ and $B$ in $\textnormal{\textsf{EN}}_M^{r\text{-}\textnormal{\textsf{meso}}}(\omega_{\infty}^{\lambda},\textnormal{\textsf{PPP}}_{\bar{\lambda}}^{\varepsilon})$. \label{num3}
      \item Every pivotal in $X_{\eta}$ has at most two router-neighbors in $\textnormal{\textsf{EN}}_M^{r\text{-}\textnormal{\textsf{meso}}}(\omega_{\eta}^{\lambda},\textnormal{\textsf{PPP}}_{\bar{\lambda}}^{\varepsilon})$.  For the collection of $\lambda$-clusters on $\mathbb{T}_M^2$ that neighbor at least two $\lambda$-closed pivotals in $X_{\eta}$, each $\lambda$-closed pivotal in $X_{\eta}$ neighbors at most two clusters in this collection.
          \label{num4}
\end{enumerate}
\end{proposition}

\subsection{Correlation length}\label{correlation}
To study near-critical percolation, it is very useful to introduce the classical concept of correlation length (also called characteristic
length) for $p$ close to $p_c$.  Roughly speaking, the system ``looks like'' critical percolation on scales smaller than correlation length, while ``notable'' super/sub-critical behavior emerges above this length.  There are several natural definitions of correlation length, and the corresponding lengths turn out to be of the same order of magnitude (see \cite[Section 7]{Nol08}).

First, on the non-renormalized lattice $\textsf{T}$, we introduce the \textbf{correlation length} that is defined in terms of the box-crossing events: For each
$\varepsilon\in(0,1/2)$ and $p<p_c=1/2$, let
\begin{equation*}
L_{\varepsilon}(p):=\inf\{R\geq 1: \mathbf{P}_p[\mbox{there is an open
left-right crossing of $[0,R]^2$}]\leq \varepsilon\},
\end{equation*}
and let $L_{\varepsilon}(p):=L_{\varepsilon}(1-p)$ for $p>p_c$.  It is well known that $L_{\varepsilon}(p)\rightarrow\infty$ when $p\rightarrow p_c$.  The particular choice of $\varepsilon$ is not important in the above
definition.  Indeed, for any $\varepsilon,\varepsilon'\in(0,1/2)$, as $p\rightarrow p_c$, we have (see, e.g., \cite[Corollary 37]{Nol08})
\begin{equation*}
L_{\varepsilon}(p)\asymp L_{\varepsilon'}(p).
\end{equation*}
We will write $L(p):=L_{0.1}(p)$ and let $L(1/2):=\infty$.

For any $\eta>0$ and $\lambda\in\mathbb{R}$, define the renormalized correlation length for $\omega_{\eta}^{\lambda}$ by
\begin{equation*}
\tilde{L}_{\eta}(\lambda):=\eta L(p(\eta,\lambda)),
\end{equation*}
where $p(\eta,\lambda)$ is the percolation parameter of $\omega_{\eta}^{\lambda}$, defined in (\ref{e42}).  As explained in the proof of \cite[Theorem 10.7]{GPS18a}, Kesten's result \cite{Kes87} implies that for any $\lambda\in \mathbb{R}\setminus\{0\}$ (see \cite[(10.7)]{GPS18a}),
\begin{equation*}
0<\liminf_{\eta\rightarrow 0}\tilde{L}_{\eta}(\lambda)\leq\limsup_{\eta\rightarrow 0}\tilde{L}_{\eta}(\lambda)<\infty.
\end{equation*}

For $R>0$, denote by $\boxminus_R$ the left-right crossing event in the quad $[0,R]^2$.  As in the discrete model, we define a notion of correlation length for $\omega_{\infty}^\lambda$.  For any $\lambda<0$, define
\begin{equation*}
\tilde{L}_{\infty}(\lambda):=\inf\{R>0:\mathbf{P}_{\infty}^{\lambda}[\boxminus_R]\leq 0.1\},
\end{equation*}
and let $\tilde{L}_{\infty}(\lambda):=\tilde{L}_{\infty}(-\lambda)$ for $\lambda>0$.  It was proved in \cite[(10.4)]{GPS18a}
that there exists a constant $C>0$ such that for all $\lambda\in\mathbb{R}\setminus\{0\}$,
\begin{equation}\label{e62}
\tilde{L}_{\infty}(\lambda)=C|\lambda|^{-4/3}.
\end{equation}
Using this formula, Theorem \ref{t4} and (\ref{e59}), it is easy to get that for any $\lambda\in\mathbb{R}\setminus\{0\}$,
\begin{equation}\label{e61}
\tilde{L}_{\eta}(\lambda)\rightarrow C|\lambda|^{-4/3}\quad\mbox{as }\eta\rightarrow 0.
\end{equation}
See also \cite[(2.26)]{BKN18} where the fact that $\tilde{L}_{\eta}(\lambda)\rightarrow\tilde{L}_{\infty}(\lambda)$ as $\eta\rightarrow 0$ was used.

We will write $\tilde{L}(\lambda):=C|\lambda|^{-4/3}$ for all $\lambda\in\mathbb{R}\setminus\{0\}$ (where $C$ is as in (\ref{e62})) and let $\tilde{L}(0)=\infty$.  We also call $\tilde{L}(\lambda)$ the correlation length for $\omega_{\eta}^{\lambda}$ since for each $\lambda\in\mathbb{R}\setminus\{0\}$ we have $\tilde{L}(\lambda)\asymp\tilde{L}_{\eta}(\lambda)$ as $\eta\rightarrow 0$ by $(\ref{e61})$.

\subsection{Classical results for planar percolation}\label{classic}
For a rectangle ${\mathrm{Rect}=[x_1,x_2]\times[y_1,y_2]}$,  a
path $(v_0,v_1,\ldots,v_k)$ of $\eta\textsf{T}$ such that $v_1,\ldots,v_{k-1}\in \mathrm{Rect}$ and the line segments $\overline{v_0v_1},\overline{v_{k-1}v_k}$ intersect respectively the left and right sides of $\mathrm{Rect}$ will be called a \textbf{left-right crossing} of $\mathrm{Rect}$.  If the path is open (resp. closed), then it will be called an \textbf{open (resp. closed) left-right crossing} of $\mathrm{Rect}$.

The so-called arm events play a central role in studying
near-critical percolation.  We write 1 and 0 for ``blue'' and
``yellow'', respectively.  A \textbf{color sequence} $\sigma=(\sigma_1,\sigma_2,\dots,\sigma_k)$ is an
element of $\{0,1\}^k$ with $k\geq 1$, and its length $|\sigma|$ is $k$.  We say that a color sequence is \textbf{polychromatic} if the
sequence contains at least one 0 and one 1.  For an annulus $A(z;r,R)$, we denote by $\mathcal
{A}_{\sigma}(z;r,R)$ the event that there exist $|\sigma|$ disjoint monochromatic paths called \textbf{arms} in $A(z;r,R)$ connecting
the two boundary pieces of $A(z;r,R)$, whose colors are those prescribed by $\sigma$, when taken in counterclockwise order.  We denote by
$\mathcal {A}_{\sigma}^+(z;r,R)$ the \textbf{upper half-plane arm event} that $\mathcal
{A}_{\sigma}(z;r,R)$ occurs with the arms contained in the half-plane $z+\{u\in\mathbb{C}:\Im u\geq 0\}$.  The lower, left and right half-plane arm events are defined similarly.  Moreover, we denote by $\mathcal{A}_{\sigma}^{3\pi/2}(z;r,R)$ the \textbf{$3\pi/2$-wedge arm event} that $\mathcal
{A}_{\sigma}(z;r,R)$ occurs with all arms contained in one of the wedges $z+\mathbb{C}\setminus\textrm{Quadrant}_i$, $i\in\{1,2,3,4\}$, where $\textrm{Quadrant}_i$ denotes the $i$-th quadrant.

We assume that the reader is familiar with the RSW (Russo-Seymour-Welsh) technology, the FKG inequality, the
BK (van den Berg-Kesten) inequality and Reimer's inequality \cite{Rei00}.  Here we
collect some classical results in critical and near-critical percolation which
will be used.  See, e.g., \cite[Section 2.2]{BKN18}, \cite[Section 2.2]{Yao24} and \cite{Nol08}.  These results are usually stated for the non-renormalized lattice $\textsf{T}$ in the literature, while they are stated here for the renormalized lattice $\eta\textsf{T}$ that is more convenient for our setup.  In the following, the expression ``for any $\eta\leq r\leq\tilde{L}(\lambda)$" must be interpreted as ``for any $r\geq\eta$" when $\lambda=0$.

\begin{enumerate}[(i)]
\item \emph{RSW bounds}. For any $k,K\geq 1$, there exists a
constant $\delta=\delta(k,K)>0$, such that for each $\lambda\in \mathbb{R}$, all $\eta>0$ small enough (depending on $\lambda$) and $\eta\leq r\leq K\tilde{L}(\lambda)$,
\begin{align*}
&\mathbf{P}_{\eta}^{\lambda}[\mbox{there is an open left-right crossing of $[0,kr]\times[0,r]$}]\geq\delta\quad\mbox{and}\\
&\mathbf{P}_{\eta}^{\lambda}[\mbox{there is a closed left-right crossing of $[0,kr]\times[0,r]$}]\geq\delta.
\end{align*}
\item \emph{Quasi-multiplicativity}.
For any color sequence $\sigma$ and $K\geq 1$, there is a constant $C=C(|\sigma|,K)>1$, such that for all $\lambda\in \mathbb{R}$, all $\eta>0$ small enough (depending on $\lambda$) and $\eta\leq r_1<r_2<r_3\leq K\tilde{L}(\lambda)$,
\begin{equation}\label{e71}
C^{-1}\mathbf{P}_{\eta}^{\lambda}[\mathcal
{A}_{\sigma}(r_1,r_3)]\leq\mathbf{P}_{\eta}^{\lambda}[\mathcal {A}_{\sigma}(r_1,r_2)]\mathbf{P}_{\eta}^{\lambda}[\mathcal
{A}_{\sigma}(r_2,r_3)]\leq C\mathbf{P}_{\eta}^{\lambda}[\mathcal
{A}_{\sigma}(r_1,r_3)].
\end{equation}
\item \emph{Stability for arm events near criticality}.
For any color sequence $\sigma$ and $K\geq 1$, there is a constant $C=C(|\sigma|,K)>1$, such that for all $\lambda\in \mathbb{R}$, all $\eta>0$ small enough (depending on $\lambda$) and $\eta\leq r<R\leq K\tilde{L}(\lambda)$,
\begin{equation}\label{e68}
C^{-1}\mathbf{P}_{\eta}^0[\mathcal
{A}_{\sigma}(r,R)]\leq\mathbf{P}_{\eta}^{\lambda}[\mathcal {A}_{\sigma}(r,R)]\leq
C\mathbf{P}_{\eta}^0[\mathcal {A}_{\sigma}(r,R)].
\end{equation}
\item \emph{Upper bound on the 6-arm exponent}.
There is a universal constant $\varepsilon_6>0$ such that for any $K\geq 1$ and polychromatic color sequence $\sigma$ with $|\sigma|=6$, there
exists a constant $C=C(K)>0$, such that for all $\lambda\in\mathbb{R}$, all $\eta>0$ small enough (depending on $\lambda$) and $\eta\leq r<R\leq K\tilde{L}(\lambda)$,
\begin{equation}\label{e67}
\mathbf{P}_{\eta}^{\lambda}[\mathcal {A}_{\sigma}(r,R)]\leq
C\left(\frac{r}{R}\right)^{2+\varepsilon_6}.
\end{equation}
\item \emph{Upper bounds on the half-plane and $3\pi/2$-wedge 3-arm exponents}.
Let $\sigma$ be any color sequence with $|\sigma|=3$.  For any $K\geq 1$, there exists a constant $C=C(K)>0$,
such that for all $\lambda\in\mathbb{R}$, all $\eta>0$ small enough (depending on $\lambda$) and $\eta\leq r<R\leq K\tilde{L}(\lambda)$,
\begin{equation}\label{e70}
\mathbf{P}_{\eta}^{\lambda}\left[\mathcal {A}_{\sigma}^+(r,R)\right]\leq
C\left(\frac{r}{R}\right)^2;
\end{equation}
similar inequalities also hold for the lower, left and right half-plane 3-arm events.  Moreover, there is a universal constant $\varepsilon_3>0$, such that for any $K\geq 1$, there exists a constant $C=C(K)>0$, such that for all $\lambda\in\mathbb{R}$, all $\eta>0$ small enough (depending on $\lambda$) and $\eta\leq r<R\leq K\tilde{L}(\lambda)$,
\begin{equation}\label{e72}
\mathbf{P}_{\eta}^{\lambda}\left[\mathcal {A}_{\sigma}^{3\pi/2}(r,R)\right]\leq
C\left(\frac{r}{R}\right)^{1+\varepsilon_3}.
\end{equation}
See, e.g., \cite[(13)]{Yao24}.
\item \emph{Upper bound on diameters of open clusters}.  Fix any $\alpha>0$ and any $M\geq 1$.  Consider percolation $\omega_{\eta}^{\lambda}$ on $\mathbb{T}_M^2$ with $\lambda<-1$.  For any $s>0$, if $\lambda<-1$ is small enough, then for all $\eta>0$ small enough, we have (see \cite[Lemma 4.3 (i)]{GPS18b})
\begin{equation}\label{e76}
\mathbf{P}_{\eta}^{\lambda}\left[\mbox{all open clusters of $\omega_{\eta}^{\lambda}$ on $\mathbb{T}_M^2$ have diameter less than $s$}\right]\geq 1-\alpha.
\end{equation}
\item \emph{Exponential decay with respect to $\tilde{L}(\lambda)$}.  There exist constants $C_1,C_2>0$, such that for all
$\lambda<0$, all $\eta>0$ small enough (depending on $\lambda$) and $R\geq\tilde{L}(\lambda)$ (see, e.g., \cite[(vi) in Section 2]{Yao19}),
\begin{align}
&\mathbf{P}_{\eta}^{\lambda}[\mbox{$\exists$ a closed circuit surrounding 0 in
$A(R,2R)$}]\geq 1-\exp\left(-C_1\frac{R}{\tilde{L}(\lambda)}\right),\label{e43}\\
&\mathbf{P}_{\eta}^{\lambda}[\mbox{the infinite closed cluster intersects $[-R,R]^2$}]\geq 1-\exp\left(-C_2\frac{R}{\tilde{L}(\lambda)}\right).\label{e44}
\end{align}
The above two inequalities also hold in the $\lambda>0$ case when ``closed" is replaced by ``open".
\end{enumerate}

\subsection{A property of large clusters}
Recall that $X$ denotes the set of pivotals of $\textsf{PPP}_{\bar{\lambda}}^{\varepsilon}$ as in Definition \ref{def8}.  In this section we prove the following lemma, which says that with high probability every large open cluster of $\omega_{\eta}^{\lambda}$ on the torus $\mathbb{T}_M^2$ contains at least two distinct sites in $X$.  This lemma will be used in the proof of Lemma \ref{l19} to show that every large open cluster has corresponding routers with high probability.

\begin{lemma}\label{l8}
Fix any $\lambda,\lambda'\in\mathbb{R}$ with $\lambda<\lambda'-1$, and fix any $M\geq 1,\rho\in(0,M/3]$ and any $\alpha\in(0,1)$.  If $\varepsilon>0$ is small enough, then for all $\eta>0$ small enough, we have
\begin{equation*}
\mathbf{P}\left[
\begin{aligned}
&\mbox{each open cluster of $\omega_{\eta}^{\lambda}$ on $\mathbb{T}_M^2$ with diameter}\\
&\mbox{at least $\rho$ contains at least two distinct sites in $X$}
\end{aligned}
\right]\geq 1-\alpha.
\end{equation*}
\end{lemma}

\begin{proof}
The proof is in large part similar to that of \cite[Proposition 3.6]{GPS18b} for volume estimates of clusters with large diameter, so we will omit some details and highlight the places that need modifications.

Take the lattice $(\rho/3)\mathbb{Z}^2$, and for each $\rho/3$-square consider the square annulus of exterior side-length $\rho$ centered around the $\rho/3$-square.  Observe that any open cluster of diameter at least $\rho$ produces a radial crossing of such a $(\rho/3,\rho)$-annulus.  Whether a given $(\rho/3,\rho)$-annulus $A_{\rho}$ is radially crossed can be decided using the radial exploration process started at any point along $\partial_1A_{\rho}$, stopped
when reaching $\partial_2A_{\rho}$ (see \cite[Figure 2.6]{GPS13}).  If $A_{\rho}$ is radially crossed, then there are two cases: either (a) there is also an open circuit in $A_{\rho}$ surrounding $\partial_1 A_{\rho}$, or (b) there is also at
least one radial dual (i.e., closed) crossing.

(a) Condition on having an open circuit in $A_{\rho}$ surrounding $\partial_1 A_{\rho}$ (this is more general than the
first of the two above cases).  Condition on the innermost such circuit, $\Gamma$.  Denote by $A_{\Gamma}$ the annulus between $\Gamma$ and $\partial_2A_{\rho}$.  Note that the percolation configuration in $A_{\Gamma}$ is undisturbed.  Let $\varepsilon\in(0,(\tilde{L}(\lambda)\wedge\rho)/100)$.
By the half-plane 3-arm exponent being 2 (see (\ref{e70})), the event that the distance between $\Gamma$ and $\partial_2A_{\rho}$ is at least $12\varepsilon$ happens with probability at least $1-O(\varepsilon/\rho)$.  Assume that this event occurs; take the set of points of $A_{\Gamma}$ whose distance from $\Gamma$ is at most $10\varepsilon$, denoted by $\tilde{A}_{\Gamma,\varepsilon}$.  Clearly $\tilde{A}_{\Gamma,\varepsilon}$ contains a collection of $K=\lfloor C_0\rho/\varepsilon\rfloor$ (where $C_0$ is an absolute constant) disjoint balls of diameter $10\varepsilon$, denoted by $\tilde{A}_i,i=1,\ldots,K$, such that all their pairwise distances are at least $2\varepsilon$; see the first picture in \cite[Figure 7]{GPS18b}.

Let $\tilde{B}_i$ be the ball of diameter $11\varepsilon$ centered at $\tilde{A}_i$, and let $B_i$ be the $\varepsilon$-square of $\varepsilon\mathbb{Z}^2$ that contains the center of $\tilde{A}_i$.  In the proof of \cite[Proposition 3.6]{GPS18b}, in order to estimate the volume of a cluster with large diameter, the authors considered the number of sites in $\tilde{A}_i$ that are connected to $\Gamma$ within ˜$\tilde{B}_i$; here, we let $Y_i(\varepsilon)$ be the number of sites in ˜$B_i\cap X$ that are connected in $\omega_{\eta}^{\lambda}$ to $\Gamma$ within ˜$\tilde{B}_i$ (where the counted sites can be $\lambda$-closed).  Recall the set of $\varepsilon$-important points $\mathscr{P}^{\varepsilon}(\omega_{\eta}^{\lambda})$ and the normalized counting measure $\mu^{\varepsilon}(\omega_{\eta}^{\lambda})$ on $\mathscr{P}^{\varepsilon}(\omega_{\eta}^{\lambda})$ defined in Section \ref{pivotal}.  By the classical
separation of arms phenomenon, RSW and FKG gluing techniques together with the quasi-multiplicativity property (see Section \ref{classic}), it is standard to prove that for any site $x$ in $B_i$ and all $\eta>0$ small enough,
\begin{equation}\label{e64}
\mathbf{P}[x\in\mathscr{P}^{\varepsilon}(\omega_{\eta}^{\lambda})\mbox{ and $x$ is connected in $\omega_{\eta}^{\lambda}$ to $\Gamma$ within ˜$\tilde{B}_i$}  \mid\Gamma]\asymp\mathbf{P}[x\in\mathscr{P}^{\varepsilon}(\omega_{\eta}^{\lambda})],
\end{equation}
and furthermore, for any pair of sites $x,y$ in $B_i$ and all $\eta>0$ small enough,
\begin{align}
&\mathbf{P}[x,y\in\mathscr{P}^{\varepsilon}(\omega_{\eta}^{\lambda}),\mbox{ and both $x$ and $y$ are connected in $\omega_{\eta}^{\lambda}$ to $\Gamma$ within ˜$\tilde{B}_i$}\mid\Gamma]\nonumber\\
&\asymp\mathbf{P}[x,y\in\mathscr{P}^{\varepsilon}(\omega_{\eta}^{\lambda})],\label{e65}
\end{align}
where the constant factors implied by the notation $\asymp$ in the two above equations are independent of everything; all the following $\asymp$ have the same meaning.  By a slight generalization of \cite[Lemma 4.5]{GPS13} to general $\lambda$ (the same proof works since $\varepsilon$ is smaller than the correlation length $\tilde{L}(\lambda)$ and the stability of multi-arm probabilities (\ref{e68}) holds),
we have
\begin{equation}\label{e63}
\mathbf{E}[\mu^{\varepsilon}(\omega_{\eta}^{\lambda})(B_i)]\asymp\varepsilon^{3/4}\quad\mbox{and}\quad \mathbf{E}[\mu^{\varepsilon}(\omega_{\eta}^{\lambda})(B_i)^2]\asymp\varepsilon^{3/2}.
\end{equation}

The normalized counting measure $\mu_{\Gamma}^{\varepsilon}(\omega_{\eta}^{\lambda})$ on the set
\begin{equation*}
\mathscr{P}_{\Gamma,i}^{\varepsilon}(\omega_{\eta}^{\lambda}):=\{x\in\mathscr{P}^{\varepsilon}(\omega_{\eta}^{\lambda})\cap B_i:\mbox{ $x$ is connected in $\omega_{\eta}^{\lambda}$ to $\Gamma$ within ˜$\tilde{B}_i$}\}
\end{equation*}
is defined by
\begin{equation*}
\mu_{\Gamma}^{\varepsilon}=\mu_{\Gamma}^{\varepsilon}(\omega_{\eta}^{\lambda}):=\sum_{x\in\mathscr{P}_{\Gamma,i}^{\varepsilon}(\omega_{\eta}^{\lambda})}\delta_x\eta^2\alpha_4(\eta,1)^{-1}.
\end{equation*}
Then, combining (\ref{e64}), (\ref{e65}) and (\ref{e63}) gives
\begin{equation}\label{e11}
\mathbf{E}[\mu_{\Gamma}^{\varepsilon}(B_i)\mid\Gamma]\asymp\varepsilon^{3/4}\quad\mbox{and}\quad \mathbf{E}[\mu_{\Gamma}^{\varepsilon}(B_i)^2\mid\Gamma]\asymp\varepsilon^{3/4}.
\end{equation}
Thus, by the Paley--Zygmund second moment inequality, there is a uniform constant $C_1>0$ such that
\begin{equation}\label{e21}
\mathbf{P}[\mu_{\Gamma}^{\varepsilon}(B_i)\geq C_1\varepsilon^{3/4}\mid\Gamma]\geq C_1.
\end{equation}
By the definition of $\textnormal{\textsf{PPP}}_{\bar{\lambda}}^{\varepsilon}$, there is a uniform constant $C_2>0$ such that
\begin{align}
\mathbf{P}[Y_i(\varepsilon)\geq 1\mid\Gamma]&=\mathbf{E}[\mathbf{E}[\mathbf{1}_{Y_i(\varepsilon)\geq 1}\mid\mu_{\Gamma}^{\varepsilon}(B_i)]\mid\Gamma]\nonumber\\
&\geq\mathbf{E}[\mathbf{1}_{\mu_{\Gamma}^{\varepsilon}(B_i)\geq C_1\varepsilon^{3/4}}\mathbf{E}[\mathbf{1}_{Y_i(\varepsilon)\geq 1}\mid\mu_{\Gamma}^{\varepsilon}(B_i),\mu_{\Gamma}^{\varepsilon}(B_i)\geq C_1\varepsilon^{3/4}]\mid\Gamma]\nonumber\\
&\geq\mathbf{E}[\mathbf{1}_{\mu_{\Gamma}^{\varepsilon}(B_i)\geq C_1\varepsilon^{3/4}}(1-\exp(-(\lambda'-\lambda)C_1\varepsilon^{3/4}))\mid\Gamma]\nonumber\\
&\geq C_1(1-\exp(-(\lambda'-\lambda)C_1\varepsilon^{3/4}))\quad\mbox{by (\ref{e21})}\nonumber\\
&\geq C_2\varepsilon^{3/4}\quad\mbox{by $\lambda'-\lambda>1$}.\label{e74}
\end{align}
Let $Y_i'(\varepsilon)$ be the number of $\lambda$-open sites in ˜$B_i\cap X$ that are connected in $\omega_{\eta}^{\lambda}$ to $\Gamma$ within ˜$\tilde{B}_i$.  It is well known that for any fixed $\lambda\in\mathbb{R}$, the percolation parameter $p(\eta,\lambda)$ of $\omega_{\eta}^{\lambda}$ (defined in (\ref{e42})) tends to $1/2$
as $\eta\rightarrow 0$.  Then the inequality (\ref{e74}) implies that there is a uniform constant $C_3>0$ such that for each $Y_i'(\varepsilon)$ and all $\eta>0$ small enough,
\begin{equation*}
\mathbf{P}[Y_i'(\varepsilon)\geq 1\mid\Gamma]\geq C_3\varepsilon^{3/4}.
\end{equation*}
Using this inequality, the independence of the variables $Y_1'(\varepsilon),\ldots,Y_K'(\varepsilon)$ conditionally on $\Gamma$ and the fact that $K=\lfloor
C_0\rho/\varepsilon\rfloor$, we get that there exist constants $C_4,C_5>0$ such that for all $\eta>0$ small enough,
\begin{align}
&\mathbf{P}[\mbox{there is at most one site of $X$ lying in the cluster containing $\Gamma$ in $\omega_{\eta}^{\lambda}$}\mid\Gamma]\nonumber\\
&\leq \mathbf{P}[Y_1'(\varepsilon)+\cdots+Y_K'(\varepsilon)\leq 1\mid\Gamma]\nonumber\\
&\leq (K+1)(1-C_3\varepsilon^{3/4})^{K-1}\leq C_4\exp(-C_5\rho\varepsilon^{-1/4}).\label{e16}
\end{align}
The above arguments imply that with probability at least ${1-O((M/\rho)^2(\varepsilon/\rho+\exp(-C_5\rho\varepsilon^{-1/4})))}$, tending to 1 as $\varepsilon\rightarrow 0$, uniformly in all small enough $\eta$, in all the at most $O((M/\rho)^2)$ annuli where case (a) occurs, there are at least two distinct sites of $X$ lying in the cluster of $\Gamma$ in $\omega_{\eta}^{\lambda}$.

(b) Condition on the second case.  For $m\in\mathbb{N}$, denote by $\mathcal{A}(m)$ the event that there are at least $m$ radially crossing open clusters and $m$ radially crossing closed clusters of $\omega_{\eta}^{\lambda}$ in $A_{\rho}$.  Thus, the second case corresponds to $\mathcal{A}(1)$.  By Reimer's inequality and RSW, there is a universal constant $C_6>0$, such that for all $m\in\mathbb{N}$ and all small $\eta>0$,
\begin{equation}\label{e20}
\mathbf{P}_{\eta}^{\lambda}[\mathcal{A}(m+1)\mid\mathcal{A}(1)]\leq(\mathbf{P}_{\eta}^{\lambda}[\mathcal{A}(1)])^m\leq\exp(-C_6m).
\end{equation}
Let $\mathcal{E}(\varepsilon)$ denote the event that for each radially crossing open cluster in $A_{\rho}$, the $12\varepsilon$-neighbourhood to its right boundary is an unexplored space.  As discussed in part (b) of the proof of \cite[Proposition 3.6]{GPS18b}, there exist universal constants $C_7,\beta>0$, such that for all small $\eta>0$,
\begin{equation}\label{e22}
\mathbf{P}_{\eta}^{\lambda}[\mathcal{E}(\varepsilon)\mid\mathcal{A}(1)]\geq 1-C_7(\varepsilon/\rho)^{\beta}.
\end{equation}
Using an analogous argument as in case (a) above and letting $C_4,C_5>0$ be as in (\ref{e16}), we get that for any $m\in\mathbb{N}$ and all small $\eta>0$,
\begin{align}
&\mathbf{P}\left[
\begin{aligned}
&\mbox{there exists a radially crossing open cluster of}\\
&\mbox{$\omega_{\eta}^{\lambda}$ in $A_{\rho}$ that contains at most one site of $X$}
\end{aligned}
\mid(\mathcal{A}(1)\cap\mathcal{E}(\varepsilon))\setminus\mathcal{A}(m+1)\right]\nonumber\\
&\leq C_4m\exp(-C_5\rho\varepsilon^{-1/4}).\label{e23}
\end{align}

Choose $m=\lfloor\rho/\varepsilon\rfloor$.  Combining (\ref{e20}), (\ref{e22}) and (\ref{e23}), we obtain that with probability at least ${1-O((M/\rho)^2(\exp(-C_6\rho/\varepsilon)+(\varepsilon/\rho)^{\beta}+(\rho/\varepsilon)\exp(-C_5\rho\varepsilon^{-1/4})))}$, tending to 1 as $\varepsilon\rightarrow 0$, uniformly in all small enough $\eta>0$,
in all the at most $O((M/\rho)^2)$ annuli where case (b) occurs, each radially crossing open cluster of $\omega_{\eta}^{\lambda}$ in such an annulus contains at least two distinct sites of $X$.

Combining the results in the above two cases (a) and (b) yields the lemma.
\end{proof}

\section{Cluster Chains}\label{clusterchain}

Chains of clusters are very useful in the study of the geometry of critical and near-critical models from statistical physics. In this section, we collect some properties of cluster chains that will be used in our proofs.  For related works using chains of clusters or cluster-boundary loops, see, for example, critical percolation on $\textsf{T}$ \cite{BN11,CN06,Yao14}, near-critical percolation on $\textsf{T}$ \cite{GPS18b,Yao24}, the critical random-cluster model on $\mathbb{Z}^2$ with cluster-weight $q\in[1,4)$ \cite{DCMT21}, and the critical FK-Ising model ($q=2$) \cite{CJN20}.

\subsection{Definitions}
For a region $D\subset\mathbb{C}$, let $D^{\eta}$ denote the union of all the hexagons of $\eta\textsf{H}$ contained in $D$, and let $\omega_{\eta}^{\lambda}\cap D$ be the restriction of $\omega_{\eta}^{\lambda}$ to $D^{\eta}$.

\begin{definition}\label{chaindef}
A \textbf{chain} of $\omega_{\eta}^{\lambda}$ is a finite sequence $(\mathcal{C}_1,y_1,\mathcal{C}_2,y_2,\ldots,y_{k-1},\mathcal{C}_k)$ such that $\mathcal{C}_1,\mathcal{C}_2,\ldots,\mathcal{C}_k$ are distinct open clusters of $\omega_{\eta}^{\lambda}$, $y_1,y_2,\ldots,y_{k-1}$ are distinct closed sites, and for each $i\in\{1,\ldots,k-1\}$, the site $y_i$ neighbors both $\mathcal{C}_i$ and $\mathcal{C}_{i+1}$.  In particular, a single open cluster is a chain.  For a region $D$, a chain of $\omega_{\eta}^{\lambda}\cap D$ is defined similarly, with the open clusters being the $D^{\eta}$-restricted version and all the closed sites of the chain lying in $D^{\eta}$.

For $k\geq 2$, a \textbf{necklace} of $\omega_{\eta}^{\lambda}$ is a finite sequence $(\mathcal{C}_1,y_1,\mathcal{C}_2,y_2,\ldots,y_{k-1},\mathcal{C}_k,y_k)$ such that $(\mathcal{C}_1,y_1,\mathcal{C}_2,y_2,\ldots,y_{k-1},\mathcal{C}_k)$ is a chain, $y_1,\ldots,y_k$ are distinct closed sites, and $y_k$ neighbors both $\mathcal{C}_k$ and $\mathcal{C}_1$.  Moreover, a single open cluster $\mathcal{C}$ is also called a necklace, and if $y$ is a closed site neighboring two distinct sites of $\mathcal{C}$, then $(\mathcal{C},y)$ is called a necklace again.  A necklace of $\omega_{\eta}^{\lambda}\cap D$ for a region $D$ is defined similarly.
\end{definition}

Let $\mathrm{Rect}$ denote a rectangle with sides parallel to the axes.  For $\omega_{\eta}^{\lambda}$, we say that there is a \textbf{left-right crossing chain} in $\mathrm{Rect}$ if there is a chain $(\mathcal{C}_1,\ldots,\mathcal{C}_k)$ of $\omega_{\eta}\cap\mathrm{Rect}$ such that the distance between $\mathcal{C}_1$ (viewed as a union of $\eta$-hexagons) and the left side is smaller than $\eta$, and the distance between $\mathcal{C}_k$ and the right side is also smaller than $\eta$.  The \textbf{top-bottom crossing chain} in $\mathrm{Rect}$ is defined similarly, but with the corresponding distance smaller than $\sqrt{3}\eta/2$ rather than $\eta$.  For $\omega_{\eta}^{\lambda},\mathrm{Rect}$ and $r>0$, define the events:
\begin{align*}
&\mathcal{H}(\mathrm{Rect},r):=\left\{
\begin{aligned}
&\mbox{there is a left-right crossing chain in $\mathrm{Rect}$, such that}\\
&\mbox{all the clusters in this chain have diameter at least $r$}
\end{aligned}
\right\},\\
&\mathcal{V}(\mathrm{Rect},r):=\left\{
\begin{aligned}
&\mbox{there is a top-bottom crossing chain in $\mathrm{Rect}$, such that}\\
&\mbox{all the clusters in this chain have diameter at least $r$}
\end{aligned}
\right\}.
\end{align*}

\subsection{Basic properties of chains}\label{chain}

\begin{lemma}[Crossing rectangles by chains of large clusters]\label{l7}
For each $\alpha>0$, there exists $\delta\in(0,1)$, such that for any $\lambda\leq 0$, any $0<R<\tilde{L}(\lambda)$, and any rectangle $\mathrm{Rect}_R$ that is a translation of $[0,3R]\times[0,R]$ or $[0,R]\times[0,3R]$, for all small enough mesh $0<\eta<\eta_0(\alpha,\lambda,R)$, we have
\begin{equation*}
\mathbf{P}_{\eta}^{\lambda}[\mathcal{H}(\mathrm{Rect}_R,\delta R)]\geq 1-\alpha\quad\mbox{and}\quad\mathbf{P}_{\eta}^{\lambda}[\mathcal{V}(\mathrm{Rect}_R,\delta R)]\geq 1-\alpha.
\end{equation*}
\end{lemma}
\begin{proof}
\cite[Theorem 7.5]{DCMT21} is a similar result for critical random-cluster model on $\mathbb{Z}^2$ with cluster-weight $q\in [1,4)$ (including critical Bernoulli bond percolation).  Its proof also works for our near-critical site percolation, since $R$ is smaller than the correlation length $\tilde{L}(\lambda)$, which implies that by choosing $\eta>0$ small enough we can use the RSW bounds, the stability for arm events near criticality and the quasi-multiplicativity (see Section \ref{classic}).
\end{proof}
\begin{remark}
Lemma \ref{l7} also follows from the proof of \cite[Lemma 4.1]{GPS18b} with straightforward modifications.  This approach is applicable only to Bernoulli percolation since it relies on the coupling of near-critical percolation and the Poisson point process $\textsf{PPP}_{\bar{\lambda}}^{\varepsilon}$, while the proof of \cite[Theorem 7.5]{DCMT21} applies to more general models.
\end{remark}

Let $A$ be an annulus.  If $\mathfrak{N}$ is a necklace of $\omega_{\eta}^{\lambda}\cap A$, then we say it surrounds the inner face of $A$ if there is a circuit surrounding the inner face of $A$ so that each site of the circuit belongs to an open cluster or a closed site of $\mathfrak{N}$.  For an annulus $A$ and $r>0$, define the event
\begin{equation*}
\mathcal{A}(\mathrm{A},r):=\left\{
\begin{aligned}
&\mbox{there is a necklace of $\omega_{\eta}^{\lambda}\cap A$ such that it surrounds the inner}\\
&\mbox{face of $A$, and all its clusters have diameter at least $r$}
\end{aligned}
\right\}.
\end{equation*}

Consider the covering of $\mathbb{T}_M^2$ by the $R$-squares given by the lattice $R\mathbb{Z}^2$, and let $\mathscr{A}_{M,R}$ be the collection of
square $(R,3R)$-annuli of exterior side-length $3R$ centered around these $R$-squares.  The following lemma is somewhat analogous to \cite[Lemma 4.2]{GPS18b}.
\begin{lemma}[Global necklace lemma]\label{l9}
For any $\alpha>0,M\geq 1$ and $0<R<M/3$, there exists $r\in(0,R)$, such that for any $\lambda\leq 0$, for all small enough $0<\eta<\eta_0(\alpha,\lambda,R,M)$, we have
\begin{equation}\label{e25}
\mathbf{P}_{\eta}^{\lambda}[\mbox{$\mathcal{A}(\mathrm{A},r)$ occurs for all $A\in\mathscr{A}_{M,R}$}]\geq 1-\alpha.
\end{equation}
\end{lemma}
\begin{proof}
For any $A\in\mathscr{A}_{M,R}$, let $A_L,A_B,A_R,A_T$ denote the left, bottom, right, and top
rectangles in $A$ which are composed of three squares of $R\mathbb{Z}^2$, respectively.  Observe that when the events $\mathcal{V}(A_L,r)$, $\mathcal{V}(A_R,r)$, $\mathcal{H}(A_B,r)$ and $\mathcal{H}(A_T,r)$ occur simultaneously, then so does $\mathcal{A}(\mathrm{A},r)$.  Note that there are $(2M/R)^2$ annuli in $\mathscr{A}_{M,R}$.  Thus, by Lemma \ref{l7} and the union bound, for any $\alpha>0,M\geq 1$ and $0<R<\tilde{L}(\lambda)\wedge M/3$, there exists $\delta\in(0,1)$, such that for any $\lambda\leq 0$, for all small enough $0<\eta<\eta_1(\alpha,\lambda,R,M)$, one has
\begin{equation}\label{e24}
\mathbf{P}_{\eta}^{\lambda}[\mbox{$\mathcal{A}(\mathrm{A},\delta R)$ occurs for all $A\in\mathscr{A}_{M,R}$}]\geq 1-\alpha.
\end{equation}
Then, for the case of $0<R<\tilde{L}(\lambda)\wedge M/3$, (\ref{e25}) follows from (\ref{e24}) by letting $r=\delta R$.

Now assume that $\tilde{L}(\lambda)\leq R<M/3$; it remains to prove this case.  Let $K$ be the smallest whole number such that $R/3^K<\tilde{L}(\lambda)$.  Observe that if $\mathcal{A}(\mathrm{A},\delta R/3^K)$ occurs for all $A\in\mathscr{A}_{M,R/3^K}$, then so does $\mathcal{A}(\mathrm{A},\delta R/3^K)$ for all $A\in\mathscr{A}_{M,R}$, since for each $A\in\mathscr{A}_{M,R}$ one can use a sequence of necklaces in the annuli of $\mathscr{A}_{M,R/3^K}$ to construct a necklace in $A$ that surrounds the inner face of $A$.  Therefore, for the case of $\tilde{L}(\lambda)\leq R<M/3$, (\ref{e25}) follows from (\ref{e24}) by letting $r=\delta R/3^K$ and taking $\eta_0(\alpha,\lambda,R,M)=\eta_1(\alpha,\lambda,R/3^K,M)$.
\end{proof}

\subsection{Chains with all closed sites in the Poisson point set $X$}\label{chain2}
Our aim in this section is to prove Lemma $\ref{l3}$ below.  Recall that $X$ is the set of pivotals of the Poisson point process $\textsf{PPP}_{\bar{\lambda}}^{\varepsilon}$ as in Definition \ref{def8}.  Lemma $\ref{l3}$ enables us to study the geometry of large clusters of $\omega_{\eta}^{\lambda'}$ by using chains of $\omega_{\eta}^{\lambda}$ which have all their closed sites in $X$ under the near-critical coupling.

\begin{lemma}\label{l3}
Fix any $\lambda'\in\mathbb{R},M\geq 1,\alpha,s\in(0,1/2)$ and $\delta\in(0,s/2)$.  Consider the near-critical coupling $(\omega_{\eta}^t)_{t\in \mathbb{R}}$ on the torus $\mathbb{T}_M^2$.  If $\lambda<\lambda'$ is very negative (depending on $\lambda',M,\alpha,s$ but independently of $\delta$) and $\varepsilon>0$ is small enough, then for all small enough $\eta>0$, with probability at least $1-\alpha$ the following holds:

For any square $B\subset\Lambda_{M-s}$, if there are two $\lambda$-clusters $\mathcal{C},\mathcal{C}'$ of diameter at least $\delta$ completely contained in $B$, and there is an open path of $\omega_{\eta}^{\lambda'}$ connecting $\mathcal{C}$ and $\mathcal{C}'$ in $B$, then there exists a chain $(\mathcal{C}_1=\mathcal{C},y_1,\mathcal{C}_2,\ldots,y_{k-1},\mathcal{C}_k=\mathcal{C}')$ of $\omega_{\eta}^{\lambda}$ contained in $B+\Lambda_s$, such that all $y_1,\ldots,y_{k-1}$ belong to $X=X(\textnormal{\textsf{PPP}}_{\bar{\lambda}}^{\varepsilon})$, with $\bar{\lambda}=(\lambda,\lambda')$.
\end{lemma}
To prove this lemma, we need to introduce a few definitions and a preliminary lemma from \cite[Section 8]{GPS18a}.
\begin{definition}[Importance of a point]
Given a percolation configuration $\omega_{\eta}^{\lambda}$ and a site $z$, let $Z(z)=Z_{\omega_{\eta}^{\lambda}}(z)$ denote the maximal radius $r$ such that the alternating 4-arm event holds from the hexagon of $z$ to distance $r$ away.  The quantity $Z(z)$ is called the \textbf{importance} of $z$ in $\omega_{\eta}^{\lambda}$.
\end{definition}

Fix $\lambda,\lambda'\in\mathbb{R}$ with $\lambda<\lambda'$.  Recall that $\bar{\lambda}=(\lambda,\lambda')$.  We will denote by $S=S_{\eta,\bar{\lambda}}$ the random set of sites on $\eta\textsf{T}$ which are updated along the near-critical dynamics ${[\lambda,\lambda']\ni t\mapsto \omega_{\eta}^{t}}$.  Let $\Omega(\omega_{\eta}^{\lambda},S)$ denote the set of percolation configurations $\tilde{\omega}_{\eta}^{\lambda}$ such that
$\tilde{\omega}_{\eta}^{\lambda}(z)=\omega_{\eta}^{\lambda}(z)$ for all $z\notin S$.  Set
\begin{equation*}
Z^S(z):=\sup_{\tilde{\omega}_{\eta}^{\lambda}\in\Omega(\omega_{\eta}^{\lambda},S)}Z_{\tilde{\omega}_{\eta}^{\lambda}}(z).
\end{equation*}

The following lemma is a slight generalization of \cite[Lemma 8.5]{GPS18a}; it has been implicitly used in the proof of \cite[Lemma 4.1]{GPS18b}.
\begin{lemma}\label{l4}
Fix any $-\infty<\lambda<\lambda'<\infty$.  For every site $z$, all small enough $\eta>0$ and every $\varepsilon$ and $r$ satisfying $2\eta<\varepsilon<2^4\varepsilon<r\leq 1$,
\begin{equation*}
\mathbf{P}[Z^S(z)\geq r,Z_{\omega_{\eta}^{\lambda}}(z)\leq\varepsilon]\leq O_{\bar{\lambda}}(1)\varepsilon^2\alpha_4(\eta,\varepsilon)\alpha_4(r,1)^{-1}.
\end{equation*}
\end{lemma}
\begin{proof}
The proof is essentially the same as that of \cite[Lemma 8.5]{GPS18a}, which relies on Lemma 8.4 of that paper.  Note that the initial configuration for the dynamics in \cite{GPS18a} is $\omega_{\eta}^0$, while the initial configuration here is $\omega_{\eta}^{\lambda}$.  Since $\bar{\lambda}$ is fixed and the scales are not larger than 1, one can use the quasi-multiplicativity (\ref{e71}) and the fact that the alternating 4-arm probabilities for $\omega_{\eta}^{\lambda}$ are comparable to the critical ones (see (\ref{e68})).  This allows us to follow the proofs of Lemmas 8.4 and 8.5 in \cite{GPS18a} to obtain Lemma \ref{l4}.
\end{proof}

\begin{proof}[Proof of Lemma \ref{l3}]
First, consider any $\lambda<\lambda'$ and $\varepsilon>0$.  Sample $\omega_{\eta}^{\lambda}$, and let $\tilde{\omega}_{\eta}^{\lambda'}$ be the configuration that is open on every site in $X$ and is equal to $\omega_{\eta}^{\lambda}$ elsewhere.  Define the event
\begin{equation*}
\mathcal{E}_1:=\left\{
\begin{aligned}
&\mbox{for any square $B\subset\Lambda_{M-s}$, if there are two $\lambda$-clusters $\mathcal{C},\mathcal{C}'$ of diameter at least $\delta$}\\
&\mbox{completely contained in $B$, and there is an open path of $\omega_{\eta}^{\lambda'}$ connecting $\mathcal{C}$ and $\mathcal{C}'$}\\
&\mbox{in $B$, then the two clusters are also connected by an open path of $\tilde{\omega}_{\eta}^{\lambda'}$ in $\tilde{B}:=B+\Lambda_{s/2}$}
\end{aligned}
\right\}
\end{equation*}
We claim that for any fixed $-\infty<\lambda<\lambda'<\infty$, $M\geq 1,\alpha,s\in(0,1/2)$ and $\delta\in(0,s/2)$, if $\varepsilon>0$ is small enough, then for all small enough $\eta>0$,
\begin{equation}\label{e1}
\mathbf{P}\left[\mathcal{E}_1\right]\geq 1-\alpha/3.
\end{equation}

To show (\ref{e1}), we shall use the proof of \cite[Proposition 8.6]{GPS18a}, with modifications from quad-crossings to connections between clusters.  Assume that the event $\mathcal{E}_1$ does not occur and $B$ is some square realizing $\neg\mathcal{E}_1$.  Let $\{z_1,z_2,\ldots,z_m\}$ denote the sites in $B^{\eta}$.  For $j=0,1,\ldots,m$, let $\omega_j$ denote the configuration that agrees with $\omega_{\eta}^{\lambda'}$ on $\{z_{j+1},\ldots,z_m\}$, and agrees with $\tilde{\omega}_{\eta}^{\lambda'}$ elsewhere.  Then there are two $\lambda$-clusters $\mathcal{C},\mathcal{C}'$ of diameter at least $\delta$ completely contained in $B$, such that there is an open path of $\omega_0\cap B$ connecting them but there are no open paths of $\omega_m\cap\tilde{B}$ connecting them.  Let $k$ be the least index $j$ such that there are no open paths of $\omega_j\cap\tilde{B}$ connecting $\mathcal{C}$ and $\mathcal{C}'$.  Then $\omega_{k-1}$ and $\omega_k$ differ only in the color of $z_k$.  Since a flip of $z_k$ modifies the connectivity between $\mathcal{C}$ and $\mathcal{C}'$ in $\tilde{B}$, and since $z_k,\mathcal{C},\mathcal{C}'\subset B$ and both $\mathcal{C}$ and $\mathcal{C}'$ have diameter at least $\delta$, it follows that $Z_{\omega_{k-1}}(z_k)\geq \delta/2$, and thus $Z^S(z_k)\geq\delta/2$.  If $z$ is any site, then in order to have $z=z_k$, we must have (i) $Z^S(z)\geq\delta/2$, (ii) $Z_{\omega_{\eta}^{\lambda}}(z)\leq 2\varepsilon$, (iii) $z\in S$, and (iv) $z\in\Lambda_{M-s/2}$. Then by Lemma \ref{l4}, the claim (\ref{e1}) follows by summing the bound for $\mathbf{P}[z_k=z]$ over all sites $z$ in $\Lambda_{M-s/2}$.

Let $\mathcal{E}_2$ denote the event that all open clusters of $\omega_{\eta}^{\lambda}$ on $\mathbb{T}_M^2$ have diameter less than $s/2$.  By (\ref{e76}), given any fixed $M\geq 1$ and $\alpha,s\in(0,1/2)$, if $\lambda$ is small enough, then for all $\eta>0$ small enough, we have
\begin{equation}\label{e2}
\mathbf{P}[\mathcal{E}_2]\geq 1-\alpha/3.
\end{equation}

Let $\mathcal{E}_3$ be the event that all the elements of $X$ are distinct and there are no neighboring sites in $X$.  By (\ref{num1}) of Proposition \ref{p6}(ii), we know that for any fixed $-\infty<\lambda<\lambda'<\infty$, $M\geq 1,\alpha\in(0,1/2)$ and $\varepsilon>0$, for all small enough $\eta>0$, we have
\begin{equation}\label{e41}
\mathbf{P}[\mathcal{E}_3]\geq 1-\alpha/3.
\end{equation}

Suppose that the event $\mathcal{E}_1\cap\mathcal{E}_2\cap\mathcal{E}_3$ occurs.  Observe that an open path $\gamma$ of $\tilde{\omega}_{\eta}^{\lambda'}\cap\tilde{B}$ which connects $\mathcal{C}$ and $\mathcal{C}'$ is composed of sub-paths in $\lambda$-clusters and sites in $X$.
Since all the $\lambda$-clusters intersecting $\gamma$ have diameter less than $s/2$, they are contained in $B+\Lambda_s$.  Moreover, there are no neighboring sites in $X\cap\gamma$.  Thus, by using the $\lambda$-clusters intersecting $\gamma$ and the sites in $X\cap\gamma$, we can extract the chain satisfying the conditions in Lemma \ref{l3}.  Then the lemma follows from this argument together with (\ref{e1}), (\ref{e2}) and (\ref{e41}).
\end{proof}

\section{Good embeddings of enhanced networks and cut-off loop ensembles}\label{routergraph}
In \cite{GPS18b}, the authors used enhanced networks to build the cut-off trees, with the embedded edges of these trees being straight line segments; such embeddings are sufficient for the topology of convergence of their models.  In this section, embedded enhanced networks will be defined more carefully than those in Section \ref{s2} and \cite{GPS18b}, since we want to use some ``good" embeddings to define cut-off loop ensembles as their boundary loops, such that whenever good embeddings are close to each other, then so are their cut-off loop ensembles, and moreover, the cut-off loop ensembles are the macroscopic approximations to the collections of cluster boundary loops for the discrete percolation configurations.

Let us start with formal definitions of embeddings of a graph on the torus $\mathbb{T}_M^2$ and on the plane.
\begin{definition}
A \textbf{$\mathbb{T}_M^2$ (resp. plane) graph} is a finite graph drawn on $\mathbb{T}_M^2$ (resp. plane) so that the vertices are all distinct
points, and every edge $(u,v)$ is represented by a simple curve with endpoints $u$ and $v$ so that each (embedded) edge minus its endpoints is disjoint from all the vertices and all the other edges.
\end{definition}

Recall that $\mathcal{Q}_M^k$ denotes the family of polygonal quads in $\mathbb{T}_M^2$ defined in Definition \ref{def9}, and $X=X(\textsf{PPP}_{\bar{\lambda}}^{\varepsilon})$ denotes the set of pivotals of $\textsf{PPP}_{\bar{\lambda}}^{\varepsilon}$ as in Definition \ref{def8}.
For the $r$-mesoscopic enhanced $\varepsilon$-network $\textsf{EN}_M^{r\text{-}\textsf{meso}}(\omega^{\lambda},\textsf{PPP}_{\bar{\lambda}}^{\varepsilon})$ defined in Definition \ref{def10}, we replace each of its router-vertex by the center point of the router, and let $V$ denote the set of these center points, with the original edges replaced accordingly by the edges between the points in $V$ and $X$.  \emph{The resulting graph will also be denoted by} $\textsf{EN}_M^{r\text{-}\textsf{meso}}(\omega^{\lambda},\textsf{PPP}_{\bar{\lambda}}^{\varepsilon})$ with a slight abuse of notation.  Thus, its vertex set is $V\cup X$ when it is nonempty.  We let $E$ denote its edge set.  For convenience, for an edge $(v,x)$ in $E$, we always let $v\in V$ and $x\in X$.  For $v\in V$, denote by $B_v$ the $r$-mesoscopic router centered at $v$, and by $\hat{B}_v$ the box centered at $v$ as in Definition \ref{subrouter}.  Note that by definition, all $\hat{B}_v$'s are strictly away from each other.  For $x\in X$, let $B_x^r$ be the box as in Definition \ref{coverings}, and for $r\leq\hat{r}$ in $2^{-\mathbb{N}}$, let $B_x^{r,\hat{r}}$ be the box of side-length $\hat{r}$ centered around $B_x^r$.

\begin{figure}
\centering
\includegraphics[height=0.48\textwidth]{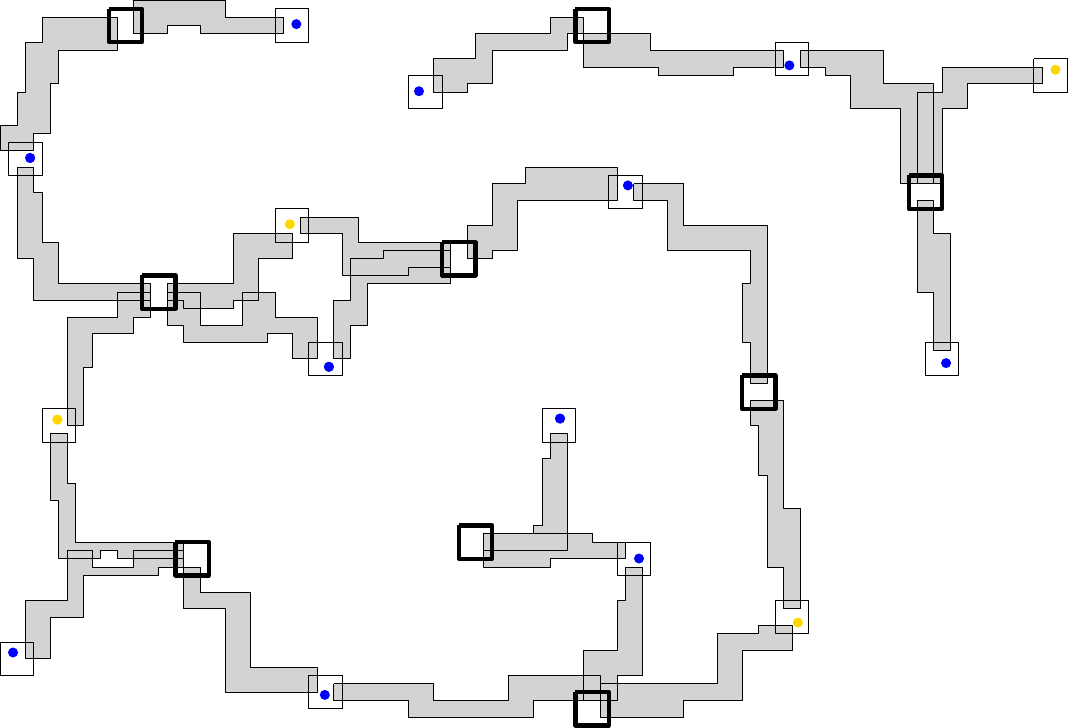}
\caption{Illustration of an $(s,k,\hat{r})$-representation of an enhanced network.  The quads in $\{Q_{v,x}\}_{(v,x)\in E}$ are drawn in light gray. The boxes in $\{\hat{B}_v\}_{v\in V}$ are shown using bold lines. The blue and yellow dots represent the $\lambda$-open and $\lambda$-closed pivotals in $X$, respectively; the squares that contain the dots represent the boxes in $\{\mbox{colored }B_x^r\}_{x\in X}$.
}\label{representation}
\end{figure}

To construct an appropriate embedding of enhanced network, we first represent it by collections of appropriate quads and boxes (see Figure \ref{representation}):
\begin{definition}[A representation of enhanced network]\label{def5}
For $k\in\mathbb{N},r\leq\hat{r}$ in $2^{-\mathbb{N}}$ and $s\geq 2\hat{r}$, an \textbf{$(s,k,\hat{r})$-representation} of a nonempty $\textsf{EN}_M^{r\text{-}\textsf{meso}}(\omega^{\lambda},\textsf{PPP}_{\bar{\lambda}}^{\varepsilon})$ is a collection of quads $\{Q_{v,x}\}_{(v,x)\in E}$ together with the collections of boxes $\{\hat{B}_v\}_{v\in V}$ and $\{\mbox{colored }B_x^r\}_{x\in X}$, which satisfies the following properties:
\begin{itemize}
  \item For each $x\in X$, the box $B_x^r$ is colored open or closed according to the color of $x$.
  \item For each $v\in V$, $\diam(\hat{B}_v)\leq s/2$, and for each $(v,x)\in E$, $\diam(Q_{v,x})\leq s$.
  \item For each $(v,x)\in E$, the quad $Q_{v,x}$ belongs to $\mathcal{Q}_M^k$.
  \item For each $(v,x)\in E$, we have $\omega^{\lambda}\in\boxminus_{Q_{v,x}}$, and, $\partial_1 Q_{v,x}$ is strictly inside $\hat{B}_v$, $\partial_3 Q_{v,x}$ is strictly inside $B_{x}^r$, and $Q_{v,x}$ is strictly away from all the boxes in ${\{B_{x'}^{r,\hat{r}}:x'\in X\setminus\{x\}\}}$ and $\{\hat{B}_{v'}: v'\in V\setminus\{v\}\}$.
  \item Any two quads $Q_{v,x}$ and $Q_{v',x'}$ with $v\neq v'$ (it is possible that $x=x'$) are disjoint.
\end{itemize}
If $\textsf{EN}_M^{r\text{-}\textsf{meso}}(\omega^{\lambda},\textsf{PPP}_{\bar{\lambda}}^{\varepsilon})$ is empty, then we let the empty set be its $(s,k,\hat{r})$-representation.
\end{definition}

As in Definitions \ref{coverings} and \ref{subrouter}, for any $r\in 2^{-\mathbb{N}}$, we let $B^r(\mathbb{T}_M^2)$ be the finite covering of $\mathbb{T}_M^2$ by overlapping $r$-squares.  For any $r\leq\hat{r}\in 2^{-\mathbb{N}}$, denote by $B^{r,\hat{r}}(\mathbb{T}_M^2)$ the collection of squares of side-length $\hat{r}$ centered at the squares of $B^r(\mathbb{T}_M^2)$.  Write $B^{r,\bullet}(\mathbb{T}_M^2):=\bigcup_{\hat{r}\geq r,\hat{r}\in 2^{-\mathbb{N}}}B^{r,\hat{r}}(\mathbb{T}_M^2)$.  In order to prove that with high probability $\textsf{EN}_M^{r\text{-}\textsf{meso}}(\omega^{\lambda},\textsf{PPP}_{\bar{\lambda}}^{\varepsilon})$ in Proposition \ref{p6}(ii) has an $(s,k,\hat{r})$-representation, we need the following lemma.
\begin{lemma}\label{l16}
Fix any $\lambda\in\mathbb{R}$, $M\geq 1$ and $r<\hat{r}$ in $2^{-\mathbb{N}}$.  For any $\alpha>0$, there exists $\delta\in(0,r/40)$, such that for all $\eta>0$ small enough, with probability at least $1-\alpha$ the following event holds.

Consider any sub-collections of squares $\mathscr{B}=\{B_1,\ldots,B_I\}\subset B^{r,\bullet}(\mathbb{T}_M^2)$, $\mathscr{B}'=\{B_1',\ldots,B_J'\}\subset B^r(\mathbb{T}_M^2)$ and $\tilde{\mathscr{B}}'=\{\tilde{B}_1',\ldots,\tilde{B}_J'\}\subset B^{r,\hat{r}}(\mathbb{T}_M^2)$ with each $\tilde{B}_j'$ centered at $B_j'$ and with all the squares of $\mathscr{B}\cup\tilde{\mathscr{B}}'$ at distance at least $r/2$ from each other.  Let $B_i^-$ be the square having the same center as $B_i$, with side-length equal to that of $B_i$ minus $2\delta$; write $\mathscr{B}^-:=\{B_1^-,\ldots,B_I^-\}$.  Define $B_j'^-$ and $\mathscr{B}'^-$ similarly.  Let $E$ be any sub-collection of ordered pairs in $\{(i,j):1\leq i\leq I,1\leq j\leq J\}$.  If
\begin{enumerate}
\item there are no open paths of $\omega_{\eta}^{\lambda}$ in $\mathbb{T}_M^2\setminus(\mathscr{B}\cup\mathscr{B}')$ connecting two distinct squares in $\mathscr{B}$, \label{b1}
\item and there exists a collection of open paths $\{\gamma_{i,j}\}_{(i,j)\in E}$ of $\omega_{\eta}^{\lambda}$ in $\mathbb{T}_M^2\setminus(\mathscr{B}\cup\mathscr{B}')$ such that each $\gamma_{i,j}$ connects $B_i$ and $B_j'$ and does not intersect any squares in $\tilde{\mathscr{B}}'\setminus\{\tilde{B}_j'\}$, and for any two distinct $(i,j),(i',j')\in E$ with $i\neq i'$ (it is possible that $j=j'$), the paths $\gamma_{i,j}$ and $\gamma_{i',j'}$ are disjoint, \label{b2}
\end{enumerate}
then there exists a collection of open paths $\{\hat{\gamma}_{i,j}\}_{(i,j)\in E}$ in $\mathbb{T}_M^2\setminus(\mathscr{B}^-\cup\mathscr{B}'^-)$ such that each $\hat{\gamma}_{i,j}$ connects $B_i^-$ and $B_j'^-$ and is at distance at least $\delta$ from all the squares in $(\mathscr{B}\setminus\{B_i\})\cup(\tilde{\mathscr{B}}'\setminus\{\tilde{B}_j'\})$, and for any two distinct $(i,j),(i',j')\in E$ with $i\neq i'$, the path $\hat{\gamma}_{i,j}$ is at distance at least $\delta$ from $\hat{\gamma}_{i',j'}$.
\end{lemma}

\begin{proof}
The proof is somewhat akin to that of \cite[Lemma 2.9]{GPS13} for the alternating 4-arm case.  Assume that both (\ref{b1}) and $(\ref{b2})$ happen.  For $(i,j)\in E$, let $\mathcal{C}_{i,j}$ denote the open cluster for $\omega_{\eta}^{\lambda}$ restricted to $\mathbb{T}_M^2\setminus(\mathscr{B}^-\cup\mathscr{B}'^-)$ which contains $\gamma_{i,j}$.  Let $\mathcal{E}$ be the event that both (\ref{b1}) and $(\ref{b2})$ occur and furthermore, for any two distinct $(i,j),(i',j')\in E$ with $i\neq i'$, $\mathcal{C}_{i,j}$ and $\mathcal{C}_{i',j'}$ are distinct, and for each $(i,j)\in E$ we can choose an open path $\hat{\gamma}_{i,j}$ in $\mathcal{C}_{i,j}$, such that it connects $B_i^-$ and $B_j'^-$ and is at distance at least $\delta$ from all the squares in $(\mathscr{B}\setminus\{B_i\})\cup(\tilde{\mathscr{B}}'\setminus\{\tilde{B}_j'\})$.

Now we make the following observation: If both (\ref{b1}) and $(\ref{b2})$ occur but $\mathcal{E}$ does not occur, then there must exist a polychromatic 3-arm event from a $2\delta$-square touching some $B\in\mathscr{B}^-\cup\mathscr{B}'^-\cup\tilde{\mathscr{B}}'$ to distance $r/4$ for $\omega_{\eta}^{\lambda}$ restricted to $\mathbb{T}_M^2\setminus B$.  Note that the sides of the squares in $\mathscr{B},\mathscr{B}'$ and $\tilde{\mathscr{B}}'$ are composed of some sides of the $O((M/r)^2)$ squares in $B^r(\mathbb{T}_M^2)$.  Then, by (\ref{e70}) and (\ref{e72}), it is standard (see the left-hand picture of Figure \ref{boundaryarm} for the case that the above $2\delta$-square is close to a corner of $B$) to show that there exist absolute constants $C_1,\varepsilon_3>0$ and a constant $C_{\lambda}>0$ depending only on $\lambda$, such that for all $\eta>0$ small enough the probability of the above event is uniformly bounded by
\begin{align*}
C_{\lambda}\left(\frac{M}{r}\right)^2\left(\left(\frac{r}{\delta}\right)\left(\frac{r}{\delta}\right)^{-2}+\sum_{k=1}^{\lfloor r/\delta\rfloor}k^{-2}\left(\frac{r}{k\delta}\right)^{-(1+\varepsilon_3)}\right)\leq C_1C_{\lambda}\left(\frac{M}{r}\right)^2\frac{\delta}{r}.
\end{align*}

\begin{figure}
\centering
\includegraphics[height=0.48\textwidth]{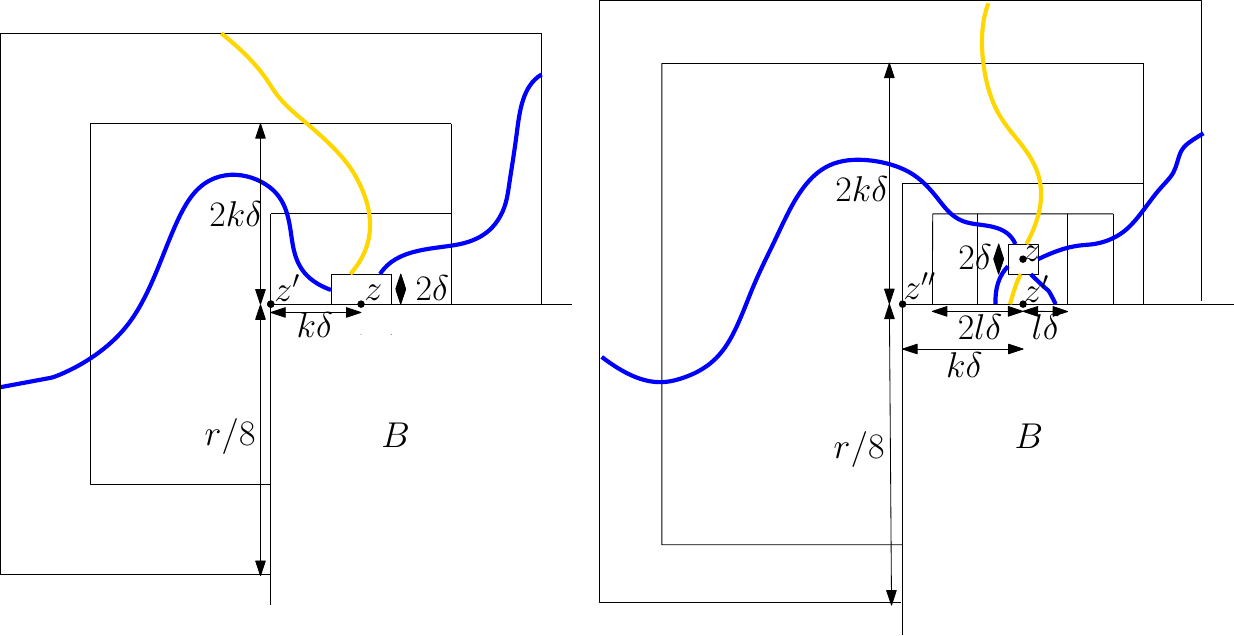}
\caption{\emph{Left}: A polychromatic 3-arm event $\mathcal{A}_{\sigma}(z;2\delta,r/4)$ in $\mathbb{T}_M^2\setminus B$ with $z$ lying in a side of $B$ and at distance $k\delta$ from the closest corner $z'$ of $B$, produces a half-plane 3-arm event in $A(z;2\delta,k\delta)\setminus B$ and a $3\pi/2$-wedge 3-arm event in $A(z';2k\delta,r/8)\setminus B$, where $3\delta\leq k\delta\leq r/20$.  \emph{Right}:  A polychromatic 6-arm event $\mathcal{A}_{\sigma}(z;\delta,l\delta)$ and a polychromatic 3-arm event $\mathcal{A}_{\sigma}(z;\delta,r/4)$ in $\mathbb{T}_M^2\setminus B$, with the latter producing a half-plane 3-arm event in $A(z';2l\delta,k\delta)\setminus B$ and a
$3\pi/2$-wedge 3-arm event in $A(z'';2k\delta,r/8)\setminus B$, where $z'$ is the closest boundary point of $B$ from $z$, $z''$ is the closest corner of $B$ from $z$, and $2\delta\leq l\delta<k\delta/2\leq r/40$.
}\label{boundaryarm}
\end{figure}
Consider the event that $\mathcal{E}$ occurs but there exist two distinct $(i,j),(i',j')\in E$ with $i\neq i'$, such that both $\hat{\gamma}_{i,j}$ and $\hat{\gamma}_{i',j'}$ intersect a $\delta$-square at distance $\delta'$ from the closest square $B\in\mathscr{B}^-\cup\mathscr{B}'^-$.  If $\delta'\geq r/8$, then we have a polychromatic 6-arm event from the $\delta$-square to distance $r/8$; if $\delta'\leq r/8$, then we have a polychromatic 3-arm event from the $\delta$-square to distance $r/4$ in $\mathbb{T}_M^2\setminus B$, and in addition a polychromatic 6-arm event from the $\delta$-square to distance $\delta'$ if $\delta'>0$.  By (\ref{e67}), (\ref{e70}) and (\ref{e72}), it is standard (see the right-hand picture of Figure \ref{boundaryarm} for the case that the above $\delta$-square is close to a corner of $B$) to show that there exist absolute constants $C_2>0,\varepsilon_0\in(0,1)$ and a constant $C_{\lambda}'>0$ depending only on $\lambda$, such that for all $\eta>0$ small enough the above event happens with probability at most
\begin{align*}
&C_{\lambda}'\left(\frac{M}{r}\right)^2\left(\left(\frac{r}{\delta}\right)^2\left(\frac{r}{\delta}\right)^{-(2+\varepsilon_0)}+\sum_{l=1}^{\lfloor r/\delta\rfloor}\sum_{k=l}^{\lfloor r/\delta\rfloor}l^{-(2+\varepsilon_0)}\left(\frac{k}{l}\right)^{-2}\left(\frac{r}{k\delta}\right)^{-(1+\varepsilon_0)}\right.\\
&\qquad\left.+\sum_{l=1}^{\lfloor r/\delta\rfloor}l\cdot l^{-(2+\varepsilon_0)}\left(\frac{r}{l\delta}\right)^{-(1+\varepsilon_0)}+\left(\frac{r}{\delta}\right)\sum_{l=1}^{\lfloor r/\delta\rfloor}l^{-(2+\varepsilon_0)}\left(\frac{r}{l\delta}\right)^{-2}\right)\leq C_2C_{\lambda}'\left(\frac{M}{r}\right)^2\left(\frac{\delta}{r}\right)^{\varepsilon_0}.
\end{align*}

In summary, for all $\eta>0$ small enough, with probability at least $1-O_{\lambda,M,r}(\delta^{\varepsilon_0})$, for any sub-collections of squares $\mathscr{B},\mathscr{B}'$ and $\tilde{\mathscr{B}}'$ in Lemma \ref{l16} which satisfy (\ref{b1}) and (\ref{b2}), we can find the desired collection of open paths $\{\hat{\gamma}_{i,j}\}_{(i,j)\in E}$.  The lemma then follows by choosing $\delta$ small enough.
\end{proof}

The following proposition is a refinement of Proposition \ref{p6}(ii).  Roughly speaking, it says that with high probability, the enhanced networks for $\omega_{\eta}^{\lambda}$ and $\omega_{\infty}^{\lambda}$ have the same representations when $\lambda$ is very negative, and moreover, any open path of $\omega_{\eta}^{\lambda}$ which connects a pair of boxes $\hat{B}_v,B_x^r$ of a representation is connected to $x$ in $B_x^{r,\hat{r}}$.

\begin{proposition}[Existence of $(s,k,\hat{r})$-representations and connections to the pivotals of $X$]\label{p5}
Fix any $\lambda'\in\mathbb{R},M\geq 1$ and $s,\alpha\in(0,1/10)$.  Then there exists $\lambda_0(M,\lambda',s,\alpha)<0$, such that for each $\lambda\leq\lambda_0$, any $\varepsilon\in(0,s)$ and $r_{\alpha}=r_{\alpha}(M,\bar{\lambda},\varepsilon)$ in Proposition \ref{p6}(ii), there exist scales $r<\hat{r}<r_{\alpha}\wedge s/2$ in $2^{-\mathbb{N}}$ and $k\in\mathbb{N}$ such that the following holds.

In any coupling with $\omega_{\eta}^{\lambda}\xrightarrow{a.s.}\omega_{\infty}^{\lambda}$ in $\mathbb{T}_M^2$, for all sufficiently small $\eta>0$ there is a coupling of $\textnormal{\textsf{PPP}}_{\bar{\lambda}}^{\varepsilon}(\omega_{\eta}^{\lambda})$ and $\textnormal{\textsf{PPP}}_{\bar{\lambda}}^{\varepsilon}(\omega_{\infty}^{\lambda})$ such that with probability at least $1-\alpha$, the following event $\mathcal{E}$ occurs: (\ref{num1})-(\ref{num4}) in Proposition \ref{p6} and (\ref{event1}), (\ref{event2}) below hold simultaneously.
\begin{enumerate}[(i)]
  \item $\textnormal{\textsf{EN}}_M^{r\text{-} \textnormal{\textsf{meso}}}(\omega_{\eta}^{\lambda},\textnormal{\textsf{PPP}}_{\bar{\lambda}}^{\varepsilon})$ and
  $\textnormal{\textsf{EN}}_M^{r\text{-}\textnormal{\textsf{meso}}}(\omega_{\infty}^{\lambda},\textnormal{\textsf{PPP}}_{\bar{\lambda}}^{\varepsilon})$ have the same non-empty set of $(s,k,\hat{r})$-representations.  Moreover, all the open clusters of $\omega_{\eta}^{\lambda}$ have diameter less than $s/2$.  \label{event1}
  \item For each pivotal $x\in X(\omega_{\eta}^{\lambda})$, every open path of $\omega_{\eta}^{\lambda}$ which connects $\partial B_x^r$ to $\partial B_x^{r,\hat{r}}$
  is connected to $x$ by an open path of $\omega_{\eta}^{\lambda}$ in $B_x^{r,\hat{r}}$. \label{event2}
\end{enumerate}
\end{proposition}
\begin{proof}
Take $r\ll\hat{r}\ll r_{\alpha}\wedge s/2$ in $2^{-\mathbb{N}}$, with $r$ and $\hat{r}$ to be fixed later.  By (\ref{e76}) and Proposition \ref{p6}, there exists $\lambda_0(M,\lambda',s,\alpha)<0$, such that for any $\lambda\leq\lambda_0$, in any coupling with $\omega_{\eta}^{\lambda}\xrightarrow{a.s.}\omega_{\infty}^{\lambda}$ in $\mathbb{T}_M^2$, for all $\eta>0$ small enough, there is a coupling of $\textnormal{\textsf{PPP}}_{\bar{\lambda}}^{\varepsilon}(\omega_{\eta}^{\lambda})$ and $\textnormal{\textsf{PPP}}_{\bar{\lambda}}^{\varepsilon}(\omega_{\infty}^{\lambda})$ such that with probability at least $1-3\alpha/4$, all the events (\ref{num1})-(\ref{num4}) in Proposition \ref{p6}(ii) occur simultaneously, and moreover, all the open clusters of $\omega_{\eta}^{\lambda}$ have diameter less than $s/2$. We denote this event by $\mathcal{E}_1$.

Write $\mathscr{B}_X:=\{B_x^r\}_{x\in X},\tilde{\mathscr{B}}_X:=\{B_x^{r,\hat{r}}\}_{x\in X}$ and $\hat{\mathscr{B}}_V:=\{\hat{B}_v\}_{v\in V}$ for short.  Define the event
\begin{equation*}
\mathcal{F}:=\left\{
\begin{aligned}
&\mbox{for each $(v,x)\in E$, there exists an open path of $\omega_{\eta}^{\lambda}$}\\
&\mbox{connecting $\hat{B}_v$ and $B_x^r$ in $\mathbb{T}_M^2\setminus(\tilde{\mathscr{B}}_{X\setminus\{x\}}\cup\hat{\mathscr{B}}_V)$}
\end{aligned}
\right\}.
\end{equation*}

If $\mathcal{E}_1$ occurs, then for every $v\in V$, there is a unique open cluster for the restriction of $\omega_{\eta}^{\lambda}$ to $\mathbb{T}_M^2\setminus(\mathscr{B}_X\cup\hat{\mathscr{B}}_V)$ which connects $\hat{B}_v$ to $B_x^r$'s for all $(v,x)\in E$.  We denote this cluster by $\mathcal{C}_v$.

Assume that $\mathcal{E}_1$ occurs but $\mathcal{F}$ does not occur.  Then there is some $(v,x)\in E$ such that there does not exist an open path of $\omega_{\eta}^{\lambda}$ connecting $\hat{B}_v$ and $B_x^r$ in $\mathbb{T}_M^2\setminus(\tilde{\mathscr{B}}_{X\setminus\{x\}}\cup\hat{\mathscr{B}}_V)$, and thus there is some $x^*\in X\setminus\{x\}$ such that the cluster $\mathcal{C}_v$ minus the square $B_{x^*}^{r,\hat{r}}$ is disconnected.  So there exist four arms of alternating colors in the annulus $B_{x^*}^{r,r_{\alpha}}\setminus B_{x^*}^{r,\hat{r}}$ connecting the two boundary pieces of this annulus, with the two open arms lying in $\mathcal{C}_v$ and connected in $B_{x^*}^{r,r_{\alpha}}\setminus B_{x^*}^r$.  Moreover, among the two closed arms there must be one that does not lie in a radially crossing closed cluster for the restriction of $\omega_{\eta}^{\lambda}$ to $B_{x^*}^{r,r_{\alpha}}\setminus B_{x^*}^r$, since otherwise the two open arms would not be connected in $B_{x^*}^{r,r_{\alpha}}\setminus B_{x^*}^r$; the $\lambda$-closed cluster restricted to $B_{x^*}^{r,r_{\alpha}}\setminus B_{x^*}^r$ which contains this closed arm is denoted by $\tilde{\mathcal{C}}$.  See the left-hand picture of Figure \ref{sixarms}.  Since $x^*$ is an $\varepsilon$-important point of $\omega_{\eta}^{\lambda}$, there are two distinct radially crossing $\lambda$-closed clusters $\tilde{\mathcal{C}}_1,\tilde{\mathcal{C}}_2$ in $B_{x^*}^{r,r_{\alpha}}\setminus B_{x^*}^r$.  Then $\tilde{\mathcal{C}},\tilde{\mathcal{C}}_1$ and $\tilde{\mathcal{C}}_2$ imply a polychromatic 6-arm event in $B_{x^*}^{r,r_{\alpha}}\setminus B_{x^*}^{r,\hat{r}}$.  By (\ref{e67}), this event happens with probability at most $C_{\lambda}(M/\hat{r})^2(\hat{r}/r_{\alpha})^{2+\varepsilon_6}$ for all $\eta>0$ small enough, where $C_{\lambda}>0$ is a constant depending only on $\lambda$.

\begin{figure}
\centering
\includegraphics[height=0.4\textwidth]{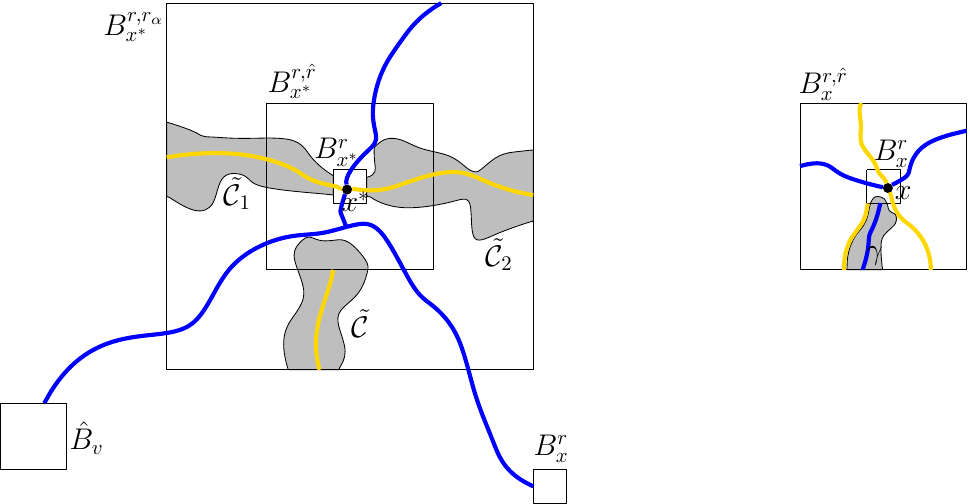}
\caption{\emph{Left}: The $\lambda$-closed clusters $\tilde{\mathcal{C}},\tilde{\mathcal{C}}_1$ and $\tilde{\mathcal{C}}_2$ (drawn in gray) produce a polychromatic 6-arm event in $B_{x^*}^{r,r_{\alpha}}\setminus B_{x^*}^{r,\hat{r}}$. \emph{Right}:  The $\lambda$-open path $\gamma$ that connects $\partial B_x^r$ to $\partial B_x^{r,\hat{r}}$, with $x\in X$ and $\gamma$ not connected by $\lambda$-open paths in $B_x^{r,\hat{r}}$ to $x$, implies a polychromatic 6-arm event from $\partial B_x^r$ to $\partial B_x^{r,\hat{r}}$.  The $\lambda$-open cluster restricted to $B_x^{r,\hat{r}}$ which contains $\gamma$ is drawn in gray.
}\label{sixarms}
\end{figure}
Assume that (\ref{event2}) does not occur.  Then there is some $x\in X$ and a $\lambda$-open path $\gamma$ connecting $\partial B_x^r$ to $\partial B_x^{r,\hat{r}}$, such that $\gamma$ is not connected by $\lambda$-open paths in $B_x^{r,\hat{r}}$ to $x$.  This implies a 6-arm event from $\partial B_x^r$ to $\partial B_x^{r,\hat{r}}$ (see the right-hand picture of Figure \ref{sixarms}):  Since $x$ is an $\varepsilon$-important point of $\omega_{\eta}^{\lambda}$, there are two disjoint open arms connecting $x$ to $\partial B_x^{r,\hat{r}}$ and belonging to distinct open clusters of $\omega_{\eta}^{\lambda}$ restricted to $B_x^{r,\hat{r}}\setminus\{x\}$.  These two arms and $\gamma$ induce three disjoint open arms from $\partial B_x^r$ to $\partial B_x^{r,\hat{r}}$ which belong to distinct open clusters of $\omega_{\eta}^{\lambda}$ restricted to $B_x^{r,\hat{r}}\setminus B_x^r$.  The three induced open arms, together with three disjoint closed arms separating them, produces a polychromatic 6-arm event from $\partial B_x^r$ to $\partial B_x^{r,\hat{r}}$.  By (\ref{e67}) again, this event happens with probability at most $C_{\lambda}'(M/r)^2(r/\hat{r})^{2+\varepsilon_6}$ for all $\eta>0$ small enough, where $C_{\lambda}'>0$ is a constant depending only on $\lambda$.

Let $\mathcal{E}_2$ be the event that $\mathcal{E}_1$, (\ref{event2}) and $\mathcal{F}$ occur simultaneously.  By the above arguments, we can choose $r<\hat{r}<r_{\alpha}\wedge s/2$ in $2^{-\mathbb{N}}$, such that for all $\eta>0$ small enough,
\begin{equation}\label{e69}
\mathbf{P}[\mathcal{E}_2]\geq 1-4\alpha/5.
\end{equation}

Take $\delta\in(0,r/40)$.  For each $x\in X$, let $B_x^-$ be the square having the same center as $B_x^r$, with side-length equal to that of $B_x^r$ minus $2\delta$.
Define $\hat{B}_v^-$ similarly for each $v\in V$.  Let $\mathcal{E}_3$ be the event that $\mathcal{E}_2$ and the following event occur simultaneously:  There exists a collection of open paths $\{\gamma_{v,x}\}_{(v,x)\in E}$ of $\omega_{\eta}^{\lambda}$ in $\mathbb{T}_M^2\setminus(\{B_x^-\}_{x\in X}\cup\{\hat{B}_v^-\}_{v\in V})$ such that each $\gamma_{v,x}$ connects $\hat{B}_v^-$ and $B_x^-$ and is at distance at least $\delta$ from all the squares in $\tilde{\mathscr{B}}_{X\setminus\{x\}}\cup\hat{\mathscr{B}}_{V\setminus\{v\}}$, and for any two distinct edges $(v,x),(v',x')\in E$ with $v\neq v'$, the path $\gamma_{v,x}$ is at distance at least $\delta$ from $\gamma_{v',x'}$.

Letting $r,\hat{r}$ be as in (\ref{e69}) and combining (\ref{e69}) and Lemma \ref{l16}, it follows that there exists $\delta=\delta(\lambda,M,r,\hat{r})\in (0,r/40)$, such that for all $\eta>0$ small enough,
\begin{equation}\label{e73}
\mathbf{P}[\mathcal{E}_3]\geq 1-5\alpha/6.
\end{equation}

Let $k$ be the minimal integer such that $2^{-k}<\delta/2^{10}$.  Observe that when $\mathcal{E}_3$ occurs and $\eta>0$ is sufficiently small, then by using $\{\gamma_{v,x}\}_{(v,x)\in E}$ we can obtain a collection of quads $\{Q_{v,x}\}_{(v,x)\in E}$, which, together with $\{\hat{B}_v\}_{v\in V}$ and $\{\mbox{colored }B_x^r\}_{x\in X}$ forms an $(s,k,\hat{r})$-representation of $\textnormal{\textsf{EN}}_M^{r\text{-}\textsf{meso}}(\omega_{\eta}^{\lambda},\textnormal{\textsf{PPP}}_{\bar{\lambda}}^{\varepsilon})$.  Then by (\ref{e73}) and (\ref{e4}), for all $\eta>0$ small enough, we have
\begin{equation*}
\mathbf{P}[\mathcal{E}_3\cap(\ref{event1})]\geq 1-\alpha,
\end{equation*}
which implies the desired proposition.
\end{proof}

Using the representations in Proposition \ref{p5}, we will now construct the same router graph for $\omega_{\eta}^{\lambda}$ and $\omega_{\infty}^{\lambda}$ with high probability; it is a $\mathbb{T}_M^2$ graph, and can be seen as a suitable embedding of $\textsf{EN}_M^{r\text{-}\textsf{meso}}(\omega^{\lambda},\textsf{PPP}_{\bar{\lambda}}^{\varepsilon})$, though this graph has vertex set $V$ while $\textsf{EN}_M^{r\text{-}\textsf{meso}}(\omega^{\lambda},\textsf{PPP}_{\bar{\lambda}}^{\varepsilon})$ has vertex set $V\cup X$.  Furthermore, we will construct the good router graph for $\omega_{\eta}^{\lambda}$.

\begin{figure}
\centering
\includegraphics[height=0.48\textwidth]{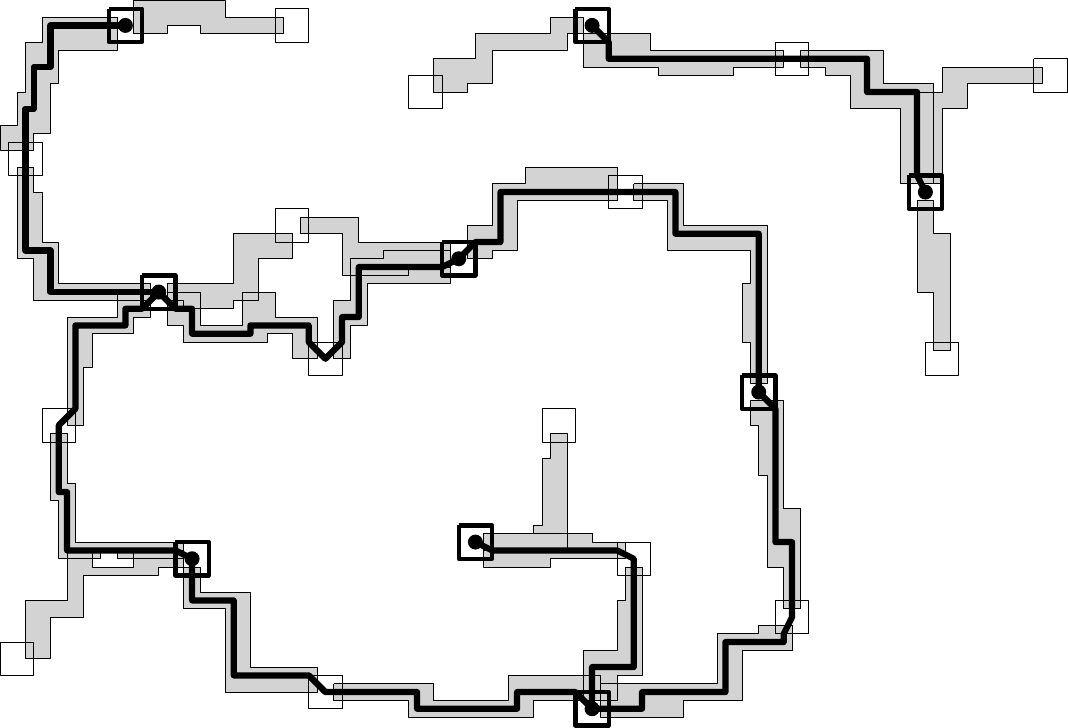}
\caption{Building the router graph $\textnormal{\textsf{G}}^{\bar{\lambda},\varepsilon}$ from the representation of Figure \ref{representation}.  Vertices of $\textnormal{\textsf{G}}^{\bar{\lambda},\varepsilon}$ are indicated by dots, and edges are indicated by black bold lines.
}\label{routerg}
\vspace{1em}
\includegraphics[height=0.48\textwidth]{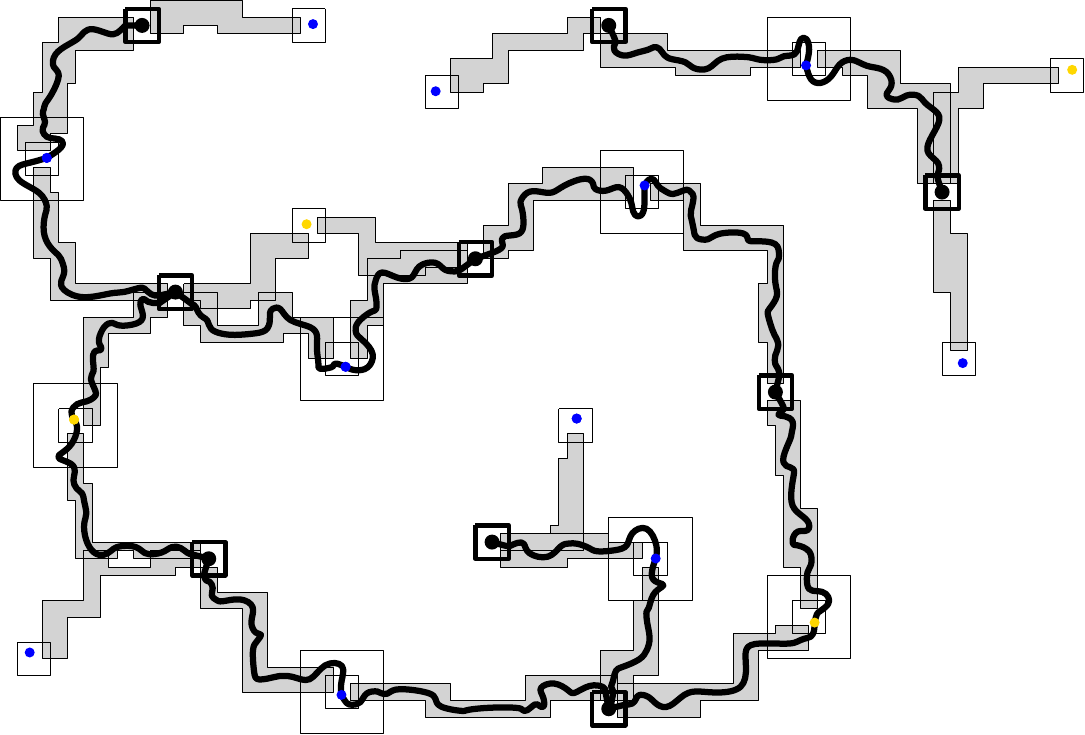}
\caption{Building the good router graph $\tilde{\textnormal{\textsf{G}}}_{\eta}^{\bar{\lambda},\varepsilon}$ from the representation of Figure \ref{representation}.  Vertices of $\tilde{\textnormal{\textsf{G}}}_{\eta}^{\bar{\lambda},\varepsilon}$ are indicated by black dots, and edges are indicated by black bold lines.  Each large square centered at a small square represents some $B_x^{r,\hat{r}}$, with $x$ being a pivotal in $X$ that corresponds to an edge (of $\tilde{\textnormal{\textsf{G}}}_{\eta}^{\bar{\lambda},\varepsilon}$).  The blue and yellow dots represent $\lambda$-open and $\lambda$-closed $\eta$-hexagons centered at the pivotals of $X$, respectively.  Each edge passes through such an $\eta$-hexagon.
}\label{routerg2}
\end{figure}

\begin{definition}[Constructing the router graphs for both $\omega_{\eta}^{\lambda}$ and $\omega_{\infty}^{\lambda}$, and the good router graph for $\omega_{\eta}^{\lambda}$]\label{def7}
Fix any $\lambda'\in\mathbb{R},M\geq 1,s,\alpha\in(0,1/10)$ and $\varepsilon\in(0,s)$.  Let $\lambda,r,\hat{r}$ and $k$ be as in
Proposition \ref{p5}.  Assume that $\textnormal{\textsf{EN}}_M^{r\text{-}\textsf{meso}}(\omega^{\lambda},\textnormal{\textsf{PPP}}_{\bar{\lambda}}^{\varepsilon})$ is nonempty and has $(s,k,\hat{r})$-representations, with all elements of $X$ distinct and having at most two router-neighbors.  We choose one of such $(s,k,\hat{r})$-representations in a deterministic way depending only on their boxes and quads, denote it by $(\{Q_{v,x}\}_{(v,x)\in E},\{\hat{B}_v\}_{v\in V},\{\mbox{colored }B_x^r\}_{x\in X})$, and use it to construct the \textbf{router graph} $\textnormal{\textsf{G}}(\omega^{\lambda},\textnormal{\textsf{PPP}}_{\bar{\lambda}}^{\varepsilon})$ in the following way (see Figure \ref{routerg}):

The vertex set is $V$.  An edge $(v,v')$ is placed between any two distinct vertices $v,v'\in V$ if and only if there exists $x\in X$ such that $(v,x),(v',x)\in E$.
If there are more than one such $x$ for $(v,v')$, then we use the following way to select one: If there are open ones, then we choose an open one in a way depending on the (open) boxes $B_x^r$'s containing such open $x$'s; otherwise, we select one of them relying on the (closed) boxes $B_x^r$'s containing such $x$'s.  Thus, each edge $(v,v')$ corresponds to a unique $x\in X$, and thus corresponds to the box $B_x^r$ and the pair of quads $Q_{v,x},Q_{v',x}$.  For each $v\in V$, let $X_v\subset X$ be the set of pivotals which correspond to the edges of the router graph incident to $v$, and let $U_v:=\hat{B}_v\cup\bigcup_{x\in X_v}Q_{v,x}$.  Let $\tilde{x}$ denote the center point of $B_x^r$.  It is clear that every edge $(v,v')$ can be embedded in $\mathbb{T}_M^2$ such that it is the concatenation of two simple curves $(v\tilde{x})$ and $(\tilde{x}v')$, with $(v\tilde{x})$ strictly inside $U_v\cup B_x^r$ and $(\tilde{x}v')$ strictly inside $U_{v'}\cup B_x^r$, making the embedded graph a $\mathbb{T}_M^2$ graph.  We use a deterministic method depending only on $(\{Q_{v,x}\}_{(v,x)\in E},\{\hat{B}_v\}_{v\in V},\{\mbox{colored }B_x^r\}_{x\in X})$ to realize the embedding.

If the above assumption is not satisfied, let $\textnormal{\textsf{G}}(\omega^{\lambda},\textnormal{\textsf{PPP}}_{\bar{\lambda}}^{\varepsilon})$ be empty.

Write $\textnormal{\textsf{G}}_{\infty}^{\bar{\lambda},\varepsilon}:=\textnormal{\textsf{G}}(\omega_{\infty}^{\lambda},\textnormal{\textsf{PPP}}_{\bar{\lambda}}^{\varepsilon})$ and $\textnormal{\textsf{G}}_{\eta}^{\bar{\lambda},\varepsilon}:=\textnormal{\textsf{G}}(\omega_{\eta}^{\lambda},\textnormal{\textsf{PPP}}_{\bar{\lambda}}^{\varepsilon})$. Let $\mathcal{E}$ be the event in Proposition \ref{p5}.  Notice that when $\mathcal{E}$ holds, we have $\textnormal{\textsf{G}}_{\eta}^{\bar{\lambda},\varepsilon}=\textnormal{\textsf{G}}_{\infty}^{\bar{\lambda},\varepsilon}$.

Assume now that $\mathcal{E}$ holds.  Then we construct the \textbf{good router graph} $\tilde{\textnormal{\textsf{G}}}_{\eta}^{\bar{\lambda},\varepsilon}$ for $\omega_{\eta}^{\lambda}$ in the following.  It has the same graph structure as $\textnormal{\textsf{G}}_{\eta}^{\bar{\lambda},\varepsilon}$ and $\textnormal{\textsf{G}}_{\infty}^{\bar{\lambda},\varepsilon}$, with its vertex set being $V$ but each of its edge $(v,v')$ being embedded in $\mathbb{T}_M^2$ in the following way.  See Figure \ref{routerg2} for an illustration.  As in the above construction, let $x\in X$ be the uniquely chosen pivotal corresponding to $(v,v')$.  Write $H_x$ for the $\eta$-hexagon centered at $x$.  Every edge $(v,v')$ is embedded in $\mathbb{T}_M^2$ such that it is the concatenation of two simple curves $(vx)$ and $(xv')$, with $(vx)$ strictly inside $U_v\cup B_x^{r,\hat{r}}$ and $(xv')$ strictly inside $U_{v'}\cup B_x^{r,\hat{r}}$, making the embedded graph a $\mathbb{T}_M^2$ graph, and moreover, $(v,v')$ intersects $\partial\hat{B}_v$ (resp. $\partial\hat{B}_{v'}$) at just a single point $w$ (resp. $w'$), and the sub-curve of $(v,v')$ from $w$ to $w'$ is strictly inside the union of $H_x$ and the open $\eta$-hexagons of $\omega_{\eta}^{\lambda}\setminus X$, with $(v,v')$ intersecting $\partial H_x$ at exactly two points.

If the event $\mathcal{E}$ does not hold, let $\tilde{\textnormal{\textsf{G}}}_{\eta}^{\bar{\lambda},\varepsilon}$ be empty.
\end{definition}

Let $G$ be a $\mathbb{T}_M^2$ graph with vertex set $V$.  A connected component of $\mathbb{T}_M^2\setminus G$ is called a \textbf{face} of $G$.  Then the boundaries of the faces of $G$ form a collection of loops in $\mathbb{T}_M^2$; such loops are called the \textbf{boundary loops} of $G$.  For a boundary loop $\ell$ of $G$, let $V(\ell)$ denote the set of vertices in $V$ that are contained in $\ell$.  For a boundary loop of a face of $G$ that is completely contained in $\Lambda_M\subset\mathbb{C}$, we orient the loop so that the points to the left of the loop are in this face.

Using router graphs, we are now going to define cut-off loop ensembles on $\Lambda_M$; the reason that they are defined on $\Lambda_M$ but not on $\mathbb{T}_M^2$ is that we want to avoid the technicalities related to non-contractible loops on $\mathbb{T}_M^2$.

\begin{figure}
\centering
\includegraphics[height=0.46\textwidth]{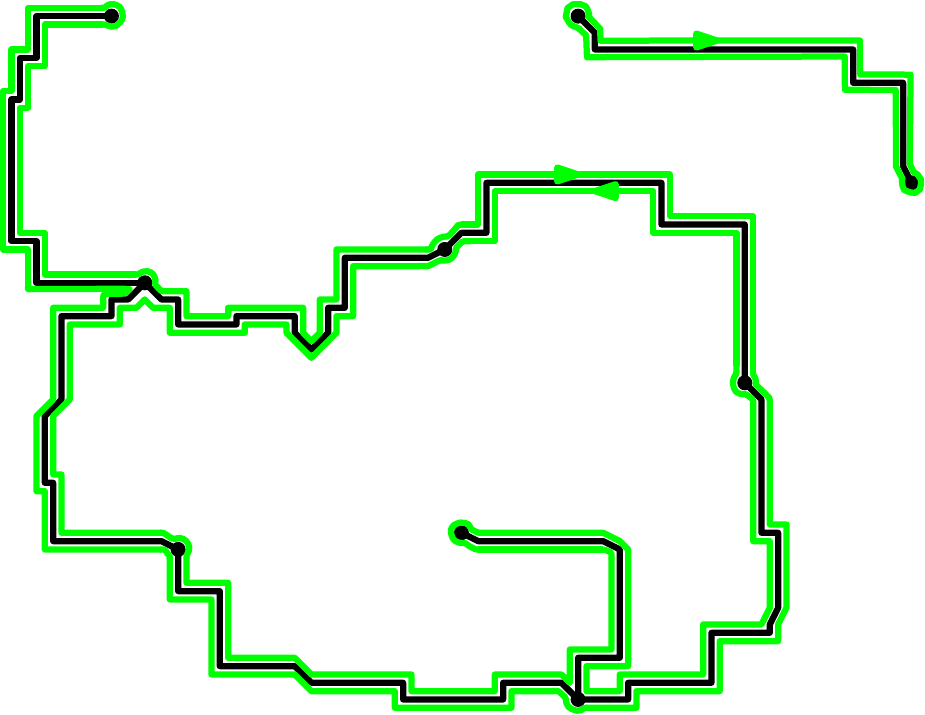}
\caption{Extracting the cut-off loop ensemble $\Gamma^{\bar{\lambda},\varepsilon}$ from the router graph of Figure \ref{routerg}.  Though the three non-trivial loops of $\Gamma^{\bar{\lambda},\varepsilon}$ are not simple here, they are represented by green simple loops for clarity.
}\label{cutoffloop}
\vspace{1em}
\includegraphics[height=0.46\textwidth]{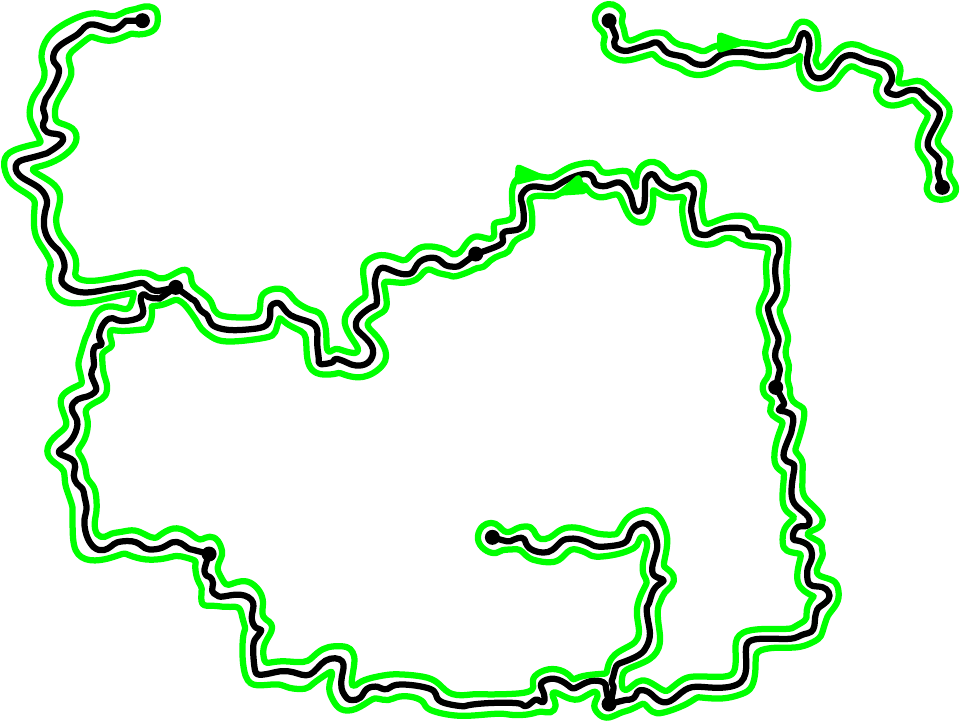}
\caption{Extracting the good cut-off loop ensemble $\tilde{\Gamma}_{\eta}^{\bar{\lambda},\varepsilon}$ from the good router graph of Figure \ref{routerg2}.
}\label{cutoffloop2}
\end{figure}

\begin{definition}[Constructing the cut-off loop ensembles for $\omega_{\eta}^{\lambda}$ and $\omega_{\infty}^{\lambda}$ and the good cut-off loop ensemble for $\omega_{\eta}^{\lambda}$]\label{def2}
Let $\textnormal{\textsf{G}}_{\infty}^{\bar{\lambda},\varepsilon}$, $\textnormal{\textsf{G}}_{\eta}^{\bar{\lambda},\varepsilon}$ and $\tilde{\textnormal{\textsf{G}}}_{\eta}^{\bar{\lambda},\varepsilon}$ be as in Definition \ref{def7}.  We denote by $\Gamma_{\infty}^{\bar{\lambda},\varepsilon}=\Gamma_{\infty}^{\bar{\lambda},\varepsilon}(\textnormal{\textsf{G}}_{\infty}^{\bar{\lambda},\varepsilon})$
the collection of all the oriented boundary loops of $\textnormal{\textsf{G}}_{\infty}^{\bar{\lambda},\varepsilon}$ which are completely contained in $\Lambda_M$ and all the trivial loops in $\Lambda_M$.  The sets of loops $\Gamma_{\eta}^{\bar{\lambda},\varepsilon}=\Gamma_{\eta}^{\bar{\lambda},\varepsilon}(\textnormal{\textsf{G}}_{\eta}^{\bar{\lambda},\varepsilon})$ and $\tilde{\Gamma}_{\eta}^{\bar{\lambda},\varepsilon}=\tilde{\Gamma}_{\eta}^{\bar{\lambda},\varepsilon}(\tilde{\textnormal{\textsf{G}}}_{\eta}^{\bar{\lambda},\varepsilon})$ are defined similarly.  See Figures \ref{cutoffloop} and \ref{cutoffloop2}.  We will call $\Gamma_{\infty}^{\bar{\lambda},\varepsilon}$ and $\Gamma_{\eta}^{\bar{\lambda},\varepsilon}$ the \textbf{cut-off loop ensembles}, and call $\tilde{\Gamma}_{\eta}^{\bar{\lambda},\varepsilon}$ the \textbf{good cut-off loop ensemble}.
\end{definition}

To study the full plane loop ensemble $\Gamma_{\eta}^{\lambda'}$, it is more convenient to compare it to $\tilde{\Gamma}_{\eta}^{\bar{\lambda},\varepsilon}$, which we will do in the next section.  Roughly speaking, Lemma \ref{l14} below says that $\tilde{\Gamma}_{\eta}^{\bar{\lambda},\varepsilon}$ is well-approximated by $\Gamma_{\infty}^{\bar{\lambda},\varepsilon}$.

For $0<R<M$, let $\Gamma_{\infty}^{\bar{\lambda},\varepsilon}(R)$ denote the collection of loops in $\Gamma_{\infty}^{\bar{\lambda},\varepsilon}$ which are contained in $\Lambda_R$.  Define $\tilde{\Gamma}_{\eta}^{\bar{\lambda},\varepsilon}(R)$ similarly.
\begin{lemma}\label{l14}
Fix any $\lambda'\in\mathbb{R},M\geq 2,s,\alpha\in(0,1/10)$ and $\varepsilon\in(0,s)$.  Assume that the event $\mathcal{E}$ in Proposition \ref{p5} occurs.  Then we have:
\begin{itemize}
  \item For each $\tilde{\ell}\in\tilde{\Gamma}_{\eta}^{\bar{\lambda},\varepsilon}(M-1)$, there exists $\ell\in\Gamma_{\infty}^{\bar{\lambda},\varepsilon}$ such that $d(\tilde{\ell},\ell)\leq 12s$.
  \item For each $\ell\in\Gamma_{\infty}^{\bar{\lambda},\varepsilon}(M-1)$, there exists $\tilde{\ell}\in\tilde{\Gamma}_{\eta}^{\bar{\lambda},\varepsilon}$ such that $d(\tilde{\ell},\ell)\leq 12s$.
\end{itemize}
\end{lemma}
\begin{proof}
Let $V$ be the vertex set of $\textnormal{\textsf{G}}_{\infty}^{\bar{\lambda},\varepsilon}$ and let $v\in V$.  As in Definition \ref{def7}, denote by $X_v\subset X$ the set of pivotals corresponding to the edges of $\textnormal{\textsf{G}}_{\infty}^{\bar{\lambda},\varepsilon}$ incident to $v$, and let $U_v:=\hat{B}_v\cup\bigcup_{x\in X_v}Q_{v,x}$.  For $x\in X_v$, write $U_{v,x}:=U_v\cup B_x^{r,\hat{r}}$, and let $F(U_{v,x})$ denote the filling of $U_{v,x}$, which is given by taking the union of $U_{v,x}$ and the components of $\mathbb{T}_{M}^2\setminus U_{v,x}$ which are homeomorphic to a disk.  For each edge $(u,v)$ of $\textnormal{\textsf{G}}_{\infty}^{\bar{\lambda},\varepsilon}$, write $F_{u,v}:=F(U_{u,x})\cup F(U_{v,x})$ with $x\in X$ corresponding to $(u,v)$.  By the definition of $(s,k,\hat{r})$-representation in Definition \ref{def5}, we have $\diam(F(U_{v,x}))\leq 3s$ and $\diam(F_{u,v})\leq 6s$.

\begin{figure}
\centering
\includegraphics[height=0.26\textwidth]{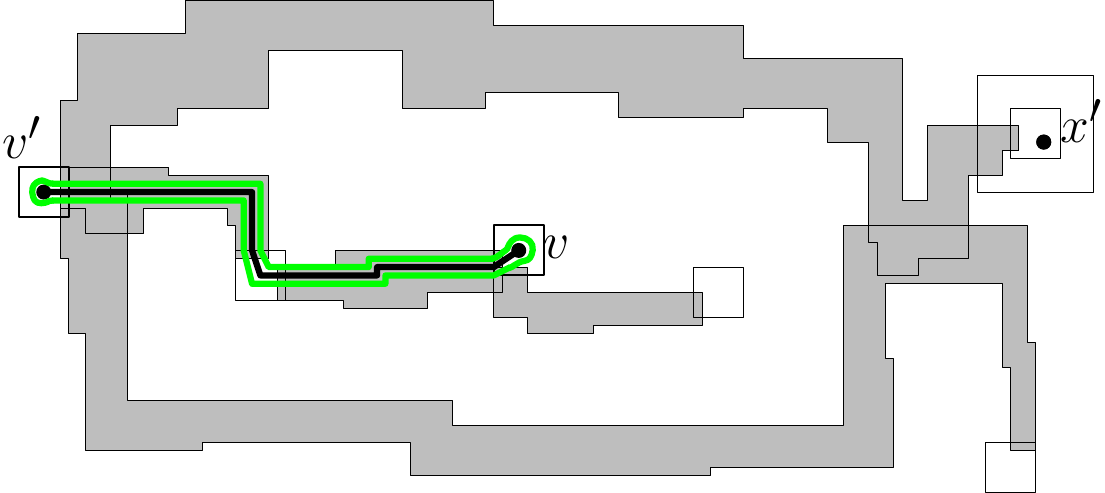}
\caption{Illustration of a vertex $v\in V$ that is contained in the filling $F(U_{v',x'})$, with $v'\in V$ and $x'\in X_{v'}$.  The green loop is an erased dangling loop attached to $v'$.
}\label{dangling}
\end{figure}

Now we construct a subgraph of $\textnormal{\textsf{G}}_{\infty}^{\bar{\lambda},\varepsilon}$:  For every $v\in V$, if there exists $v'\in V\setminus\{v\}$ and $x'\in X_{v'}$ such that $v$ is contained in $F(U_{v',x'})$ (see Figure \ref{dangling}), then we remove $v$ and the edges of $\textnormal{\textsf{G}}_{\infty}^{\bar{\lambda},\varepsilon}$ incident to $v$.  The resulting graph is denoted by $\textsf{G}_{\infty}'$.  The subgraph $\tilde{\textsf{G}}_{\eta}'$ of $\tilde{\textnormal{\textsf{G}}}_{\eta}^{\bar{\lambda},\varepsilon}$ is constructed similarly.  Note that $\textsf{G}_{\infty}'$ and $\tilde{\textsf{G}}_{\eta}'$ have the same (embedded) vertex set, denoted $V'$, and their graph structures coincide.  Observe that each edge $(u,v)$ of $\textsf{G}_{\infty}'$ and the corresponding edge $\widetilde{(u,v)}$ of $\tilde{\textsf{G}}_{\eta}'$ are homotopic relative to their endpoints in $F_{u,v}$, and there are no vertices of $V'\setminus\{u,v\}$ contained in $F_{u,v}$; thus, $(u,v)$ can be continuously deformed into $\widetilde{(u,v)}$ in $F_{u,v}$ such that in the deformation, $u,v$ are fixed and at any time the deformed edge does not intersect any vertex in $V'\setminus\{u,v\}$.  (So $\textsf{G}_{\infty}'$ and $\tilde{\textsf{G}}_{\eta}'$ are isotopic with fixed vertices.)

Let $\Gamma_{\infty}'$ (resp. $\tilde{\Gamma}_{\eta}'$) be the set of trivial loops in $\Lambda_{M-1}$ and the boundary loops $\ell$ of $\textsf{G}_{\infty}'$ (resp. $\tilde{\textsf{G}}_{\eta}'$) satisfying $V'(\ell)\subset\Lambda_{M-1}$.  Note that each non-trivial loop $\ell$ of $\Gamma_{\infty}'$ induces a cyclic ordering $(e_1,\ldots,e_n)$ of the oriented edges $e_1,\ldots,e_n$ of $\textsf{G}_{\infty}'$ traveled along $\ell$; we can identify $\ell$ as this ordering and write $\ell=(e_1,\ldots,e_n)$.  For an oriented edge $\langle u,v\rangle$ of $\textsf{G}_{\infty}'$, let $\widetilde{\langle u,v\rangle}$ be the corresponding oriented edge of $\tilde{\textsf{G}}_{\eta}'$ from $u$ to $v$.

The above argument implies that there is a bijection between the two sets of non-trivial loops of $\Gamma_{\infty}'$ and $\tilde{\Gamma}_{\eta}'$, such that each non-trivial loop $(e_1,\ldots,e_n)$ in $\Gamma_{\infty}'$ is mapped to the non-trivial loop $(\tilde{e}_1,\ldots,\tilde{e}_n)$ in $\tilde{\Gamma}_{\eta}'$.  Since $(u,v),\widetilde{(u,v)}\subset F_{u,v}$ and $\diam(F_{u,v})\leq 6s$ for each edge $(u,v)$ of $\textsf{G}_{\infty}'$, we have
\begin{equation}\label{e66}
d_{\mathcal {L}}(\Gamma_{\infty}',\tilde{\Gamma}_{\eta}')\leq\max\{\diam(F_{u,v}):(u,v)\mbox{ is an edge of $\textsf{G}_{\infty}'$}\}\leq 6s.
\end{equation}

If $\gamma$ is a non-trivial sub-loop of a boundary loop of $\textsf{G}_{\infty}^{\bar{\lambda},\varepsilon}$ and there exists $v'\in V'$ lying in $\gamma$ such that $V(\gamma)\subset F(U_{v',x'})$ for some $x'\in X_{v'}$ (thus all the vertices of $V(\gamma)\setminus\{v'\}$ are removed from $\textsf{G}_{\infty}^{\bar{\lambda},\varepsilon}$ in the construction of $\textsf{G}_{\infty}'$), then we call $\gamma$ an \textbf{erased dangling loop} of $\textsf{G}_{\infty}^{\bar{\lambda},\varepsilon}$ attached to $v'$.  See Figure \ref{dangling}.  Notice that we also have $\gamma\subset F(U_{v',x'})$, and thus $\diam(\gamma)\leq\diam(F(U_{v',x'}))\leq 3s$.  The erased dangling loops of $\tilde{\textsf{G}}_{\eta}^{\bar{\lambda},\varepsilon}$ are defined similarly.

Denote by $\Gamma_{\infty}''$ (resp. $\tilde{\Gamma}_{\eta}''$) the set of trivial loops in $\Lambda_{M-1}$ and the boundary loops $\ell$ of $\textsf{G}_{\infty}^{\bar{\lambda},\varepsilon}$ (resp. $\tilde{\textsf{G}}_{\eta}^{\bar{\lambda},\varepsilon}$) satisfying $V(\ell)\subset\Lambda_{M-1}$ or $\emptyset\neq(V(\ell)\cap V')\subset\Lambda_{M-1}$.  We will distinguish three types of loops $\ell$ of $\Gamma_{\infty}''$ and construct a surjection $\varphi:\Gamma_{\infty}''\rightarrow\Gamma_{\infty}'$ as follows:
\begin{itemize}
  \item $\ell$ is also a loop in $\Gamma_{\infty}'$.  Let $\varphi(\ell)=\ell$.
  \item $\ell$ is non-trivial and has no vertices in $V'$ (thus $\ell$ is removed as a subgraph in the construction of $\textsf{G}_{\infty}'$ from $\textsf{G}_{\infty}^{\bar{\lambda},\varepsilon}$ and is contained in $F(U_{v,x})$ for some $v\in V'$ and $x\in X_v$).  Let $\varphi(\ell)$ be the trivial loop at a vertex in $V(\ell)$.
  \item $\ell$ is a natural concatenation of a loop $\ell'\in\Gamma_{\infty}'$ and a sequence of erased dangling loops of $\textsf{G}_{\infty}^{\bar{\lambda},\varepsilon}$ attached to the vertices in $V'(\ell')$.  Let $\varphi(\ell)=\ell'$.
\end{itemize}
It is easy to see that $d(\ell,\varphi(\ell))\leq\max\{\diam(F(U_{v,x})):v\in V,x\in X_v\}\leq 3s$ for every $\ell\in\Gamma_{\infty}''$, and thus
\begin{equation}\label{e31}
d_{\mathcal {L}}(\Gamma_{\infty}'',\Gamma_{\infty}')\leq 3s.
\end{equation}
Similarly, we have
\begin{equation}\label{e32}
d_{\mathcal {L}}(\tilde{\Gamma}_{\eta}'',\tilde{\Gamma}_{\eta}')\leq 3s.
\end{equation}
Combining (\ref{e31}), (\ref{e32}) and (\ref{e66}), we get
\begin{equation}\label{e33}
d_{\mathcal {L}}(\Gamma_{\infty}'',\tilde{\Gamma}_{\eta}'')\leq 12s.
\end{equation}

If $(u,v)$ is an edge of a loop $\ell$ in $\Gamma_{\infty}''$, then $u,v\in\Lambda_{M-1}$ and $(u,v)\subset F_{u,v}$,
or there exist $w\in V'\cap\Lambda_{M-1}$ and $x\in X_w$ such that $(u,v)\subset F(U_{w,x})$, which implies that $(u,v)\subset\Lambda_M$ since $\diam(F_{u,v})\leq 6s<1$ and $F(U_{w,x})\leq 3s<1$.  Therefore, $\Gamma_{\infty}^{\bar{\lambda},\varepsilon}(M-1)\subset\Gamma_{\infty}''\subset\Gamma_{\infty}^{\bar{\lambda},\varepsilon}$.  Similarly, we have $\tilde{\Gamma}_{\eta}^{\bar{\lambda},\varepsilon}(M-1)\subset\tilde{\Gamma}_{\eta}''\subset\tilde{\Gamma}_{\eta}^{\bar{\lambda},\varepsilon}$.
The lemma follows immediately from this argument and (\ref{e33}).
\end{proof}

\section{Proof of the main theorem}\label{mainproof}
Our goal in this section is to prove Theorem \ref{t1}.  Using cluster chains, we will show that for $\lambda'$ fixed and $\lambda<0$ much smaller than $\lambda'$, for $M$ large enough and $\varepsilon>0$ small enough, with high probability $\Gamma_{\eta}^{\lambda'}$ is well-approximated by the good cut-off loop ensemble $\tilde{\Gamma}_{\eta}^{\bar{\lambda},\varepsilon}$ for all small enough $\eta>0$.  This approximation result in the discrete setting, combined with Lemma \ref{l14} and Proposition \ref{p5} in the previous section, will allow us to prove the existence and uniqueness of the scaling limit of $\Gamma_{\eta}^{\lambda'}$.

We start with some notation.  Let $V$ denote the vertex set of $\tilde{\textnormal{\textsf{G}}}_{\eta}^{\bar{\lambda},\varepsilon}$.  Let $X(\tilde{\textnormal{\textsf{G}}}_{\eta}^{\bar{\lambda},\varepsilon})$ denote the set of pivotals in the Poisson point set $X$ which correspond to the edges of $\tilde{\textnormal{\textsf{G}}}_{\eta}^{\bar{\lambda},\varepsilon}$ in the construction (see Definition \ref{def7}).  For a path $\gamma$ of $\tilde{\textnormal{\textsf{G}}}_{\eta}^{\bar{\lambda},\varepsilon}$ (which is seen as a subgraph of $\tilde{\textnormal{\textsf{G}}}_{\eta}^{\bar{\lambda},\varepsilon}$), we let $X(\gamma)$ be the set of pivotals in $X(\tilde{\textnormal{\textsf{G}}}_{\eta}^{\bar{\lambda},\varepsilon})$ corresponding to the edges of $\gamma$.

In this section, we often view a cluster of $\omega_{\eta}^{\lambda}$ as a union of hexagons of $\eta\textsf{H}$.  For a site $x$ of $\eta\textsf{T}$, write $H_x$ for the $\eta$-hexagon centered at $x$.  Given a bounded set $A\subset\mathbb{C}$, recall that the filling of $A$, denoted $F(A)$, is the complement of the unbounded connected component of $\mathbb{C}\setminus A$.  Similarly, the filling $F(\mathcal{C})$ of a contractible cluster $\mathcal{C}$ of $\omega_{\eta}^{\lambda}$ on $\mathbb{T}_M^2$ is the complement of the non-contractible component of $\mathbb{T}_M^2\setminus\mathcal{C}$.

For a chain $\mathfrak{C}=(\mathcal{C}_1,y_1,\mathcal{C}_2,\ldots,y_{n-1},\mathcal{C}_n)$ of $\omega_{\eta}^{\lambda}$ with $\mathcal{C}_1,\ldots,\mathcal{C}_n$ contractible, we write
\begin{equation*}
U(\mathfrak{C}):=F(\mathcal{C}_1)\cup H_{y_1}\cup F(\mathcal{C}_2)\cup\cdots\cup H_{y_{n-1}}\cup F(\mathcal{C}_n).
\end{equation*}

Recall that the vertices of $\tilde{\textnormal{\textsf{G}}}_{\eta}^{\bar{\lambda},\varepsilon}$ are the center points of the routers related to the enhanced network $\textsf{EN}_M^{r\text{-}\textsf{meso}}(\omega_{\eta}^{\lambda},\textsf{PPP}_{\bar{\lambda}}^{\varepsilon})$, and these routers are $r$-boxes defined in Section \ref{s2}.  In this section, the vertices of $\tilde{\textnormal{\textsf{G}}}_{\eta}^{\bar{\lambda},\varepsilon}$ are simply called \textbf{routers} for convenience.  By definition, for each router $v\in V$, there exist open circuits of $\omega_{\eta}^{\lambda}$ surrounding $\hat{B}_v$ in the square annulus with inner face $\hat{B}_v$ and outer radius $\rho_M$, with a unique open cluster of $\omega_{\eta}^{\lambda}$ (resp. $\omega_{\eta}^{\lambda}\setminus X$) containing all such circuits; we say that this cluster \textbf{has the router} $v$.  Note that an open cluster of $\omega_{\eta}^{\lambda}$ may have more than one routers, but an open cluster of $\omega_{\eta}^{\lambda}\setminus X$ has at most one router when (\ref{num2}) of Proposition \ref{p6}(ii) holds.

For a path $\gamma$ of $\tilde{\textnormal{\textsf{G}}}_{\eta}^{\bar{\lambda},\varepsilon}$, let $\mathscr{C}(\gamma)$ denote the collection of all the clusters of $\omega_{\eta}^{\lambda}$ which have routers in $\gamma$, and write
\begin{equation*}
U(\gamma):=\left(\bigcup_{\mathcal{C}\in\mathscr{C}(\gamma)}F(\mathcal{C})\right)\cup\left(\bigcup_{x\in X(\gamma)}H_x\right).
\end{equation*}

A connected component of $\tilde{\textnormal{\textsf{G}}}_{\eta}^{\bar{\lambda},\varepsilon}$ will also be called a \textbf{cluster}.  A cluster of $\tilde{\textnormal{\textsf{G}}}_{\eta}^{\bar{\lambda},\varepsilon}$ is called non-trivial if it has at least two routers.

\subsection{Geometric properties of $\tilde{\textnormal{\textsf{G}}}_{\eta}^{\bar{\lambda},\varepsilon}$ and $\tilde{\Gamma}_{\eta}^{\bar{\lambda},\varepsilon}$}

We now gather some geometric facts concerning $\tilde{\textnormal{\textsf{G}}}_{\eta}^{\bar{\lambda},\varepsilon}$ and $\tilde{\Gamma}_{\eta}^{\bar{\lambda},\varepsilon}$ which will be used to compare $\tilde{\Gamma}_{\eta}^{\bar{\lambda},\varepsilon}$ with $\Gamma_{\eta}^{\lambda'}$.  For $R>0$, let $\Gamma_{\eta}^{\lambda'}(R)$ denote the collection of loops of $\Gamma_{\eta}^{\lambda'}$ which are completely contained in $\Lambda_R$.
\begin{lemma}\label{l15}
Assume that the event $\mathcal{E}$ in Proposition \ref{p5} holds.  Then $\tilde{\textnormal{\textsf{G}}}_{\eta}^{\bar{\lambda},\varepsilon}$ and $\tilde{\Gamma}_{\eta}^{\bar{\lambda},\varepsilon}$ have the following properties.
\begin{enumerate}[(i)]
\item Any loop of $\Gamma_{\eta}^{\lambda'}(M-\eta)$ with diameter larger than $s/2$ cannot intersect the non-trivial loops of $\tilde{\Gamma}_{\eta}^{\bar{\lambda},\varepsilon}$ and the non-trivial clusters of $\tilde{\textnormal{\textsf{G}}}_{\eta}^{\bar{\lambda},\varepsilon}$. \label{nointersect}
\item Let $\mathcal{C}$ be an open cluster of $\omega_{\eta}^{\lambda}$ that has routers.  Then all the routers of $\mathcal{C}$ lie in the interior of $F(\mathcal{C})$, and any two routers of $\mathcal{C}$ are connected by a path of $\tilde{\textnormal{\textsf{G}}}_{\eta}^{\bar{\lambda},\varepsilon}$ which is contained (as a curve) in the interior of $F(\mathcal{C})$, with all the routers of the path belonging to $\mathcal{C}$.  For each $x\in X(\tilde{\textnormal{\textsf{G}}}_{\eta}^{\bar{\lambda},\varepsilon})$ neighboring $\mathcal{C}$, the edge of $\tilde{\textnormal{\textsf{G}}}_{\eta}^{\bar{\lambda},\varepsilon}$ that corresponds to $x$ has a router of $\mathcal{C}$.

    If an edge of $\tilde{\textnormal{\textsf{G}}}_{\eta}^{\bar{\lambda},\varepsilon}$ intersects $F(\mathcal{C})$ with $x\in X(\tilde{\textnormal{\textsf{G}}}_{\eta}^{\bar{\lambda},\varepsilon})$ corresponding to the edge, then one of the following cases occurs:
    \begin{itemize}
      \item the edge lies strictly inside $F(\mathcal{C})$;
      \item the edge intersects $\partial F(\mathcal{C})$ at a single point in $\partial H_x$, and one router of the edge belongs to $\mathcal{C}$ while the other one belongs to some $\lambda$-cluster $\mathcal{C}'$ outside $F(\mathcal{C})$, with the edge strictly inside $F(\mathcal{C})\cup H_x\cup F(\mathcal{C}')$;
      \item the edge intersects $\partial F(\mathcal{C})$ at exactly two points in $\partial H_x$, such that the piece of the edge connecting the two points lies in $H_x$, and both routers of the edge belong to $\mathcal{C}$, with the edge strictly inside $F(\mathcal{C})\cup H_x$.
    \end{itemize}\label{topology}
\item If $\mathcal{C}$ and $\mathcal{C}'$ are two distinct open clusters of $\omega_{\eta}^{\lambda}$ which have routers, and there is a $\lambda$-closed pivotal in $X$ neighboring both $\mathcal{C}$ and $\mathcal{C}'$, then there is a $\lambda$-closed pivotal in $X(\tilde{\textnormal{\textsf{G}}}_{\eta}^{\bar{\lambda},\varepsilon})$ neighboring both of these clusters. \label{neighbor}
\item  If $\mathfrak{C}=(\mathcal{C}_1,x_1,\mathcal{C}_2,\ldots,x_{j-1},\mathcal{C}_j)$ is a chain of $\omega_{\eta}^{\lambda}$ so that $\{x_1,\ldots,x_{j-1}\}\subset X$ and both $\mathcal{C}_1$ and $\mathcal{C}_j$ have routers, then all the open clusters of $\mathfrak{C}$ have routers, and furthermore, there exist distinct $x_1',\ldots,x_{j-1}'\in X(\tilde{\textnormal{\textsf{G}}}_{\eta}^{\bar{\lambda},\varepsilon})$ such that $\mathfrak{C}'=(\mathcal{C}_1,x_1',\mathcal{C}_2,\ldots,x_{j-1}',\mathcal{C}_j)$ is a chain, and for any routers $v$ of $\mathcal{C}_1$ and $v'$ of $\mathcal{C}_j$, there is a path from $v$ to $v'$ in $\tilde{\textnormal{\textsf{G}}}_{\eta}^{\bar{\lambda},\varepsilon}$ such that every router in this path belongs to an open cluster in $\mathfrak{C}$, with the path (viewed as a curve) lying in the interior of $U(\mathfrak{C}')$.
    \label{routerpath}
\end{enumerate}
\end{lemma}

\begin{proof}
(\ref{nointersect}) By the construction of $\tilde{\textnormal{\textsf{G}}}_{\eta}^{\bar{\lambda},\varepsilon}$ in Definition \ref{def7}, all of the edges of $\tilde{\textnormal{\textsf{G}}}_{\eta}^{\bar{\lambda},\varepsilon}$ are contained in the union of $\bigcup_{v\in V}\hat{B}_v$ and the interior of the union of all the $\lambda'$-open $\eta$-hexagons on $\mathbb{T}_M^2$.  This implies that
any loop of $\Gamma_{\eta}^{\lambda'}(M-\eta)$ lying outside $\bigcup_{v\in V}\hat{B}_v$ cannot intersect the non-trivial clusters of $\tilde{\textnormal{\textsf{G}}}_{\eta}^{\bar{\lambda},\varepsilon}$.  By definition, for each $v\in V$, there exist open circuits of $\omega_{\eta}^{\lambda}$ surrounding $\hat{B}_v$ in the square annulus with inner face $\hat{B}_v$ and outer radius $\rho_M$.  Since $\mathcal{E}$ holds and thus all the open clusters of $\omega_{\eta}^{\lambda}$ on $\mathbb{T}_M^2$ have diameter less than $s/2$, any loop of $\Gamma_{\eta}^{\lambda'}(M-\eta)$ with diameter larger than $s/2$ must lie outside $\bigcup_{v\in V}\hat{B}_v$.  Then (\ref{nointersect}) follows from this and the above argument.

(\ref{topology}) Similarly to the above, each router of $\mathcal{C}$ is surrounded by an open circuit of $\mathcal{C}$, thus all the routers of $\mathcal{C}$ lie in the interior of $F(\mathcal{C})$.  For each router $v$, let $\mathcal{C}_v$ be the open cluster of $\omega_{\eta}^{\lambda}\setminus X$ that has the router $v$.  To prove the second sentence, suppose that $\mathcal{C}$ has at least two routers; the case of one router is trivial.  Let $v$ and $v'$ be any two distinct routers of $\mathcal{C}$.  Since (\ref{num2}) of Proposition \ref{p6} holds, we have $\mathcal{C}_v\neq\mathcal{C}_{v'}$, and moreover, we can find a sequence $(\mathcal{C}_1=\mathcal{C}_v,x_1,\mathcal{C}_2,x_2,\ldots,x_{j-1},\mathcal{C}_j=\mathcal{C}_{v'})$ such that the following holds:
\begin{itemize}
  \item $\mathcal{C}_1,\mathcal{C}_2,\ldots,\mathcal{C}_j$ are distinct open clusters of $\mathcal{C}\setminus X$, and for each $1\leq i\leq j$, $\mathcal{C}_i$ has a unique router $v_i\in V$, with $v_1,\ldots,v_j$ distinct.  Notice that $v_1=v$ and $v_j=v'$.
  \item $x_1,x_2,\ldots,x_{j-1}$ are distinct pivotals in $X\cap\mathcal{C}$, and for each $1\leq i\leq j-1$, $x_i$ neighbors both $\mathcal{C}_i$ and $\mathcal{C}_{i+1}$.
  \item In $\textnormal{\textsf{EN}}_M^{r\text{-}\textnormal{\textsf{meso}}}(\omega_{\eta}^{\lambda},\textnormal{\textsf{PPP}}_{\bar{\lambda}}^{\varepsilon})$, the router $v_1$ (resp. $v_j$) neighbors $x_1$ (resp. $x_{j-1}$),  and for each $2\leq i\leq j-1$, the router $v_i$ neighbors $x_{i-1}$ and $x_i$.
\end{itemize}
By the construction of $\tilde{\textnormal{\textsf{G}}}_{\eta}^{\bar{\lambda},\varepsilon}$ in Definition \ref{def7}, there is a path $v_1,\ldots,v_j$ in $\tilde{\textnormal{\textsf{G}}}_{\eta}^{\bar{\lambda},\varepsilon}$, such that for each $1\leq i\leq j-1$, there exists $x_i'\in X(\tilde{\textnormal{\textsf{G}}}_{\eta}^{\bar{\lambda},\varepsilon})\cap\mathcal{C}$ (it is possible that $x_i'=x_i$) so that the edge $(v_i,v_{i+1})$ of $\tilde{\textnormal{\textsf{G}}}_{\eta}^{\bar{\lambda},\varepsilon}$ lies strictly inside $F(\mathcal{C}_i)\cup H_{x'}\cup F(\mathcal{C}_{i+1})\subset F(\mathcal{C})$.  Then the second sentence follows.  By (\ref{num2}) of Proposition \ref{p6} again, it clear that there is an open cluster of $\mathcal{C}\setminus X$ and a pivotal $x'\in X$ with $x'\neq x$, such that this cluster has a unique router and neighbors both $x$ and $x'$, with the router neighboring $x$ in $\textnormal{\textsf{EN}}_M^{r\text{-}\textnormal{\textsf{meso}}}(\omega_{\eta}^{\lambda},\textnormal{\textsf{PPP}}_{\bar{\lambda}}^{\varepsilon})$.  Then the construction of $\tilde{\textnormal{\textsf{G}}}_{\eta}^{\bar{\lambda},\varepsilon}$ implies the third sentence.

Now let us show the rest of (\ref{topology}).  Let $(u,v)$ denote the edge in the statement.  If $x$ lies outside $F(\mathcal{C})$, then $x$ must neighbors $F(\mathcal{C})$.  Indeed, otherwise, because $(u,v)\cap F(\mathcal{C})\neq\emptyset$, one of $\hat{B}_u$ and $\hat{B}_v$, say $\hat{B}_u$, would intersect $F(\mathcal{C})$ without lying in $F(\mathcal{C})$, and thus $F(\mathcal{C})$ would be surrounded by an open circuit of $\mathcal{C}_u$ in the square annulus with inner face $\hat{B}_u$ and outer radius $\rho_M$, which together with the fact that $\mathcal{C}$ has routers would imply that there exists $x'\in X$ lying in or neighboring $\mathcal{C}$ so that $x'$ is surrounded by the above open circuit, contradicting the fact that $u$ is at distance at least $2\rho$ from $x'$.  Thus, at least one of $\mathcal{C}_u$ and $\mathcal{C}_v$ belongs to $\mathcal{C}$, since otherwise $x$ would neighbor three routers in $\textnormal{\textsf{EN}}_M^{r\text{-}\textnormal{\textsf{meso}}}(\omega_{\eta}^{\lambda},\textnormal{\textsf{PPP}}_{\bar{\lambda}}^{\varepsilon})$ with one of them belonging to $\mathcal{C}$, a contradiction.  By the above argument, there are three cases:
\begin{itemize}
  \item $x$ lies strictly inside $F(\mathcal{C})$.
  \item $x$ neighbors $F(\mathcal{C})$, and exactly one of $\mathcal{C}_u$ and $\mathcal{C}_v$ belongs to $\mathcal{C}$.
  \item $x$ neighbors $F(\mathcal{C})$, and both $\mathcal{C}_u$ and $\mathcal{C}_v$ belong to $\mathcal{C}$.
\end{itemize}
The above cases correspond to the three cases in (\ref{topology}), respectively.  In the first case, it is clear that $F(\mathcal{C}_u)\cup H_x\cup F(\mathcal{C}_v)\subset F(\mathcal{C})$.  Then, by the construction of $\tilde{\textnormal{\textsf{G}}}_{\eta}^{\bar{\lambda},\varepsilon}$, $(u,v)$ lies strictly inside $F(\mathcal{C}_u)\cup H_x\cup F(\mathcal{C}_v)$, and thus it lies strictly inside $F(\mathcal{C})$.  Similarly, the other two cases in (\ref{topology}) follow from the construction of $\tilde{\textnormal{\textsf{G}}}_{\eta}^{\bar{\lambda},\varepsilon}$.

(\ref{neighbor}) Let $x$ be a $\lambda$-closed pivotal in $X$ neighboring both $\mathcal{C}$ and $\mathcal{C}'$.  Since (\ref{num2}) in Proposition \ref{p6} holds, it follows that there is an open cluster $\mathcal{C}_v$ of $\mathcal{C}\setminus X$ (resp. $\mathcal{C}_{v'}'$ of $\mathcal{C}'\setminus X$), such that it has a unique router $v$ (resp. $v'$) and neighbors $x$, with the router neighboring $x$ in $\textnormal{\textsf{EN}}_M^{r\text{-}\textnormal{\textsf{meso}}}(\omega_{\eta}^{\lambda},\textnormal{\textsf{PPP}}_{\bar{\lambda}}^{\varepsilon})$.  By our construction of $\tilde{\textnormal{\textsf{G}}}_{\eta}^{\bar{\lambda},\varepsilon}$, there is an edge $(v,v')$ in $\tilde{\textnormal{\textsf{G}}}_{\eta}^{\bar{\lambda},\varepsilon}$ crossing through
an $\eta$-hexagon $H_{x'}$ centered at a closed pivotal $x'\in X(\tilde{\textnormal{\textsf{G}}}_{\eta}^{\bar{\lambda},\varepsilon})$ (it is possible that $x'=x$) neighboring both $\mathcal{C}_v$ and $\mathcal{C}_{v'}'$.  It is clear that $x'$ neighbors both $\mathcal{C}$ and $\mathcal{C}'$.

(\ref{routerpath})  If $\mathfrak{C}$ has only one open cluster of $\omega_{\eta}^{\lambda}$, then (\ref{routerpath}) follows from the second sentence of (\ref{topology}).  Now assume that $j\geq 2$.  Since (\ref{num4}) of Proposition \ref{p6} holds, it is easy to see that for each $2\leq i\leq j-1$, $\mathcal{C}_i$ is the unique open cluster of $\omega_{\eta}^{\lambda}$ on $\mathbb{T}_M^2$ that neighbors both $x_{i-1}$ and $x_i$.  Then, (\ref{num2}) of Proposition \ref{p6} implies that for each $2\leq i\leq j-1$, $\mathcal{C}_i$ has routers.  Thus, all the open clusters of $\mathfrak{C}$ have routers.  Because of (\ref{neighbor}), it follows that for each $1\leq i\leq j-1$, there exists a $\lambda$-closed pivotal $x_i'$ in $X(\tilde{\textnormal{\textsf{G}}}_{\eta}^{\bar{\lambda},\varepsilon})$ neighboring both $\mathcal{C}_i$ and $\mathcal{C}_{i+1}$.  Moreover, $x_1',\ldots,x_{j-1}'$ must be distinct, since otherwise there would exist some $x_i'$ neighboring at least three open clusters of $\omega_{\eta}^{\lambda}$ that have routers, and thus by the third sentence of (\ref{topology}), $x_i'$ would have three router-neighbors in $\textnormal{\textsf{EN}}_M^{r\text{-}\textnormal{\textsf{meso}}}(\omega_{\eta}^{\lambda},\textnormal{\textsf{PPP}}_{\bar{\lambda}}^{\varepsilon})$, contradicting (\ref{num4}) of Proposition \ref{p6}.  Therefore, $\mathfrak{C}'=(\mathcal{C}_1,x_1',\mathcal{C}_2,\ldots,x_{j-1}',\mathcal{C}_j)$ is indeed a chain.

By the third sentence of (\ref{topology}) again, for each $1\leq i\leq j-1$, the edge of $\tilde{\textnormal{\textsf{G}}}_{\eta}^{\bar{\lambda},\varepsilon}$ which corresponds to $x_i'$ has a router of $\mathcal{C}_i$ and a router of $\mathcal{C}_{i+1}$, denoted by $v_i^+$ and $v_{i+1}^-$, respectively.  By the second paragraph of (\ref{topology}), each edge $(v_i^+,v_{i+1}^-)$ of $\tilde{\textnormal{\textsf{G}}}_{\eta}^{\bar{\lambda},\varepsilon}$ lies strictly inside $F(\mathcal{C}_i)\cup H_{x_i'}\cup F(\mathcal{C}_{i+1})$.  Write $v_1^-=v$ and $v_j^+=v'$.  By the second sentence of (\ref{topology}), for each $1\leq i\leq j$, there exists a path $\gamma_i$ from $v_i^-$ to $v_i^+$ in $\tilde{\textnormal{\textsf{G}}}_{\eta}^{\bar{\lambda},\varepsilon}$, such that $\gamma_i$ is contained in the interior of $F(\mathcal{C}_i)$, with all the routers of $\gamma_i$ belonging to $\mathcal{C}_i$.  By concatenating $\gamma_1,(v_1^+,v_2^-),\gamma_2,\ldots,(v_{j-1}^+,v_j^-),\gamma_j$, we get the desired path.
\end{proof}

\subsection{Approximation of $\Gamma_{\eta}^{\lambda'}$ by $\tilde{\Gamma}_{\eta}^{\bar{\lambda},\varepsilon}$ on $\Lambda_M$}

We start with some notation.  Let $M\geq 8$ and $0<\delta<s/2<1/20$.  Let $\mathcal{A}(\mathrm{A},\delta)$ and $\mathscr{A}_{M,s}$ be the event and the collection of annuli defined above Lemma \ref{l9}, respectively.  Denote by $\Lambda_s^{\circ}$ the interior of $\Lambda_s$.  Write
\begin{align*}
&\mathcal{F}_1:=\{\mbox{each open cluster of $\omega_{\eta}^{\lambda}$ on $\mathbb{T}_M^2$ with diameter at least $\delta$ has routers in $\tilde{\textnormal{\textsf{G}}}_{\eta}^{\bar{\lambda},\varepsilon}$}\},\\
&\mathcal{F}_2:=\{\mbox{for $\omega_{\eta}^{\lambda}$ on $\mathbb{T}_M^2$, $\mathcal{A}(\mathrm{A},\delta)$ occurs for all $A\in\mathscr{A}_{M,s}$}\},\quad\mbox{and}\\
&\mathcal{F}_3:=\left\{
\begin{aligned}
&\mbox{for any square $B\subset\Lambda_{M-s}$, if there are two $\lambda$-clusters $\mathcal{C},\mathcal{C}'$ of diameter at least $\delta$}\\
&\mbox{completely contained in $B$, and there is an open path of $\omega_{\eta}^{\lambda'}$ connecting $\mathcal{C}$ and $\mathcal{C}'$}\\
&\mbox{in $B$, then there is a chain $\mathfrak{C}=(\mathcal{C},x_1,\ldots,x_m,\mathcal{C}')$ of $\omega_{\eta}^{\lambda}$ such that $U(\mathfrak{C})\subset B+\Lambda_s^{\circ}$}\\
&\mbox{and $\{x_1,\ldots,x_m\}\subset X(\tilde{\textnormal{\textsf{G}}}_{\eta}^{\bar{\lambda},\varepsilon})$, and moreover, for any routers $v$ of $\mathcal{C}$ and $v'$ of $\mathcal{C}'$,}\\
&\mbox{there is a path from $v$ to $v'$ in $\tilde{\textnormal{\textsf{G}}}_{\eta}^{\bar{\lambda},\varepsilon}$ such that every router in this path belongs to}\\
&\mbox{an open cluster in $\mathfrak{C}$, with this path (viewed as a curve) lying in the interior of $U(\mathfrak{C})$
}
\end{aligned}
\right\}.
\end{align*}

Our first lemma says that with high probability some ``good" event $\mathcal{F}$ occurs; we will compare $\Gamma_{\eta}^{\lambda'}$ to $\tilde{\Gamma}_{\eta}^{\bar{\lambda},\varepsilon}$ conditioned on $\mathcal{F}$.
\begin{lemma}\label{l19}
Fix any $\lambda'\in\mathbb{R},M\geq 8$ and $\alpha,s\in(0,1/10)$.  Let $\mathcal{E}$ be the event in Proposition \ref{p5}.  If $\lambda<\lambda'$ is very negative, $\delta\in(0,s/2)$ is very small, and $\varepsilon\in(0,\delta)$ is very small, then for all $\eta>0$ small enough, for the event $\mathcal{F}:=\mathcal{E}\cap\mathcal{F}_1\cap\mathcal{F}_2\cap\mathcal{F}_3$ we have
\begin{equation*}
\mathbf{P}[\mathcal{F}]>1-2\alpha.
\end{equation*}
\end{lemma}
\begin{proof}
Write
\begin{equation*}
\mathcal{F}_3':=\left\{
\begin{aligned}
&\mbox{for any square $B\subset\Lambda_{M-s}$, if there are two $\lambda$-clusters $\mathcal{C},\mathcal{C}'$ of diameter}\\
&\mbox{at least $\delta$ completely contained in $B$, and there is an open path of $\omega_{\eta}^{\lambda'}$}\\
&\mbox{connecting $\mathcal{C}$ and $\mathcal{C}'$ in $B$, then there is a chain $\mathfrak{C}=(\mathcal{C},y_1,\ldots,y_m,\mathcal{C}')$}\\
&\mbox{of $\omega_{\eta}^{\lambda}$ such that $U(\mathfrak{C})\subset B+\Lambda_{s/2}$ and $\{y_1,\ldots,y_m\}\subset X$}
\end{aligned}
\right\}.
\end{equation*}
Notice that if $\mathcal{E}$ holds and if an open cluster of $\omega_{\eta}^{\lambda}$ on $\mathbb{T}_M^2$ with diameter in $[\delta,s/2)$ contains at least two sites in $X$, then this cluster has routers in $\tilde{\textnormal{\textsf{G}}}_{\eta}^{\bar{\lambda},\varepsilon}$.  Then, combining Proposition \ref{p5}, Lemmas \ref{l8}, \ref{l9} and \ref{l3}, it follows that if $\lambda<\lambda'$ is very negative, $\delta\in(0,s/2)$ is very small, and $\varepsilon\in(0,\delta)$ is very small, then for all $\eta>0$ small enough,
\begin{equation*}
\mathbf{P}[\mathcal{E}\cap\mathcal{F}_1\cap\mathcal{F}_2\cap\mathcal{F}_3']>1-2\alpha.
\end{equation*}
This together with Lemma \ref{l15}(\ref{routerpath}) gives the desired result.
\end{proof}

Let $\mathscr{B}^s:=\{[j,j+1)s\times[k,k+1)s:(j,k)\in \mathbb{Z}^2\}$ be the set of lattice $s$-boxes.  For $R>0$ and $z\in\mathbb{C}$, let $B_R(z):=z+[-R/2,R/2)\times[-R/2,R/2)$ denote the $R$-box centered at $z$.  In order to compare $\Gamma_{\eta}^{\lambda'}$ to $\tilde{\Gamma}_{\eta}^{\bar{\lambda},\varepsilon}$, we use the following method to divide a large loop of $\Gamma_{\eta}^{\lambda'}$ into segments, and will compare them to the corresponding segments of a loop in $\tilde{\Gamma}_{\eta}^{\bar{\lambda},\varepsilon}$ in the next lemma.

\begin{figure}
\centering
\includegraphics[height=0.35\textwidth]{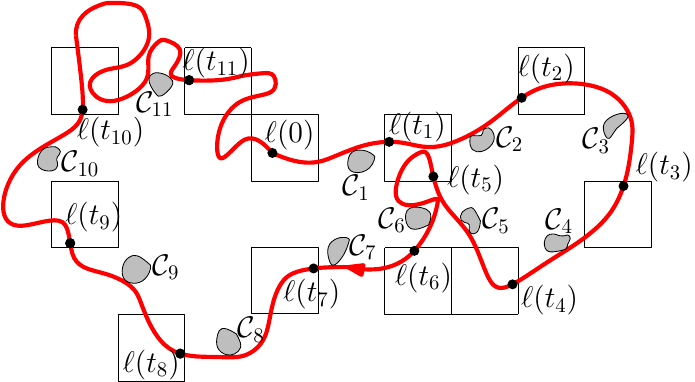}
\caption{A clockwise loop $\ell$ of $\Gamma_{\eta}^{\lambda'}$ with its associated $s$-boxes in Definition \ref{def3} and $\lambda$-clusters in Lemma \ref{l2}.  The black dots represent $\ell(t_0=0),\ell(t_1),\ldots,\ell(t_N)$, the squares represent $B_s(z_0),B_s(z_1),\ldots,B_s(z_N)$, and the gray regions represent $\mathcal{C}_1,\ldots,\mathcal{C}_N$, where $N=11$ and $B_s(z_1)=B_s(z_5)$.
}\label{boxes}
\end{figure}

\begin{definition}[Dividing a large loop of $\Gamma_{\eta}^{\lambda'}$ into segments by a sequence of $s$-boxes]\label{def3}
Given a loop $\ell$ of $\Gamma_{\eta}^{\lambda'}$ with $\diam(\ell)\geq 10s$, we associate to $\ell$ a sequence of $s$-boxes in $\mathscr{B}^s$ and divide $\ell$ into segments in the following way.  Write $t_0:=0$, and let $B_s(z_0)$ be the box of $\mathscr{B}^s$
containing $\ell(0)$.  We start at the point $\ell(0)$ and move along $\ell$ in its orientation.  Let $\ell(t_1)$ be the first vertex (of $\eta\textsf{H}$) outside $B_{3s}(z_0)$, and let $B_s(z_1)$ be the box of $\mathscr{B}^s$ containing $\ell(t_1)$.  Let $\ell(t_2)$ be the first vertex outside $B_{3s}(z_1)$ after time $t_1$, and let $B_s(z_2)$ be the box of $\mathscr{B}^s$ containing $\ell(t_2)$.  The process stops after $N=N(\ell)$ steps, and we have $\ell([t_N,1])\subset B_{3s}(z_N)$ and $B_s(z_0)\subset B_{3s}(z_N)$.  Then we obtain a sequence of boxes $B_s(z_0),B_s(z_1),\ldots,B_s(z_N)$ in $\mathscr{B}^s$ and a sequence of segments $\ell([0,t_1]),\ldots,\ell([t_{N-1},t_N]),\ell([t_N,1])$ of $\ell$.  See Figure \ref{boxes} for an example.
\end{definition}

For $R>0$, let $\Gamma_{\eta}^{\lambda'}(R)$ (resp. $\tilde{\Gamma}_{\eta}^{\bar{\lambda},\varepsilon}(R)$) be the collection of loops in $\Gamma_{\eta}^{\lambda'}$ (resp. $\tilde{\Gamma}_{\eta}^{\bar{\lambda},\varepsilon}$) which are entirely contained in $\Lambda_R$.

\begin{lemma}[Each large loop in $\Gamma_{\eta}^{\lambda'}(3M/4)$ being close to a loop in $\tilde{\Gamma}_{\eta}^{\bar{\lambda},\varepsilon}$]\label{l2}
Suppose that the event $\mathcal{F}$ in Lemma \ref{l19} holds.  Then for each loop $\ell\in\Gamma_{\eta}^{\lambda'}(3M/4)$ of diameter at least $10s$, there is a sequence of $\lambda$-clusters $\mathcal{C}_1,\ldots,\mathcal{C}_N$ satisfying the following two conditions (see Figures \ref{boxes} and \ref{case1chain}; the notation in Definition \ref{def3} was used):

\begin{enumerate}[(i)]
\item For each $1\leq i\leq N$, $\mathcal{C}_i$ touches the segment $\ell([t_{i-1},t_i])$, with $\mathcal{C}_i\subset B_{4s}(z_{i-1})$ and $\diam(\mathcal{C}_i)\geq\delta$.\label{con1}
\item Write $\mathcal{C}_{N+1}:=\mathcal{C}_1$.  Let
\begin{equation*}
B_i:=
\left\{
\begin{aligned}
&B_{8s}((z_{i-1}+z_i)/2)\qquad\mbox{ if $1\leq i\leq N-1$,}\\
&B_{9s}(z_{N-1}/4+3z_N/4)\quad\mbox{if $i=N$.}
\end{aligned}
\right.
\end{equation*}
For each $1\leq i\leq N$, there exists a chain $\mathfrak{C}_i=(\mathcal{C}_{i,1}=\mathcal{C}_i,x_{i,1},\mathcal{C}_{i,2},x_{i,2},\ldots,x_{i,n_i-1},\mathcal{C}_{i,n_i}=\mathcal{C}_{i+1})$ of $\omega_{\eta}^{\lambda}$, such that $x_{i,1},\ldots,x_{i,n_i-1}\in X(\tilde{\textnormal{\textsf{G}}}_{\eta}^{\bar{\lambda},\varepsilon})$ and $U(\mathfrak{C}_i)\subset B_i$.  Moreover, for any router $v_i$ of $\mathcal{C}_i$, there is a path from $v_i$ to $v_{i+1}$ in $\tilde{\textnormal{\textsf{G}}}_{\eta}^{\bar{\lambda},\varepsilon}$, such that all the routers in this path belong to the open clusters in $\mathfrak{C}_i$, with the path (viewed as a curve) lying in the interior of $U(\mathfrak{C}_i)$.
\label{con2}
\end{enumerate}

Let $\tilde{\mathcal{C}}=\tilde{\mathcal{C}}(\ell)$ denote the cluster of $\tilde{\textnormal{\textsf{G}}}_{\eta}^{\bar{\lambda},\varepsilon}$ which contains the routers of $\mathcal{C}_1,\ldots,\mathcal{C}_N$.  If $\ell$ is clockwise, then $\ell$ surrounds $\tilde{\mathcal{C}}$ and we let $\tilde{\ell}$ be the clockwise boundary loop of $\tilde{\mathcal{C}}$; if $\ell$ is counterclockwise, then there is a unique counterclockwise boundary loop $\tilde{\ell}$ of $\tilde{\mathcal{C}}$ surrounding $\ell$ in $\Lambda_M$.  Furthermore, we have
\begin{equation}\label{e17}
d(\ell,\tilde{\ell})\leq 9s.
\end{equation}
\end{lemma}

\begin{proof}
Since $\mathcal{F}_2$ occurs, it follows that if a percolation interface of $\omega_{\eta}^{\lambda'}$ radially crosses an annulus $A\in\mathscr{A}_{M,s}$, then there is an open cluster of $\omega_{\eta}^{\lambda}\cap A$ of diameter at least $\delta$ touching this interface.  It is clear that for each $1\leq i\leq N$, $\ell([t_{i-1},t_i])$ radially crosses $B_{3s}(z_{i-1})\setminus B_s(z_{i-1})$.  The above argument and the fact that all open clusters of $\omega_{\eta}^{\lambda}$ on $\mathbb{T}_M^2$ have diameter less than $s/2$ (since $\mathcal{E}$ holds) imply that there is a sequence of $\lambda$-clusters $\mathcal{C}_1,\ldots,\mathcal{C}_N$ satisfying condition (\ref{con1}).  For each $1\leq i\leq N-1$, by using the $\lambda'$-open hexagons touching the interface $\ell([t_{i-1},t_{i+1}])$ one can find a $\lambda'$-open path connecting $\mathcal{C}_i$ and $\mathcal{C}_{i+1}$ in $B_{4s}(z_{i-1})\cup B_{4s}(z_i)\subset B_{6s}((z_{i-1}+z_i)/2)$.  Moreover, one can also find a $\lambda'$-open path connecting
$\mathcal{C}_N$ and $\mathcal{C}_1$ in $B_{4s}(z_{N-1})\cup B_{4s}(z_N)\cup B_{4s}(z_0)\subset B_{7s}(z_{N-1}/4+3z_N/4)$.  Then, the occurrence of $\mathcal{F}_1\cap\mathcal{F}_3$ implies that the sequence $\mathcal{C}_1,\ldots,\mathcal{C}_N$ also satisfies condition (\ref{con2}).

It remains to prove the last paragraph of the lemma.  For $1\leq i\leq N$, let $\ell(t_i^*)$ be some touching point of $\mathcal{C}_i$ with $\ell([t_{i-1},t_i])$, and let $\gamma_i$ be the oriented piece of $\ell$ from time $t_i^*$ to time $t_{i+1}^*$, with $t_{N+1}^*:=t_1^*$.  By Lemma \ref{l15}(\ref{nointersect}) and $\diam(\ell)\geq 10s$, we have $\ell\cap\tilde{\mathcal{C}}=\emptyset$.  Now we distinguish two cases:

\textbf{Case 1: $\ell$ is clockwise.}  It is clear that $F(\mathcal{C}_i)$ is surrounded by $\ell$ for each $1\leq i\leq N$.  By Lemma \ref{l15}(\ref{topology}),
all the routers of $\mathcal{C}_i$ is surrounded by $\partial F(\mathcal{C}_i)$.  Thus, $\ell$ surrounds all the routers of $\mathcal{C}_i$, which lie in $\tilde{\mathcal{C}}$.  By the Jordan curve theorem and $\ell\cap\tilde{\mathcal{C}}=\emptyset$, we deduce that $\ell$ surrounds $\tilde{\mathcal{C}}$, and thus $\ell$ surrounds the clockwise boundary loop $\tilde{\ell}$ of $\tilde{\mathcal{C}}$.

\begin{figure}
\centering
\includegraphics[height=0.25\textwidth]{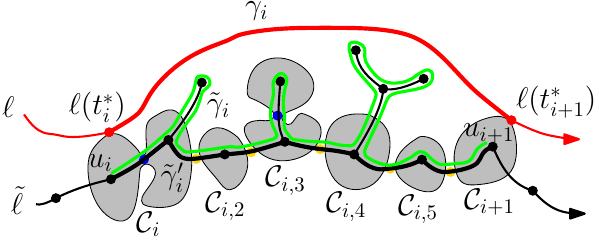}
\caption{A piece of a clockwise loop $\ell$ of $\Gamma_{\eta}^{\lambda'}$, a piece of $\tilde{\ell}$ and the chain $\mathfrak{C}_i$ in Lemma \ref{l2}.  Black dots and black lines represent the vertices and edges of $\tilde{\textnormal{\textsf{G}}}_{\eta}^{\bar{\lambda},\varepsilon}$ in $\tilde{\ell}$, respectively.  Red (resp. green, black) heavy line represents $\gamma_i$ (resp., $\tilde{\gamma}_i,\tilde{\gamma}_i'$) in the proof of Lemma \ref{l2}.  Blue (resp. yellow) dots represent $\lambda$-open (resp. $\lambda$-closed) pivotals in $X(\tilde{\textnormal{\textsf{G}}} _{\eta}^{\bar{\lambda},\varepsilon})$.
}\label{case1chain}
\end{figure}

The facts that $\tilde{\mathcal{C}}$ intersects all $\mathcal{C}_i$'s, $\ell$ surrounds $\tilde{\ell}$, and each $\mathcal{C}_i$ touches $\ell$, imply that $\tilde{\ell}$ intersects each $\mathcal{C}_i$.  Let $\mathcal{R}_i$ be the component of $F(\mathcal{C}_i)\setminus\tilde{\ell}$ which contains $\ell(t_i^*)$.  Then, there exists a router $u_i$ of $\tilde{\ell}$ on $\partial\mathcal{R}_i$ (see Figure \ref{case1chain}), since otherwise there would exist a segment of an edge of $\tilde{\ell}$ lying in $F(\mathcal{C}_i)$ and connecting two distinct points in $\partial F(\mathcal{C}_i)$, contradicting the second paragraph of Lemma \ref{l15}(\ref{topology}).  Let $\hat{u}_i$ be a prime end on $\partial\mathcal{R}_i$ that corresponds to $u_i$, and denote by $\tilde{\gamma}_i$ the oriented segment of $\tilde{\ell}$ from $\hat{u}_i$ to $\hat{u}_{i+1}$.

By condition (\ref{con2}), we can find a self-avoiding oriented path $\tilde{\gamma}_i'$ from $u_i$ to $u_{i+1}$ in $\tilde{\mathcal{C}}$, such that all the routers in $\tilde{\gamma}_i'$ belong to the open clusters in $\mathfrak{C}_i$, with $\tilde{\gamma}_i'\subset U(\mathfrak{C}_i)$.  Then, it is easy to see that
\begin{equation}\label{e29}
\tilde{\gamma}_i\subset F(\mathcal{C}_i\cup\mathcal{C}_{i+1}\cup\gamma_i\cup\tilde{\gamma}_i')\subset F(U(\mathfrak{C}_i)\cup\gamma_i)\subset B_i.
\end{equation}
See Figure \ref{case1chain} for an illustration.  We can construct an increasing reparametrization of $\tilde{\ell}$ on $[0,1]$ so that $\tilde{\ell}(t_i^*)=\hat{u}_i$ for each $1\leq i\leq N$, and thus $\tilde{\gamma}_i$ is the oriented segment of $\tilde{\ell}$ from time $t_i^*$ to time $t_{i+1}^*$.  Then (\ref{e17}) follows from (\ref{e29}) immediately.

\textbf{Case 2: $\ell$ is counterclockwise.}  Let $\tilde{\gamma}_i'$ be the oriented path from $v_i$ to $v_{i+1}$ in condition (\ref{con2}).  We concatenate the paths $\tilde{\gamma}_1',\ldots,\tilde{\gamma}_N'$ to obtain an oriented circuit $\tilde{\ell}'$ of $\tilde{\mathcal{C}}$.  Since $M\geq 8, s<1/10, \ell\subset \Lambda_{3M/4}$ and $\gamma_i,\tilde{\gamma}_i'\subset B_i$ for each $i$, we have $\tilde{\ell}'\subset\Lambda_M$.  It is clear that $\tilde{\ell}'$ lies outside $\ell$.  To show that $\tilde{\ell}'$ surrounds $\ell$, it remains to show that for a fixed point $o$ surrounded by $\ell$, the winding number $\wind(\tilde{\ell}',o)$ of $\tilde{\ell}'$ around $o$ is equal to 1.  Let $(\ell(t_i^*)v_i)$ be a curve from $\ell(t_i^*)$ to $v_i$ in $F(\mathcal{C}_i)$ and let $(v_i\ell(t_i^*))$ be its reversal.  Let $\tilde{\gamma}_i''$ be the concatenation of the curves $(\ell(t_i^*)v_i),\tilde{\gamma}_i'$ and $(v_{i+1}\ell(t_{i+1}^*))$. Concatenating the curves $\tilde{\gamma}_1'',\ldots,\tilde{\gamma}_N''$, we obtain a loop $\tilde{\ell}''$.  It is obvious that
\begin{equation}\label{e19}
\wind(\tilde{\ell}'',o)=\wind(\tilde{\ell}',o).
\end{equation}
Since $\gamma_i,\tilde{\gamma}_i''\subset B_i,\diam(B_i)\leq 9s,\diam(\ell)\geq 10s$ and $\tilde{\gamma}_i''$ does not intersect the domain surrounded by $\ell$, it follows that $\wind(\tilde{\gamma}_i'',o)=\wind(\gamma_i,o)$.  Therefore,
\begin{equation}\label{e18}
\wind(\tilde{\ell}'',o)=\sum_{i=1}^N\wind(\tilde{\gamma}_i'',o)=\sum_{i=1}^N\wind(\gamma_i,o)=\wind(\ell,o)=1.
\end{equation}
Combining (\ref{e19}) and (\ref{e18}), we get that $\tilde{\ell}'$ surrounds $\ell$.  Then, since $\tilde{\ell}'$ is a circuit of $\tilde{\mathcal{C}}$ and $\tilde{\ell}'\subset\Lambda_M$, there is a unique counterclockwise boundary loop $\tilde{\ell}$ of $\tilde{\mathcal{C}}$ surrounding $\ell$ in $\Lambda_M$.

The facts that $\tilde{\mathcal{C}}$ intersects all $\mathcal{C}_i$'s, $\tilde{\ell}$ surrounds $\ell$, and each $\mathcal{C}_i$ touches $\ell$, imply that $\tilde{\ell}$ intersects each $\mathcal{C}_i$.  The rest of the argument is analogous to that in Case 1.
\end{proof}

The above lemma shows that for each large loop in $\Gamma_{\eta}^{\lambda'}(3M/4)$, one can find a loop in $\tilde{\Gamma}_{\eta}^{\bar{\lambda},\varepsilon}$ close to it; the following is a result for the other direction.

\begin{lemma}[Each large loop in $\tilde{\Gamma}_{\eta}^{\bar{\lambda},\varepsilon}(5M/8)$ being close to a loop in $\Gamma_{\eta}^{\lambda'}$]\label{l18}
Suppose that the event $\mathcal{F}$ in Lemma \ref{l19} holds.  Then for each loop $\tilde{\ell}\in\tilde{\Gamma}_{\eta}^{\bar{\lambda},\varepsilon}(5M/8)$ of
diameter at least $137s$ and the cluster $\mathcal{C}(\tilde{\ell})$ of $\omega_{\eta}^{\lambda'}$ having the routers in $\tilde{\ell}$, the following holds.  If $\tilde{\ell}$ is clockwise, then $\mathcal{C}(\tilde{\ell})$ is contained in $\Lambda_{3M/4}$, and its clockwise boundary loop $\ell$ surrounds $\tilde{\ell}$ and satisfies
\begin{equation}\label{e30}
d(\ell,\tilde{\ell})\leq 9s;
\end{equation}
if $\tilde{\ell}$ is counterclockwise, then there exists a counterclockwise boundary loop $\ell$ of $\mathcal{C}(\tilde{\ell})$ that is surrounded by $\tilde{\ell}$, such that (\ref{e30}) holds.
\end{lemma}
\begin{proof}
We first define the event
\begin{equation*}
\mathcal{G}:=\left\{
\begin{aligned}
&\mbox{each path in $\tilde{\textnormal{\textsf{G}}}_{\eta}^{\bar{\lambda},\varepsilon}$ of diameter at least $3s$ has a router such that}\\
&\mbox{the diameter of the $\lambda$-cluster that has this router is at least $\delta$}
\end{aligned}
\right\},
\end{equation*}
and explain that $\mathcal{G}\supset\mathcal{F}$.  Let $\tilde{\gamma}$ be a path in $\tilde{\textnormal{\textsf{G}}}_{\eta}^{\bar{\lambda},\varepsilon}$ of diameter at least $3s$.  It is clear that $\diam(U(\tilde{\gamma}))\geq\diam(\tilde{\gamma})\geq 3s$.  Since $\mathcal{F}_2$ occurs, there is some $A\in\mathscr{A}_{M,s}$ such that $U(\tilde{\gamma})$ radially crosses a necklace of $\omega_{\eta}^{\lambda}\cap A$, with the necklace surrounding the inner face of $A$ and with all the restricted $\lambda$-clusters of the necklace having diameter at least $\delta$.  So, there is a restricted $\lambda$-cluster of the necklace lying in a cluster $\mathcal{C}_0$ of $\mathscr{C}(\tilde{\gamma})$, since otherwise by the occurrence of $\mathcal{F}_1$, there would exist a pivotal $x\in X(\tilde{\textnormal{\textsf{G}}}_{\eta}^{\bar{\lambda},\varepsilon})$ neighboring at least three $\lambda$-clusters (one contains a restricted cluster of the necklace and the other two are from $\mathscr{C}(\tilde{\gamma})$, see Figure \ref{necklacepath}) such that all these clusters have routers, and thus by the third sentence of Lemma \ref{l15}(\ref{topology}) the edge of $\tilde{\textnormal{\textsf{G}}}_{\eta}^{\bar{\lambda},\varepsilon}$ that corresponds to $x$ has three routers, a contradiction.  Obviously, $\diam(\mathcal{C}_0)\geq\delta$, which gives $\mathcal{G}\supset\mathcal{F}$.

\begin{figure}
\centering
\includegraphics[height=0.35\textwidth]{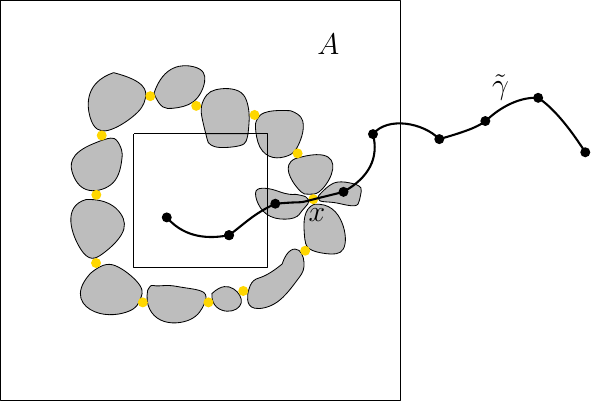}
\caption{A path $\tilde{\gamma}$ of $\tilde{\textnormal{\textsf{G}}}_{\eta}^{\bar{\lambda},\varepsilon}$ crossing a necklace of $\omega_{\eta}^{\lambda}\cap A$, with each restricted $\lambda$-cluster of the necklace having diameter at least $\delta$ and not lying in any cluster of $\mathscr{C}(\tilde{\gamma})$.  Such a scenario would contradict the occurrence of $\mathcal{F}$.
}\label{necklacepath}
\end{figure}

Since $\diam(\tilde{\ell})>3s$ and $\mathcal{G}$ holds, there is a router $u\in\tilde{\ell}$ and a cluster $\mathcal{C}_u\in\mathscr{C}(\tilde{\ell})$ of diameter at least $\delta$ having the router $u$.  Let $\tilde{\mathcal{C}}(\tilde{\ell})$ denote the cluster of $\tilde{\textnormal{\textsf{G}}}_{\eta}^{\bar{\lambda},\varepsilon}$ which contains the routers of $\tilde{\ell}$.  Now we distinguish the following two cases.

\textbf{Case 1: $\tilde{\ell}$ is clockwise.}
First, let us show that $\mathcal{C}(\tilde{\ell})\subset\Lambda_{3M/4}$.  Assume, by way of contradiction, that this is not the case.  Since $\tilde{\ell}\subset\Lambda_{5M/8}$ and all the open clusters of $\omega_{\eta}^{\lambda}$ on $\mathbb{T}_M^2$ have diameter less than $s/2$, we get that all the clusters of $\mathscr{C}(\tilde{\ell})$ are contained in $\Lambda_{5M/8+s/2}$.  Recall that $M\geq 8$ and $s<1/10$.  Then there would exist a $\lambda'$-open path $\gamma$ in $\mathcal{C}(\tilde{\ell})\cap\Lambda_{3M/4-2s}$, such that $\gamma$ touches some cluster in $\mathscr{C}(\tilde{\ell})$ and radially crosses some annulus of $\mathscr{A}_{M,s}$ contained in $A(5M/8+s,3M/4-2s)$.  By considering a necklace in this annulus and using a similar argument as that in the first paragraph of the proof, it follows that $\gamma$ would intersect or touch a $\lambda$-cluster $\mathcal{C}_v$ outside $\tilde{\ell}$ in $\Lambda_{3M/4-s}$, such that $\mathcal{C}_v$ has diameter at least $\delta$ and has a router $v$.  So, $\mathcal{C}_u$ and $\mathcal{C}_v$ are connected by a $\lambda'$-open path in $\Lambda_{3M/4-s}$.  Since $\mathcal{F}_3$ holds, the routers $u$ and $v$ are connected in $\tilde{\textnormal{\textsf{G}}}_{\eta}^{\bar{\lambda},\varepsilon}$,
and thus a router outside $\tilde{\ell}$ is connected to $\tilde{\ell}$ in $\tilde{\textnormal{\textsf{G}}}_{\eta}^{\bar{\lambda},\varepsilon}$, contradicting the fact that $\tilde{\ell}$ is the clockwise boundary loop of a cluster of $\tilde{\textnormal{\textsf{G}}}_{\eta}^{\bar{\lambda},\varepsilon}$ in $\Lambda_{5M/8}$.

Therefore, $\mathcal{C}(\tilde{\ell})\subset\Lambda_{3M/4}$.  Since $\tilde{\ell}$ lies in the interior of $U(\tilde{\ell})$ and $U(\tilde{\ell})\subset F(\mathcal{C}(\tilde{\ell}))$, we get that $\ell$ surrounds $\tilde{\ell}$, and thus $\diam(\ell)\geq\diam(\tilde{\ell})\geq 137s$.  Let $\mathcal{C}_1,\ldots,\mathcal{C}_N$ and $\tilde{\mathcal{C}}(\ell)$ be as in Lemma \ref{l2}.  Because all the $\lambda$-clusters $\mathcal{C}_u,\mathcal{C}_1,\ldots,\mathcal{C}_N$ are contained in $\mathcal{C}(\tilde{\ell})$, $u\in\tilde{\ell}$ and $\mathcal{F}_3$ holds, it follows that all the routers of $\mathcal{C}_u,\mathcal{C}_1,\ldots,\mathcal{C}_N$ are connected in $\tilde{\mathcal{C}}(\tilde{\ell})$, and thus $\tilde{\mathcal{C}}(\tilde{\ell})=\tilde{\mathcal{C}}(\ell)$.  Then by Lemma \ref{l2} we obtain (\ref{e30}).

\begin{figure}
\centering
\includegraphics[height=0.3\textwidth]{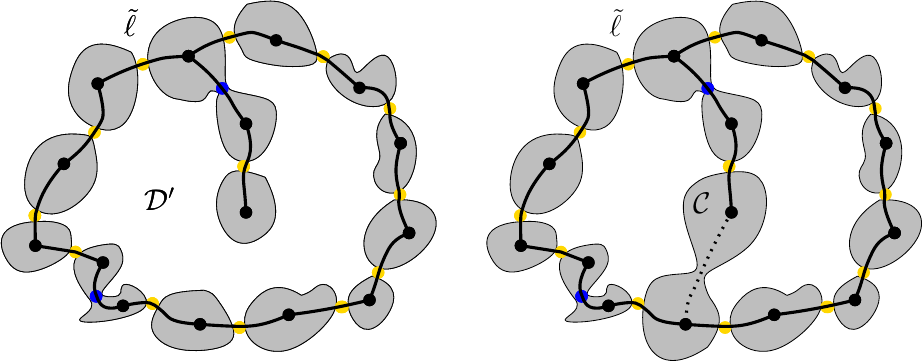}
\caption{\emph{Left}: The domain $\mathcal{D}'$. \emph{Right}: A cluster $\mathcal{C}$ in $\mathscr{C}(\tilde{\ell})$ so that $\mathcal{D}\setminus F(\mathcal{C})$ is not simply connected.  Since $\mathcal{F}$ occurs, such a scenario would produce a path (represented by the dotted line) of $\tilde{\textnormal{\textsf{G}}}_{\eta}^{\bar{\lambda},\varepsilon}$ connecting two distinct routers of $\mathcal{C}$ which belong to $\tilde{\ell}$, with the path lying in $\mathcal{D}$ except its starting and ending points, a contradiction.
}\label{jordan}
\end{figure}

\textbf{Case 2: $\tilde{\ell}$ is counterclockwise.}
Let $\mathcal{D}$ be the connected component of $\mathbb{C}\setminus\tilde{\ell}$ that is surrounded by $\tilde{\ell}$.  We first prove that $\mathcal{D}':=\mathcal{D}\setminus U(\tilde{\ell})$ is a Jordan domain; see the left-hand side of Figure \ref{jordan}.  Let $\mathcal{C}$ be a cluster in $\mathscr{C}(\tilde{\ell})$.  By Lemma \ref{l15}(\ref{topology}), all the routers of $\mathcal{C}$ are connected by paths of $\tilde{\textnormal{\textsf{G}}}_{\eta}^{\bar{\lambda},\varepsilon}$ in $F(\mathcal{C})$.  From this and the second paragraph of Lemma \ref{l15}(\ref{topology}) for an edge of $\tilde{\textnormal{\textsf{G}}}_{\eta}^{\bar{\lambda},\varepsilon}$ intersecting $F(\mathcal{C})$, together with the fact that $\mathcal{D}$ is simply connected, we deduce that $\mathcal{D}\setminus F(\mathcal{C})$ is simply connected, since otherwise there would exist a path of
$\tilde{\textnormal{\textsf{G}}}_{\eta}^{\bar{\lambda},\varepsilon}$ connecting two distinct routers of $\mathcal{C}$, such that the two routers lie in $\tilde{\ell}$ and the path lies in $\mathcal{D}$ except its starting and ending points, contradicting the definition of $\mathcal{D}$; see the right-hand side of Figure \ref{jordan}.  Then, it follows that $\mathcal{D}\setminus\bigcup_{\mathcal{C}\in\mathscr{C}(\tilde{\ell})}F(\mathcal{C})$ is simply connected since
for any two clusters $\mathcal{C},\mathcal{C}'\in\mathscr{C}(\tilde{\ell})$, either $F(\mathcal{C})$ and $F(\mathcal{C}')$ have empty intersection, or one of them is contained in the other.  Observe that for each $x\in X(\tilde{\ell})$, either the hexagon $H_x$ is in the interior of $F(\mathcal{C})$ for some $\mathcal{C}\in\mathscr{C}(\tilde{\ell})$, or $H_x$ contains a boundary piece of $\mathcal{D}\setminus \bigcup_{\mathcal{C}\in\mathscr{C}(\tilde{\ell})}F(\mathcal{C})$ which starts at a point in $\partial F(\mathcal{C})$ and ends at a point in $\partial F(\mathcal{C}')$ for some $\mathcal{C},\mathcal{C}'\in\mathscr{C}(\tilde{\ell})$ (it is possible that $\mathcal{C}=\mathcal{C}'$), and in the latter case there is no cluster in $\mathscr{C}(\tilde{\ell})\setminus\{\mathcal{C},\mathcal{C}'\}$ touching $H_x$, since otherwise $x$ would neighbor at least three routers in $\textnormal{\textsf{EN}}_M^{r\text{-}\textnormal{\textsf{meso}}}(\omega_{\eta}^{\lambda},\textnormal{\textsf{PPP}}_{\bar{\lambda}}^{\varepsilon})$, contradicting (\ref{num4}) of Proposition \ref{p6}.  The above observation and the fact that $\mathcal{D}\setminus\bigcup_{\mathcal{C}\in\mathscr{C}(\tilde{\ell})}F(\mathcal{C})$ is simply connected imply that $\mathcal{D}'$ is also simply connected and thus is a Jordan domain since $\mathcal{D}'$ is the interior of a union of $\eta$-hexagons.  So, $\partial\mathcal{D'}$ is a simple loop.

For any point $a\in\tilde{\ell}$, the above argument implies that there is some $\mathcal{C}\in\mathscr{C}(\tilde{\ell})$ such that $a\in F(\mathcal{C})$ and $F(\mathcal{C})$ touches $\partial\mathcal{D}'$, or there is some $x\in X(\tilde{\ell})$ such that $a\in H_x$ and $H_x$ touches $\partial\mathcal{D}'$.  Thus $\diam(\mathcal{D}')\geq\diam(\tilde{\ell})-s\geq 136s$ since all $\lambda$-clusters on $\mathbb{T}_M^2$ have diameter less than $s/2$.  We orient $\partial\mathcal{D}'$ counterclockwise.  It is obvious that we can partition $\partial\mathcal{D}'$ into $m$ oriented segments $\gamma_1,\ldots,\gamma_m$ with $m\geq 6$, such that the ending point of $\gamma_i$ is the starting point of $\gamma_{i+1}$ for each $1\leq i\leq m$, and all of $\gamma_1,\ldots,\gamma_{m-1}$ have diameter $42s$ while $\gamma_m$ has diameter in $[42s,84s)$.

Now, we claim that for each $1\leq i\leq m$, there is a $\lambda'$-closed path connecting $\gamma_i$ and $\gamma_{i+2}$ (we adopt here cyclic notation) in $\overline{\mathcal{D}'}$, and thus all such paths belong to a $\lambda'$-closed cluster, whose counterclockwise boundary loop is also a counterclockwise boundary loop of $\mathcal{C}(\tilde{\ell})$, denoted by $\ell$.

To verify the claim we just need to prove that for any pair $i,j\in\{1,\ldots,m\}$ with $j\neq i-1,i,i+1$, there are no $\lambda'$-open paths connecting $\gamma_i$ and $\gamma_j$ in $\overline{\mathcal{D}'}$.  Assume for a contradiction that there is a self-avoiding $\lambda'$-open path $\gamma$ in $\overline{\mathcal{D}'}$ connecting some $\gamma_i$ and some $\gamma_j$ with $j\neq i-1,i,i+1$.  Then we distinguish two cases:

\textbf{Case (a):} $\diam(\gamma)<9s$.  Observe that for any segment $\gamma'$ of $\partial\mathcal{D}'$ with $\diam(\gamma')\geq 5s$, there exists a sub-path $\tilde{\gamma}$ of $\tilde{\ell}$ with $\diam(\tilde{\gamma})\geq 3s$, such that any cluster in $\mathscr{C}(\tilde{\gamma})$ that is not contained in the hull of another cluster in $\mathscr{C}(\tilde{\gamma})$ touches $\gamma'$; since $\mathcal{G}$ holds, it follows that there is some cluster in $\mathscr{C}(\tilde{\gamma})$ of diameter at least $\delta$ touching $\gamma'$.  Choose two segments $\gamma_i'\subset\gamma_i,\gamma_j'\subset\gamma_j$ of diameter $5s$ such that $\gamma$ connects $\gamma_i'$ and $\gamma_j'$.  Then the above observation implies that there are $\mathcal{C}_i,\mathcal{C}_j\in\mathscr{C}(\tilde{\ell})$ of diameter at least $\delta$, such that $\mathcal{C}_i$ touches $\gamma_i'$ and $\mathcal{C}_j$ touches $\gamma_j'$.  Note that $\diam(\mathcal{C}_i),\diam(\mathcal{C}_j)< s/2,\diam(\gamma_i')=\diam(\gamma_j')=5s$ and $\diam(\gamma)<9s$.  Thus, there is a square $B\subset\Lambda_{5M/8+s}$ of side-length $21s$ so that $\mathcal{C}_i,\mathcal{C}_j\subset B$ and $\mathcal{C}_i,\mathcal{C}_j$ are connected by a $\lambda'$-open path in $B$.  Since $\mathcal{F}_3$ holds, there is a chain $\mathfrak{C}=(\mathcal{C}_i,x_1,\ldots,x_n,\mathcal{C}_j)$ of $\omega_{\eta}^{\lambda}$ such that $x_1,\ldots,x_n\in X(\tilde{\textnormal{\textsf{G}}}_{\eta}^{\bar{\lambda},\varepsilon})$ and
$U(\mathfrak{C})\subset B+\Lambda_s$, and there is a path $\tilde{\gamma}_{ij}$ of $\tilde{\textnormal{\textsf{G}}}_{\eta}^{\bar{\lambda},\varepsilon}$ from a router of $\mathcal{C}_i$ in $\tilde{\ell}$ to a router of $\mathcal{C}_j$ in $\tilde{\ell}$, such that all the routers in $\tilde{\gamma}_{ij}$ belong to the $\lambda$-clusters in $\mathfrak{C}$, with $\tilde{\gamma}_{ij}$ lying in the interior of $U(\mathfrak{C})$.  Notice that $U(\mathfrak{C})$ lies outside $\mathcal{D}'$, since otherwise $\tilde{\gamma}_{ij}$ would intersect $\mathcal{D}'$, contradicting the definition of $\tilde{\ell}$.  Therefore, the hull $F(U(\mathfrak{C})\cup\gamma_i'\cup\gamma_j'\cup\gamma)$ contains $\gamma_{i+1}$ or $\gamma_{i-1}$.  But the diameter of this hull is less than $42s$, which would imply that the diameter of $\gamma_{i+1}$ or $\gamma_{i-1}$ is less than $42s$, a contradiction.

\begin{figure}
\centering
\includegraphics[height=0.4\textwidth]{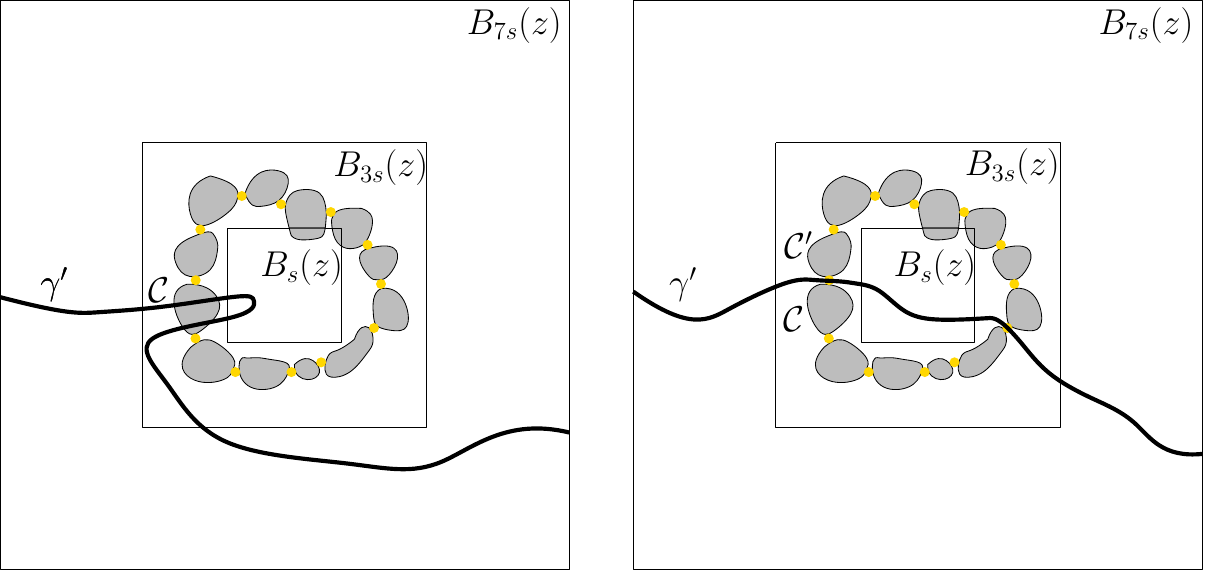}
\caption{Illustration of the segment $\gamma'$ of $\gamma$ intersecting a necklace, as in Case (b) in the proof of Lemma \ref{l18}.  \emph{Left}: The $\lambda$-cluster $\mathcal{C}$ is contained in $\overline{\mathcal{D}'}$, has diameter at least $\delta$, and intersects $\gamma'$ (subcase 1). \emph{Right}: The $\lambda$-clusters $\mathcal{C},\mathcal{C}'$ are contained in $B_{4s}(z)$, have diameter at least $\delta$, and are separated by $\gamma'$ in $B_{7s}(z)$, with a site in $\gamma'\cap B_{3s}(z)$ neighboring both of them.  Both of $\mathcal{C}$ and $\mathcal{C}'$ are in $\mathscr{C}(\tilde{\ell})$ (subcase 2), or at least one of them is contained in $\overline{\mathcal{D}'}$ (subcase 1).
}\label{necklacepath2}
\end{figure}
\textbf{Case (b):} $\diam(\gamma)\geq 9s$.  Clearly, we can find a box $B_s(z)$ of $\mathscr{B}^s$ centered at $z$ and a segment $\gamma'$ of $\gamma$, such that $\gamma'$ starts and ends in $\partial B_{7s}(z)$ with $\gamma'\cap B_s(z)\neq \emptyset$, and $\gamma'$ is contained in the interior of $B_{7s}(z)$ except for its starting and ending points.  Since $\mathcal{F}_2$ holds and all the $\lambda$-clusters on $\mathbb{T}_M^2$ have diameter less than $s/2$, by considering a necklace in the annulus of $\mathscr{A}_{M,s}$ centered at $B_s(z)$ we get that there are two subcases (see Figure \ref{necklacepath2}): 1. There exists a $\lambda$-cluster $\mathcal{C}$ of diameter at least $\delta$, such that $\mathcal{C}$ is contained in $\overline{\mathcal{D}'}$ and intersects or touches $\gamma'$.  Thus, there is a $\lambda'$-open path in $\gamma\cup U(\tilde{\ell})$ connecting $\mathcal{C}$ and $\mathcal{C}_u$.  Then by the occurrence of $\mathcal{F}_1\cap\mathcal{F}_3$, $\mathcal{C}$ has a router $v$ connected to $\tilde{\ell}$ in $\tilde{\textnormal{\textsf{G}}}_{\eta}^{\bar{\lambda},\varepsilon}$.  This and the fact that $v\in\mathcal{D}$ contradict the definition of $\mathcal{D}$.  2. There are two distinct clusters $\mathcal{C},\mathcal{C}'\in\mathscr{C}(\tilde{\ell})$ of diameter at least $\delta$, such that they are contained in $B_{4s}(z)$ and are separated by $\gamma'$ in $B_{7s}(z)$, with a site in $\gamma'\cap B_{3s}(z)$ neighboring both of them.  Then, since $\mathcal{F}_1\cap\mathcal{F}_3$ holds, there is a chain $\mathfrak{C}=(\mathcal{C},x_1,\ldots,x_n,\mathcal{C}')$ of $\omega_{\eta}^{\lambda}$ such that $U(\mathfrak{C})\subset B_{6s}(z)$, and moreover, a router of $\mathcal{C}$ in $\tilde{\ell}$ is connected to a router of $\mathcal{C}'$ in $\tilde{\ell}$ by a path of $\tilde{\textnormal{\textsf{G}}}_{\eta}^{\bar{\lambda},\varepsilon}$ in $U(\mathfrak{C})$.  But this path would cross $\gamma$ since $\gamma$ separates $\mathcal{C}$ and $\mathcal{C}'$ in $B_{7s}(z)$, and thus would intersect $\mathcal{D}$, a contradiction.

So our claim follows.  It is obvious that $\ell$ is surrounded by $\tilde{\ell}$ and touches all $\gamma_1,\ldots,\gamma_m$.  Thus, $\diam(\ell)>\diam(\mathcal{D}')-42s-84s\geq 10s$.  Let $\tilde{\mathcal{C}}(\ell)$ be as in Lemma \ref{l2}.  Similarly to the argument at the end of Case 1, by using Lemma \ref{l2} we obtain $\tilde{\mathcal{C}}(\tilde{\ell})=\tilde{\mathcal{C}}(\ell)$ and, moreover, $\tilde{\ell}$ is the unique counterclockwise boundary loop of $\tilde{\mathcal{C}}(\tilde{\ell})$ surrounding $\ell$, and (\ref{e30}) holds.
\end{proof}

The following result combines Lemmas \ref{l2} and \ref{l18}, with a slight generalization to small loops.
\begin{lemma}\label{l5}
Fix any $\lambda'\in \mathbb{R}$ and any $M\geq 8$.  For any $\alpha>0$, if $\lambda<\lambda'$ is very negative and $\varepsilon>0$ is small enough, then for all small enough $\eta>0$, we have
\begin{equation}\label{e10}
\mathbf{P}\left[
\begin{aligned}
&\mbox{for each $\ell\in\Gamma_{\eta}^{\lambda'}(3M/4)$, there exists $\tilde{\ell}\in\tilde{\Gamma}_{\eta}^{\bar{\lambda},\varepsilon}$, s.t. $d(\ell,\tilde{\ell})<\alpha$;}\\
&\mbox{for each $\tilde{\ell}\in\tilde{\Gamma}_{\eta}^{\bar{\lambda},\varepsilon}(5M/8)$, there exists $\ell\in\Gamma_{\eta}^{\lambda'}$, s.t. $d(\ell,\tilde{\ell})<\alpha$}
\end{aligned}
\right]>1-\alpha.
\end{equation}
\end{lemma}
\begin{proof}
Note that for each $\ell\in\Gamma_{\eta}^{\lambda'}(3M/4)$ with diameter smaller than $\alpha$, there is a trivial loop $\tilde{\ell}$ of $\tilde{\Gamma}_{\eta}^{\bar{\lambda},\varepsilon}$ at a point of $\ell$, such that $d(\ell,\tilde{\ell})<\alpha$.  Combining this with
Lemmas \ref{l19} and \ref{l2}, we get that for any $\alpha>0$, if $\lambda<\lambda'$ is very negative and $\varepsilon>0$ is small enough, then for all small enough $\eta>0$,
\begin{equation}\label{e15}
\mathbf{P}[\mbox{for each $\ell\in\Gamma_{\eta}^{\lambda'}(3M/4)$, there exists $\tilde{\ell}\in\tilde{\Gamma}_{\eta}^{\bar{\lambda},\varepsilon}$, s.t. $d(\ell,\tilde{\ell})<\alpha$}]>1-\alpha/2.
\end{equation}

It remains to prove the second part.  Using RSW (see Section \ref{classic}) and FKG, it is easy to show that for any fixed $\alpha>0$ and all small enough $\eta>0$, with probability at least $1-\alpha/4$, for each point $z\in\Lambda_{5M/8}$, there exists a loop in $\Gamma_{\eta}^{\lambda'}$ with diameter smaller than $\alpha/4$, such that the distance between this loop and $z$ is smaller than $\alpha/4$, and thus for each $\tilde{\ell}\in\tilde{\Gamma}_{\eta}^{\bar{\lambda},\varepsilon}(5M/8)$ with $\diam(\tilde{\ell})<\alpha/2$, there exists $\ell\in\Gamma_{\eta}^{\lambda'}$ with $\diam(\ell)<\alpha/4$, such that $d(\ell,\tilde{\ell})<\alpha$.  This together with Lemmas \ref{l19} and \ref{l18} implies that for any $\alpha>0$, if $\lambda<\lambda'$ is very negative and $\varepsilon>0$ is small enough, then for all small enough $\eta>0$,
\begin{equation}\label{e34}
\mathbf{P}[\mbox{for each $\tilde{\ell}\in\tilde{\Gamma}_{\eta}^{\bar{\lambda},\varepsilon}(5M/8)$, there exists $\ell\in\Gamma_{\eta}^{\lambda'}$, s.t. $d(\ell,\tilde{\ell})<\alpha$}]>1-\alpha/2.
\end{equation}
Combining (\ref{e15}) and (\ref{e34}) gives (\ref{e10}).
\end{proof}

\subsection{Proof of Theorem \ref{t1}}
For notational clarity, we will use $\omega_{\eta}^{\lambda,M}$ and $\omega_{\infty}^{\lambda,M}$ for $\omega_{\eta}^{\lambda}$ on the torus $\mathbb{T}_M^2$ and its scaling limit, respectively.  Similarly, we will use $\Gamma_{\infty}^{\bar{\lambda},\varepsilon,M}$ and $\tilde{\Gamma}_{\eta}^{\bar{\lambda},\varepsilon,M}$
for $\Gamma_{\infty}^{\bar{\lambda},\varepsilon}$ and $\tilde{\Gamma}_{\eta}^{\bar{\lambda},\varepsilon}$ defined in Definition \ref{def2}, respectively.  Moreover, $\textsf{PPP}_{\bar{\lambda}}^{\varepsilon,M}$ will denote $\textsf{PPP}_{\bar{\lambda}}^{\varepsilon}$ in Definition \ref{def8}.

\begin{proof}[Proof of Theorem \ref{t1}]
Using Proposition \ref{p5} and Lemmas \ref{l14} and \ref{l5}, the proof is classical; a similar strategy was used in \cite[Section 5]{GPS18b} and \cite[Section 9]{GPS18a}.

Fix any $\lambda'\in\mathbb{R}\setminus\{0\}$ and let $j\geq 2$.  By Lemma \ref{l5}, there exists $\lambda_j'<\lambda'$, such that for each $\lambda<\lambda_j'$ there exists $0<\varepsilon_j=\varepsilon_j(\lambda)<2^{-j}/100$, such that for all $\eta>0$ sufficiently small,
\begin{equation}\label{e7}
\mathbf{P}\left[
\begin{aligned}
&\mbox{for each $\ell\in\Gamma_{\eta}^{\lambda'}(2^{j})$, there exists $\tilde{\ell}\in\tilde{\Gamma}_{\eta}^{\bar{\lambda},\varepsilon_j,2^{j+1}}$, s.t. $d(\ell,\tilde{\ell})<2^{-j}$;}\\
&\mbox{for each $\tilde{\ell}\in\tilde{\Gamma}_{\eta}^{\bar{\lambda},\varepsilon_j,2^{j+1}}(2^{j}+2^{-j})$, there exists $\ell\in\Gamma_{\eta}^{\lambda'}$, s.t. $d(\ell,\tilde{\ell})<2^{-j}$}
\end{aligned}
\right]>1-2^{-j}.
\end{equation}

By Lemma \ref{l14} and Proposition \ref{p5}, there exists $\lambda_j<\lambda_j'$, such that for $\bar{\lambda}_j=(\lambda_j,\lambda')$ and $\varepsilon_j=\varepsilon_j(\lambda_j)$ there is a coupling between $(\omega_{\eta}^{\lambda_j,2^{j+1}},\textsf{PPP}_{\bar{\lambda}_j}^{\varepsilon_j,2^{j+1}})$ and $(\omega_{\infty}^{\lambda_j,2^{j+1}},\textsf{PPP}_{\bar{\lambda}_j}^{\varepsilon_j,2^{j+1}})$, so that for all $\eta>0$ sufficiently small,
\begin{equation}\label{e8}
\mathbf{P}\left[
\begin{aligned}
&\mbox{for each $\tilde{\ell}\in\tilde{\Gamma}_{\eta}^{\bar{\lambda}_j,\varepsilon_j,2^{j+1}}(2^j+2^{-j})$, there exists $\ell_{\infty}\in\Gamma_{\infty}^{\bar{\lambda}_j,\varepsilon_j,2^{j+1}}$, s.t. $d(\tilde{\ell},\ell_{\infty})<2^{-j}$;}\\
&\mbox{for each $\ell_{\infty}\in\Gamma_{\infty}^{\bar{\lambda}_j,\varepsilon_j,2^{j+1}}(2^j)$, there exists $\tilde{\ell}\in\tilde{\Gamma}_{\eta}^{\bar{\lambda}_j,\varepsilon_j,2^{j+1}}$, s.t. $d(\tilde{\ell},\ell_{\infty})<2^{-j}$}
\end{aligned}
\right]>1-2^{-j}.
\end{equation}

Combining (\ref{e7}) and (\ref{e8}) and using the triangle inequality, we know that there is a coupling between $\omega_{\eta}^{[\lambda_j,\lambda']}$ and $(\omega_{\infty}^{\lambda_j,2^{j+1}},\textsf{PPP}_{\bar{\lambda}_j}^{\varepsilon_j,2^{j+1}})$, such that for all $\eta>0$ sufficiently small,
\begin{equation}\label{e39}
\mathbf{P}\left[
\begin{aligned}
&\mbox{for each $\ell\in\Gamma_{\eta}^{\lambda'}(2^j)$, there exists $\ell_{\infty}\in\Gamma_{\infty}^{\bar{\lambda}_j,\varepsilon_j,2^{j+1}}$, s.t. $d(\ell,\ell_{\infty})<2^{-j+1}$;}\\
&\mbox{for each $\ell_{\infty}\in\Gamma_{\infty}^{\bar{\lambda}_j,\varepsilon_j,2^{j+1}}(2^j)$, there exists $\ell\in\Gamma_{\eta}^{\lambda'}$, s.t. $d(\ell,\ell_{\infty})<2^{-j+1}$}
\end{aligned}
\right]>1-2^{-j+1}.
\end{equation}

We shall extend the estimate (\ref{e39}) to the case of the full plane.  Recall the metrics $\Delta(\cdot,\cdot)$ and $D(\cdot,\cdot)$ defined in Section \ref{loop}.  Note that there is an absolute constant $C>1$, such that for all $M\geq 1$ the $\Delta$-distance between any two points in $\hat{\mathbb{C}}\setminus\Lambda_M$ is smaller than $C/M$, and for any two loops $\ell,\ell'$ on $\mathbb{C}$, we have $D(\ell,\ell')\leq Cd(\ell,\ell')$.

By the half-plane 3-arm exponent being 2 (see (\ref{e70})), it is easy to get that there is a constant $C_1>0$ such that for each $j\geq 2$ and all $\eta>0$ small enough,
\begin{equation}\label{e37}
\mathbf{P}_{\eta}^{\lambda'}[\mbox{for all $\ell\in\Gamma_{\eta}^{\lambda'}(2^j+2^{-j+1})$ intersecting $A(2^j-2^{-j+1},2^j+2^{-j+1})$, $\diam(\ell)<1$}]>1-C_12^{-j}.
\end{equation}
Denote by $\bullet(z)$ the trivial loop at $z\in\hat{\mathbb{C}}$.  Observe that when the events in (\ref{e37}) and (\ref{e39}) occur simultaneously, then the following holds: If a loop $\ell$ of $\Gamma_{\eta}^{\lambda'}(2^j)$ intersects ${A(2^j-2^{-j+1},2^j)}$, then $\ell\subset A(2^j-2^{-j+1}-1,2^j)$ and thus $D(\ell,\bullet(2^j))<C2^{-j+1}$; if a loop $\ell_{\infty}$ of $\Gamma_{\infty}^{\bar{\lambda}_j,\varepsilon_j,2^{j+1}}(2^j)$ intersects $A(2^j-2^{-j+1},2^j)$ and there exists $\ell\in\Gamma_{\eta}^{\lambda'}$ such that $\ell\nsubseteq\Lambda_{2^j}$ and $d(\ell,\ell_{\infty})<2^{-j+1}$, then $\ell\subset A(2^j-1,2^j+2^{-j+1})$ and $\ell_{\infty}\subset A(2^j-2^{-j+1}-1,2^j)$, and thus $D(\ell_{\infty},\bullet(\infty))<C2^{-j+1}$.  Hence, combining (\ref{e37}) and (\ref{e39}), we have
\begin{equation}\label{e38}
\mathbf{P}[D_{\mathcal {L}}(\Gamma_{\eta}^{\lambda'}(2^j)\cup\{\bullet(\infty)\},\Gamma_{\infty}^{\bar{\lambda}_j,\varepsilon_j,2^{j+1}}(2^j))<C2^{-j+1}]>1-(2+C_1)2^{-j}.
\end{equation}

By (\ref{e43}), there exists $C_2>0$ such that for all $j\geq 2$ and all $\eta>0$ small enough, with probability at least $1-\exp(-C_22^j)$, there exists a closed or open circuit of $\omega_{\eta}^{\lambda'}$ surrounding 0 in $A(2^{j-1},2^j)$, and thus all the loops of $\Gamma_{\eta}^{\lambda'}$ that intersect $\Lambda_{2^j}$ but do not lie in $\Lambda_{2^j}$ must lie outside $\Lambda_{2^{j-1}}$.  This implies that
\begin{equation*}
\mathbf{P}_{\eta}^{\lambda'}[D_{\mathcal {L}}(\Gamma_{\eta}^{\lambda'},\Gamma_{\eta}^{\lambda'}(2^j)\cup\{\bullet(\infty)\})<C2^{-j+1}]>1-\exp(-C_22^j).
\end{equation*}
Combined with (\ref{e38}), this leads to
\begin{equation}\label{e36}
\mathbf{P}[D_{\mathcal {L}}(\Gamma_{\eta}^{\lambda'},\Gamma_{\infty}^{\bar{\lambda}_j,\varepsilon_j,2^{j+1}}(2^j))<C2^{-j+2}]>1-(2+C_1)2^{-j}-\exp(-C_22^j).
\end{equation}
Thus, we can now couple all the loop ensembles $\Gamma_{\infty}^{\bar{\lambda}_j,\varepsilon_j,2^{j+1}}(2^j)$ to $\Gamma_{\eta}^{\lambda'}$ one-by-one and use the triangular inequality to obtain that, there exists an absolute constant $C_3>0$ such that for all $j\geq 2$,
\begin{equation}\label{e9}
\mathbf{P}[D\texttt{}_{\mathcal {L}}(\Gamma_{\infty}^{\bar{\lambda}_j,\varepsilon_j,2^{j+1}}(2^j),\Gamma_{\infty}^{\bar{\lambda}_i,\varepsilon_i,2^{i+1}}(2^i))<C2^{-j+3}\mbox{ for all $i\geq j$}]>1-2^{-C_3j}.
\end{equation}
Using the Borel-Cantelli lemma for the events in (\ref{e9}) shows that the sequence
$\Gamma_{\infty}^{\bar{\lambda}_j,\varepsilon_j,2^{j+1}}(2^j)$ is almost surely a Cauchy sequence in $(\mathcal {L},D_{\mathcal {L}})$.  This space is complete, hence there is an almost sure limit, denoted by $\Gamma_{\infty}^{\lambda'}$.  This limit may a priori depend on the sequences $\{\bar{\lambda}_j\},\{\varepsilon_j\}$ and on the coupling.  However, using (\ref{e36}) and the triangle inequality again, it follows that for any fixed $\alpha>0$ and all $\eta>0$ small enough,
\begin{equation*}
\mathbf{P}[D_{\mathcal {L}}(\Gamma_{\eta}^{\lambda'},\Gamma_{\infty}^{\lambda'})<\alpha]>1-\alpha.
\end{equation*}
Therefore, in this coupling, $\Gamma_{\eta}^{\lambda'}$ converges in probability to $\Gamma_{\infty}^{\lambda'}$ in $(\mathcal {L},D_{\mathcal {L}})$, which implies the theorem.
\end{proof}

\section{Basic properties of the limiting loop ensemble $\Gamma_{\infty}^{\lambda}$}\label{property}

In this section we prove Theorems \ref{t7} and \ref{t2}.

\begin{proof}[Proof of Theorem \ref{t7}]
In the proof, we will assume that $\lambda'\in\mathbb{R}\setminus\{0\}$ is fixed and represents $\lambda$ in the theorem.

(\ref{rotation}) It is clear that $\Gamma_{\eta}^{\lambda'}$ is invariant under lattice shifts of $\eta\textsf{T}$.  Then the fact that $\Gamma_{\infty}^{\lambda'}$ is invariant under translations follows from this and Theorem \ref{t1}.  Now let us show the rotational invariance of $\Gamma_{\infty}^{\lambda'}$, whose proof is similar to that for the rotational invariance of the scaling limit of minimal spanning trees in \cite[Section 5.2]{GPS18b}.  Let $f_{\theta}:\mathbb{C}\rightarrow\mathbb{C}$ be the rotation by angle $\theta$.  By (5.7) of \cite{GPS18b}, we have for full plane configurations,
\begin{equation}\label{e45}
f_{\theta}(\omega_{\infty}^{\lambda},\textsf{PPP}_{\bar{\lambda}}^{\varepsilon}(\omega_{\infty}^{\lambda}))\overset{d}{=}(\omega_{\infty}^{\lambda},
\textsf{PPP}_{\bar{\lambda}}^{\varepsilon,\theta}(\omega_{\infty}^{\lambda})),
\end{equation}
where $\textsf{PPP}_{\bar{\lambda}}^{\varepsilon,\theta}$ is constructed using a $\theta$-rotated grid to define $\varepsilon$-importance.

Write $\Lambda_M^{\theta}:=f_{\theta}(\Lambda_M)$, and let $\mathbb{T}_M^{2,\theta}$ denote the torus that is realized as the subdomain $f_{\theta}([-M,M)^2)$ of $\mathbb{C}$ with a periodic boundary condition.  Similarly to the construction of $\Gamma^{\bar{\lambda},\varepsilon,M}$ in Definition \ref{def2}, we can use $\omega^{\lambda}$ and $\textsf{PPP}_{\bar{\lambda}}^{\varepsilon,\theta}(\omega^{\lambda})$ on $\mathbb{T}_M^{2,\theta}$ together with the $\theta$-rotated $r$-grid and the $\theta$-rotated $(s,k,\hat{r})$-representation to construct a loop ensemble in $\Lambda_M^{\theta}$, denoted by $\Gamma^{\bar{\lambda},\varepsilon,M,\theta}$.

By the construction of $\Gamma^{\bar{\lambda},\varepsilon,2M}$ we get that when $\lambda<0$ is very negative and both $\varepsilon$ and $r$ are small, with high probability $\Gamma^{\bar{\lambda},\varepsilon,2M}|_{\Lambda_M}$ only depends on $(\omega^{\lambda}|_{\Lambda_{2M}},\textsf{PPP}_{\bar{\lambda}}^{\varepsilon}(\omega^{\lambda}|_{\Lambda_{2M}}))$, not on the torus structure.  By the proof of Theorem \ref{t1}, we can choose a large $M$ and a very negative $\lambda$, together with small enough $\varepsilon,r$ and $\eta$, such that with high probability $\Gamma_{\infty}^{\bar{\lambda},\varepsilon,2M}|_{\Lambda_M}$ is close to $\Gamma_{\infty}^{\lambda'}$ and $\Gamma_{\eta}^{\lambda'}$.  Similarly, one can obtain that, with high probability $\Gamma_{\infty}^{\bar{\lambda},\varepsilon,2M,\theta}|_{\Lambda_M^{\theta}}$ is close to $\Gamma_{\infty}^{\lambda'}$ and $\Gamma_{\eta}^{\lambda'}$.  Using (\ref{e45}), we have $f_{\theta}(\Gamma_{\infty}^{\bar{\lambda},\varepsilon,2M}|_{\Lambda_M})\overset{d}{=}\Gamma_{\infty}^{\bar{\lambda},\varepsilon,2M,\theta}|_{\Lambda_M^{\theta}}$.  Thus,
$f_{\theta}(\Gamma_{\infty}^{\bar{\lambda},\varepsilon,2M}|_{\Lambda_M})$, which is close to $f_{\theta}(\Gamma_{\infty}^{\lambda'})$ with high probability, is also close to $\Gamma_{\infty}^{\lambda'}$ with high probability under some coupling.  After taking all the limits we get that
$f_{\theta}(\Gamma_{\infty}^{\lambda'})$ and $\Gamma_{\infty}^{\lambda'}$ have the same distribution.

(\ref{scaling}) The proof is analogous to that of scale invariance of the scaling limit of minimal spanning trees in \cite[Section 5.2]{GPS18b}, though the property here is scale covariance.  Let $X_{\bar{\lambda}}^{\varepsilon}(\omega_{\infty}^{\lambda,M})$ denote the set of pivotals in $\textsf{PPP}_{\bar{\lambda}}^{\varepsilon}(\omega_{\infty}^{\lambda,M})$.  By (5.5) and (5.6) in \cite{GPS18b}, we have
\begin{equation}\label{e46}
f_{\alpha}(\omega_{\infty}^{\lambda,M})\overset{d}{=}\omega_{\infty}^{\alpha^{-3/4}\lambda,\alpha^{-3/4}M},\quad f_{\alpha}(X_{\bar{\lambda}}^{\varepsilon}(\omega_{\infty}^{\lambda,M}))\overset{d}{=}
X_{\alpha^{-3/4}\bar{\lambda}}^{\alpha\varepsilon}(\omega_{\infty}^{\alpha^{-3/4}\lambda,\alpha^{-3/4}M}).
\end{equation}
Our construction of $\Gamma_{\infty}^{\bar{\lambda},\varepsilon,M}$ in Definition \ref{def2} uses the graph $\textnormal{\textsf{G}}_{\infty}^{\bar{\lambda},\varepsilon,M}$ constructed in Definition \ref{def7}, which relies on a chosen $(s,k,\hat{r})$-representation of $\textsf{EN}_M^{r\text{-}\textsf{meso}}(\omega_{\infty}^{\lambda,M},X_{\bar{\lambda}}^{\varepsilon}(\omega_{\infty}^{\lambda,M}))$.  One can slightly change our definition of $\Gamma_{\infty}^{\bar{\lambda},\varepsilon,M}$ to make the construction equivariant under scalings.  More precisely, for a given configuration $(f_{\alpha}(\omega_{\infty}^{\lambda,M}),f_{\alpha}(X_{\bar{\lambda}}^{\varepsilon}(\omega_{\infty}^{\lambda,M})))$ on $\mathbb{T}_{\alpha M}^2$, we let the chosen representation of $\textsf{EN}_{\alpha M}^{\alpha r\text{-}\textsf{meso}}(f_{\alpha}(\omega_{\infty}^{\lambda,M}),f_{\alpha}(X_{\bar{\lambda}}^{\varepsilon}(\omega_{\infty}^{\lambda,M})))$ be the $f_{\alpha}$-scaling of the chosen representation of $\textsf{EN}_M^{r\text{-}\textsf{meso}}(\omega_{\infty}^{\lambda,M},X_{\bar{\lambda}}^{\varepsilon}(\omega_{\infty}^{\lambda,M}))$, and let the loop ensemble constructed from the former representation be the $f_{\alpha}$-scaling of the loop ensemble constructed from the latter.  We denote such modified version of $\Gamma_{\infty}^{\bar{\lambda},\varepsilon,M}$ by $\hat{\Gamma}_{\infty}^{\bar{\lambda},\varepsilon,M}$.  Then (\ref{e46})
implies that
\begin{equation}\label{e47}
f_{\alpha}(\hat{\Gamma}_{\infty}^{\bar{\lambda},\varepsilon,M})\overset{d}{=}\hat{\Gamma}_{\infty}^{\alpha^{-3/4}\bar{\lambda},\alpha\varepsilon,\alpha M}.
\end{equation}
Our proof of Theorem \ref{t1} also works for $\hat{\Gamma}_{\infty}^{\bar{\lambda},\varepsilon,M}$, and one can obtain that there exist sequences $\{\lambda_j\}_{j\geq 2},\{\varepsilon_j\}_{j\geq 2}$ such that, as $j\rightarrow\infty$, $\hat{\Gamma}_{\infty}^{(\lambda_j,\lambda'),\varepsilon_j,2^{j+1}}(2^j)$ converges to $\Gamma_{\infty}^{\lambda'}$ (i.e., the scaling limit  of $\Gamma_{\eta}^{\lambda'}$), and $\hat{\Gamma}_{\infty}^{\alpha^{-3/4}(\lambda_j,\lambda'),\alpha\varepsilon_j,\alpha 2^{j+1}}(\alpha 2^j)$ converges to $\Gamma_{\infty}^{\alpha^{-3/4}\lambda'}$.  From this and (\ref{e47}) we have $f_{\alpha}(\Gamma_{\infty}^{\lambda'})\overset{d}{=}\Gamma_{\infty}^{\alpha^{-3/4}\lambda'}$.

(\ref{reversing}) By the definition of $(\omega_{\eta}^{\lambda'})_{\lambda'\in\mathbb{R}}$, $\omega_{\eta}^{\lambda'}$ has the same distribution as the color-switching of $\omega_{\eta}^{-\lambda'}$.  Then, the definition of $\Gamma_{\eta}^{\lambda'}$ gives that $f(\Gamma_{\eta}^{\lambda'})$ and $\Gamma_{\eta}^{-\lambda'}$ have the same distribution for the orientation-reversing transformation $f$.  From this and Theorem \ref{t1} we get (\ref{reversing}) immediately.
\end{proof}

\begin{proof}[Proof of Theorem \ref{t2}]
The proofs of (\ref{local}), (\ref{noinfinite}), (\ref{triple}), (\ref{parity}) and (\ref{largest}) are standard, while (\ref{chaining}) needs a somewhat new approach; here we give the proofs of all these properties for completeness.

Since Theorem \ref{t1} allows us to transfer discrete results to the continuum case, it will be convenient to rely on the discrete setting.  Note that Theorem \ref{t1} implies that there is a coupling $\mathbf{P}$ of $(\Gamma_{\eta}^{\lambda})_{\eta>0}$ and $\Gamma_{\infty}^{\lambda}$ such that $\Gamma_{\eta}^{\lambda}$
converges pointwise in $(\mathcal{L},D_{\mathcal{L}})$ towards $\Gamma_{\infty}^{\lambda}$.  We shall work with $\mathbf{P}$ throughout the proof.

(\ref{local})  Note that if a loop of diameter larger than $\varepsilon$ intersects an $\varepsilon/3$-square, then it radially crosses the square annulus of exterior side-length $\varepsilon$ centered around this square.  By Theorem \ref{t1}, RSW and the BK inequality, there exists $p_0\in(0,1)$, such that for any $0<\varepsilon<\tilde{L}(\lambda)$, $k\in\mathbb{N}$ and $z\in\mathbb{C}$,
\begin{align*}
&\mathbf{P}[\mbox{there exist at least $k$ loops of $\Gamma_{\infty}^{\lambda}$ radially crossing $A(z;\varepsilon/6,\varepsilon/2)$}]\\
&\leq\limsup_{\eta\rightarrow 0}\mathbf{P}[\mbox{there exist at least $k$ loops of $\Gamma_{\eta}^{\lambda}$ radially crossing $A(z;\varepsilon/5,\varepsilon/3)$}]\\
&\leq\limsup_{\eta\rightarrow 0}\mathbf{P}[\mbox{there exist at least $k$ disjoint open arms radially crossing $A(z;\varepsilon/5,\varepsilon/3)$}]\\
&\leq p_0^k.
\end{align*}
Therefore, for each $\varepsilon/3$-square of the lattice $(\varepsilon/3)\mathbb{Z}^2$, there a.s. exist finitely many loops of $\Gamma_{\infty}^{\lambda}$ of diameter larger than $\varepsilon$ intersecting the square, which implies (\ref{local}).

(\ref{noinfinite}) Using (\ref{e43}) and Theorem \ref{t1}, one can deduce that almost surely, for all $k\in\mathbb{N}$ there is no infinite loop of $\Gamma_{\infty}^{\lambda}$ intersecting $\Lambda_k$.  This implies that $\Gamma_{\infty}^{\lambda}$ a.s. has no infinite loop.

(\ref{chaining}) The finite chaining property for $\Gamma_{\infty}^{\lambda=0}$ was given in \cite[Theorem 2, item 3]{CN06}, whose proof relies on the construction of CLE$_6$ by pasting together portions of SLE$_6$ paths and the fact that any given disk is almost surely surrounded by a loop in the full-plane CLE$_6$.  However, this method does not work for $\lambda\neq 0$, since analogous construction for $\Gamma_{\infty}^{\lambda}$ has not been established and the above loop-surrounding property does not hold.  Instead, our proof uses some basic properties of cluster chains.  In the following we assume that $\lambda<0$, and the $\lambda>0$ case follows from the $\lambda<0$ case by a color-switching trick.

Recall that in Definition \ref{chaindef} we have defined chains of open clusters.  Here, for simplicity, we call a finite sequence of distinct open (resp. closed) clusters $\mathcal{C}_1,\ldots,\mathcal{C}_k$ in $\omega_{\eta}^{\lambda}$ a \textbf{chain of open (resp. closed) clusters}, if for each $i\in\{1,\ldots,k-1\}$, there is a closed (resp. open) site neighboring both $\mathcal{C}_i$ and $\mathcal{C}_{i+1}$.  For $0<\varepsilon<\delta<1$ and $K,N,M\in\mathbb{N}$, let $\mathcal{F}(K,\varepsilon,\delta,N,M)$ be the event that any two open clusters of $\omega_{\eta}^{\lambda}\cap\Lambda_{N}$ of diameter at least $\delta N$ are connected by a chain of at most $K$ open clusters of $\omega_{\eta}^{\lambda}\cap\Lambda_{N}$, each of diameter at least $\varepsilon N$, and moreover, any two closed clusters of $\omega_{\eta}^{\lambda}\cap\Lambda_{M}$ of diameter at least $\delta N$ are connected by a chain of at most $K$ closed clusters of $\omega_{\eta}^{\lambda}\cap\Lambda_{M}$, each of diameter at least $\varepsilon N$.  We claim that for any fixed $N,M\in\mathbb{N},\delta\in(0,1)$ and any $\alpha>0$, there exist $K\in\mathbb{N}$ and $\varepsilon\in(0,\delta)$ such that for all $\eta>0$ small enough,
\begin{equation}\label{e48}
\mathbf{P}_{\eta}^{\lambda}[\mathcal{F}(K,\varepsilon,\delta,N,M)]\geq 1-\alpha.
\end{equation}
The proof of (\ref{e48}) is the same as that of \cite[Corollary 7.6]{DCMT21} for critical random-cluster model on $\mathbb{Z}^2$ (including critical Bernoulli bond percolation), and Corollary 7.6 relies on \cite[Theorem 7.5]{DCMT21}; the proofs of these two results in \cite{DCMT21} used RSW, estimates of arm events and quasi-multiplicativity, which are available for our near-critical site percolation $\omega_{\eta}^{\lambda}$ constrained in $\Lambda_{N}$ or $\Lambda_{M}$, with $\lambda,N$ and $M$ fixed.  We omit the details.

Analogously to the definition of cluster chains, a \textbf{(discrete) chain of loops} is a finite sequence of distinct loops $\ell_1,\ldots,\ell_k$ in $\Gamma_{\eta}^{\lambda}$, such that for each $i\in\{1,\ldots,k-1\}$, there is a site of $\eta\textsf{T}$ neighboring both $\ell_i$ and $\ell_{i+1}$.  In particular, a single loop is also a chain.  Similarly, a \textbf{(continuum) chain of loops} is a finite sequence of distinct loops $\ell_1,\ldots,\ell_k$ in $\Gamma_{\infty}^{\lambda}$, such that for each $i\in\{1,\ldots,k-1\}$, $\ell_i\cap\ell_{i+1}\neq\emptyset$.

For a chain $\mathscr{L}$ of loops, let $|\mathscr{L}|$ be the number of loops in $\mathscr{L}$.  For $1\leq N<M$ and two loops $\ell,\ell'$ of $\Gamma_{\eta}^{\lambda}$ contained in $\Lambda_M$, we define
\begin{equation*}
\mathbf{N}(\ell,\ell';\varepsilon,N,M):=\min\left\{|\mathscr{L}|:\begin{aligned}
&\mbox{$\mathscr{L}$ is a chain of loops from $\ell$ to $\ell'$ in $\Gamma_{\eta}^{\lambda}$, s.t. every loop}\\
&\mbox{of $\mathscr{L}$ is contained in $\Lambda_{M}$ and has diameter at least $\varepsilon N$}
\end{aligned}
\right\}.
\end{equation*}
\begin{figure}
\centering
\includegraphics[height=0.8\textwidth]{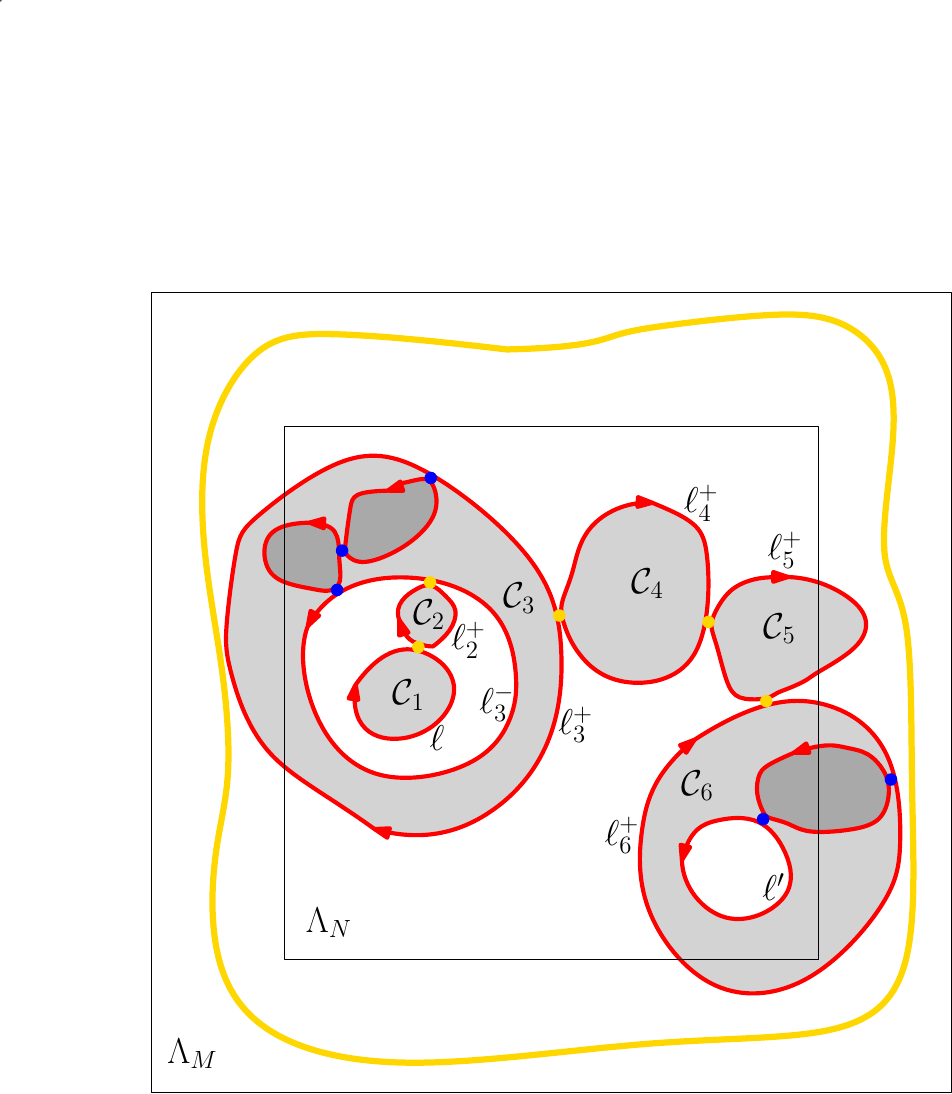}
\caption{Constructing a chain of loops from $\ell$ to $\ell'$ via a chain of open clusters and chains of closed clusters in the proof of Theorem \ref{t2}(\ref{chaining}).  The open clusters in the chain $(\mathcal{C}_1,\ldots,\mathcal{C}_{k=6})$ are drawn in light gray; closed clusters are drawn in dark gray.  The blue and yellow dots represent the open and closed sites neighboring two loops in the loop-chain, respectively.  A closed circuit surrounding 0 in the annulus $A(N,M)$ is shown using a yellow line.}\label{finitechaining}
\end{figure}

For $1\leq N<M$, let $\textrm{Circuit}(N,M)$ denote the event that there is a closed circuit of $\omega_{\eta}^{\lambda}$ surrounding 0 in $A(N,M)$.  By (\ref{e43}), for any $N\in\mathbb{N}$ and any $\alpha>0$, there exists an integer $M>N$, such that for all $\eta>0$ small enough,
\begin{equation}\label{e49}
\mathbf{P}_{\eta}^{\lambda}[\textrm{Circuit}(N,M)]\geq 1-\alpha.
\end{equation}

Assume now that the event $\mathcal{F}(K,\varepsilon,\delta,N,M)\cap \textrm{Circuit}(N,M)$ holds.  Let us explain that if $\ell$ and $\ell'$ are two distinct loop of $\Gamma_{\eta}^{\lambda}$ in $\Lambda_N$ with diameter at least $\delta N$, then
\begin{equation}\label{e50}
\mathbf{N}(\ell,\ell';\varepsilon,N,M)\leq K^2+3K.
\end{equation}
Suppose that $\ell$ (resp. $\ell'$) is a boundary loop of the open cluster $\mathcal{C}$ (resp. $\mathcal{C}'$); it is possible that $\mathcal{C}=\mathcal{C}'$.  Then there is a chain $(\mathcal{C}_1=\mathcal{C},\mathcal{C}_2,\ldots,\mathcal{C}_k=\mathcal{C}')$ of $k\leq K$ open clusters of $\omega_{\eta}^{\lambda}$, each of diameter at least $\varepsilon N$ and fully contained in $\Lambda_{M}$.  For $1\leq i\leq k$, let $\ell_i^+$ denote the clockwise boundary loop of $\mathcal{C}_i$.  Suppose that there are $m$ clusters in $\{\mathcal{C}_2,\ldots,\mathcal{C}_k\}$ surrounding $\mathcal{C}_1$, which we denote by $\mathcal{C}_{i_1},\ldots,\mathcal{C}_{i_{m}}$ when $m\geq 1$, with $1<i_1<\cdots<i_m\leq k$; suppose that there are $n$ clusters in $\{\mathcal{C}_1,\ldots,\mathcal{C}_{k-1}\}$ surrounding $\mathcal{C}_k$, which we denote by $\mathcal{C}_{i_{m+1}},\ldots,\mathcal{C}_{i_{m+n}}$ when $n\geq 1$, with $i_m\leq i_{m+1}<\cdots<i_{m+n}<k$.  For $1\leq j\leq m$ (resp. $m+1\leq j\leq m+n$), let $\ell_{i_j}^{-}$ denote the counterclockwise boundary loop of $\mathcal{C}_{i_j}$ that surrounds $\mathcal{C}_1$ (resp. $\mathcal{C}_k$).  See Figure \ref{finitechaining} for an illustration.  Write $i_0:=1$ and $i_{m+n+1}:=k$.  We make the following simple observations:
\begin{itemize}
  \item When $m\geq 1$, for any $0\leq j\leq m-1$, $(\ell_{i_j}^+,\ell_{i_j+1}^+,\ldots,\ell_{i_{j+1}-1}^+,\ell_{i_{j+1}}^-)$ is a chain of loops, each of diameter at least $\varepsilon N$ and contained in $\Lambda_{M}$.
  \item $(\ell_{i_m}^+,\ell_{i_m+1}^+,\ldots,\ell_{i_{m+1}-1}^+,\ell_{i_{m+1}}^+)$ is a chain of loops, each of diameter at least $\varepsilon N$ and contained in $\Lambda_{M}$.
  \item When $n\geq 1$, for any $m+1\leq j\leq m+n$, $(\ell_{i_j}^-,\ell_{i_j+1}^+,\ldots,\ell_{i_{j+1}-1}^+,\ell_{i_{j+1}}^+)$ is a chain of loops, each of diameter at least $\varepsilon N$ and contained in $\Lambda_{M}$.
\end{itemize}
Moreover, since for each $1\leq j\leq m+n$ both $\ell_{i_j}^{-}$ and $\ell_{i_j}^{+}$ have diameter at least $\delta N$, there is a chain of at most $K$ closed clusters connecting the two closed clusters having $\ell_{i_j}^{-}$ and $\ell_{i_j}^{+}$ as their boundary loops, each of diameter at least $\varepsilon N$.  Note that all the closed clusters in this chain except the starting and ending ones are between $\ell_{i_j}^{-}$ and $\ell_{i_j}^{+}$, and do not surround  $\ell_{i_j}^{-}$.  See Figure \ref{finitechaining} again.  Thus, $\ell_{i_j}^{-},\ell_{i_j}^{+}$ together with the counterclockwise boundary loops of the closed clusters between $\ell_{i_j}^{-}$ and $\ell_{i_j}^{+}$ in this chain produces a chain of loops that ensures $\mathbf{N}(\ell_{i_j}^{-},\ell_{i_j}^{+};\varepsilon,N,M)\leq K$.  If $\ell$ is a counterclockwise boundary loop of $\mathcal{C}_1$, then using the same argument as above we obtain $\mathbf{N}(\ell,\ell_1^+;\varepsilon,N,M)\leq K$.  Similarly, we have $\mathbf{N}(\ell',\ell_k^+;\varepsilon,N,M)\leq K$.  The above arguments imply that
\begin{align*}
\mathbf{N}(\ell,\ell';\varepsilon,N,M)&\leq K+\mathbf{N}(\ell,\ell_1^+;\varepsilon,N,M)+\mathbf{N}(\ell',\ell_k^+;\varepsilon,N,M)\\
&\quad+\mathbf{1}_{\{m+n\geq 1\}}\sum_{j=1}^{m+n}\mathbf{N}(\ell_{i_j}^-,\ell_{i_j}^+;\varepsilon,N,M)\\
&\leq K+(m+n+2)K\leq K^2+3K,
\end{align*}
which gives (\ref{e50}).  Combining (\ref{e48}), (\ref{e49}) and (\ref{e50}), we get that for any fixed $N\in\mathbb{N}$, $\delta\in(0,1)$ and $\alpha>0$, there exist $\varepsilon\in(0,\delta)$ and $M,K\in \mathbb{N}$ with $M>N$, such that for all small enough $\eta>0$, we have
\begin{equation*}
\mathbf{P}_{\eta}^{\lambda}\left[\begin{aligned}
&\mbox{for any two loops $\ell,\ell'$ in $\Gamma_{\eta}^{\lambda}$ of diameter at least $\delta N$ and}\\
&\mbox{contained in $\Lambda_N$, we have $\mathbf{N}(\ell,\ell';\varepsilon,N,M)\leq K^2+3K$}
\end{aligned}
\right]\geq 1-2\alpha.
\end{equation*}
Using this inequality and Theorem \ref{t1}, it is easy to obtain that for any fixed $N\in\mathbb{N}$ and $\delta\in(0,1)$,
\begin{equation*}
\mathbf{P}\left[\begin{aligned}
&\mbox{for any two loops $\ell,\ell'$ in $\Gamma_{\infty}^{\lambda}$ of diameter at least $\delta N$ and contained}\\
&\mbox{in $\Lambda_N$, there is a chain of loops connecting $\ell$ and $\ell'$ in $\Gamma_{\infty}^{\lambda}$}
\end{aligned}
\right]=1.
\end{equation*}
This and (\ref{noinfinite}) imply the finite chaining property of $\Gamma_{\infty}^{\lambda}$ immediately.

(\ref{triple}) Using Theorem \ref{t1} and the 6-arm exponent being larger than 2 (see (\ref{e67})), property (\ref{triple}) follows from the proof for the critical case, see the proof of \cite[Theorem 2, item 1]{CN06}.  Indeed, a triple point would imply for $\omega_{\eta}^{\lambda}$ the occurrence of a 6-arm event for an annulus centered at this point whose ratio of inner to outer radius goes to zero as $\eta\rightarrow 0$, leading to a contradiction.

(\ref{parity}) This property for the critical case was discussed just below \cite[Proposition 2.7]{GPS13}.  The proof is analogous to that of (\ref{triple}).  Indeed, by Theorem \ref{t1} and the 6-arm exponent being larger than 2, if two nested loops of $\Gamma_{\infty}^{\lambda}$ touch each other and have the same orientation, then a touching point would imply for $\omega_{\eta}^{\lambda}$ the occurrence of a 6-arm event for an annulus centered at this point whose ratio of inner to outer radius goes to zero as $\eta\rightarrow 0$, leading to a contradiction.  Similar argument also applies to the case that two non-nested loops of $\Gamma_{\infty}^{\lambda}$ touch each other and have opposite orientation.

(\ref{largest}) This property essentially follows from the proof of \cite[Theorem 1.1(3)]{CJM09}.  Without loss of generality, assume that $z=0$ and $\lambda<0$.  Consider the sequence of dyadic nested annuli centered at 0 with diameter smaller than the correlation length $\tilde{L}(\lambda)$.  By the near-critical RSW and FKG inequality together with Theorem \ref{t1}, we know that there exists a constant $p_0>0$, such for each annulus above, with probability at least $p_0$ there is a loop of $\Gamma_{\infty}^{\lambda}$ surrounding 0 in the annulus.  Then, the Borel-Cantelli lemma implies that there a.s. exist infinitely many nested loops of $\Gamma_{\infty}^{\lambda}$ surrounding 0 with diameters going to zero.  Moreover, such family of nested loops is a.s. countable, since otherwise there would exist $N\in\mathbb{N}$ such that with positive probability there exist uncountably many loops of $\Gamma_{\infty}^{\lambda}$ of diameter larger than $1/N$ that intersect $\Lambda_{N}$, contradicting (\ref{local}).

By (\ref{e43}) and (\ref{e44}), for any $\alpha>0$ there exists $R=R(\lambda,\alpha)$ such that for all $\eta>0$ small enough, with probability at least $1-\alpha$, there is a closed circuit of $\omega_{\eta}^{\lambda}$ surrounding 0 in $A(R,2R)$, with this circuit lying in the infinite closed cluster.  Note that the occurrence of this event implies that all the loops of $\Gamma_{\eta}^{\lambda}$ that surround 0 are contained in $\Lambda_{2R}$.  Thus, with probability at least $1-\alpha$, all the loops of $\Gamma_{\infty}^{\lambda}$ that surround 0 are contained in $\Lambda_{3R}$.  Letting $\alpha\rightarrow 0$, we get that there is a.s. an outermost loop in the loops of $\Gamma_{\infty}^{\lambda}$ which surround 0.
\end{proof}

\section{Convergence of exploration paths and loop ensembles in a disk with Dobrushin and monochromatic boundary conditions}\label{disk}
In this section we prove Theorem \ref{t3}, using the convergence result for $\Gamma_{\eta}^{\lambda}$ in the whole plane, Theorem \ref{t1}.  Note that in the critical case, Camia and Newman first proved in \cite{CN07} that the percolation exploration paths converge to an SLE$_6$ curve in a disk with Dobrushin boundary conditions, and then in \cite{CN06} they used this result to prove the convergence of the collection of cluster boundary loops in a disk with monochromatic boundary conditions, and finally got the convergence of $\Gamma_{\eta}^{\lambda=0}$ in the whole plane.  Here we somewhat reverse the steps of the above strategy.  The reason that we can do this is that we already have the near-critical quad-crossing scaling limit result proved in \cite{GPS18a}.

Before proving Theorem \ref{t3}, we need to define a distance between two finite collections $\mathcal {F},\mathcal {F}'$ of curves or loops in the plane:
\begin{equation*}
dist(\mathcal {F},\mathcal {F}'):=\min_{\phi}\max_{\gamma\in\mathcal {F}}d(\gamma,\phi(\gamma)),
\end{equation*}
where the infimum is taken over all bijections $\phi:\mathcal {F}\rightarrow\mathcal {F}'$.  If $|\mathcal {F}|\neq|\mathcal {F}'|$, then we
we let the distance be infinite.

\begin{proof}[Proof of Theorem \ref{t3}]
(\ref{sle}) Right after Theorem 7 in \cite{CN06}, Camia and Newman used ``excursions" from the full-plane CLE$_{6}$ to generate a chordal SLE$_6$ curve in a Jordan domain.  To deduce that the generated curve is indeed an SLE$_6$ curve, they used similar excursion-construction for critical percolation exploration paths, together with the result that these paths converge to an SLE$_6$ curve, which was proved in \cite{Smi01,CN07}.  To show Theorem \ref{t3}(\ref{sle}), we will apply such excursion-construction for both $\gamma_{\eta}^{\lambda}(\mathbb{D},a,b)$ and the limiting curve, combined with Theorem \ref{t1}.  Notice that we cannot apply Camia-Newman's method directly, since their argument relies on the convergence result for critical exploration paths, while the near-critical analog of this result is what we need to prove.  Let us mention that analogous excursion-construction was also used in \cite[Section 6]{BDCH16} for the critical Ising model.

We start with some definitions.  Given a loop $\ell\in\Gamma_{\infty}^{\lambda}$ intersecting $\partial\mathbb{D}$, we call a non-trivial subcurve of $\ell$ a \textbf{(continuum) excursion} if its endpoints are on $\partial\mathbb{D}$ and the remaining part is contained in $\mathbb{D}$.  An excursion of $\ell$ is directed by following the orientation of $\ell$.   Using the half-plane 3-arm exponent being 2 together with Theorem \ref{t1}, it is standard to deduce that a.s. every continuum excursion has distinct starting and ending points, and distinct continuum excursions have no endpoints in common.  Let $\mathscr{E}_{\infty}$ denote the collection of all excursions of the loops in $\Gamma_{\infty}^{\lambda}$ that intersect $\partial\mathbb{D}$.  For $\delta>0$, let $\mathscr{E}_{\infty}(\delta)$ denote the collection of all excursions of $\mathscr{E}_{\infty}$ with diameter larger than $\delta$.  Similarly to the proof of Theorem \ref{t2}(\ref{local}), by using RSW and the BK inequality it is easy to show that for a fixed $\delta>0$, there are a.s. only finitely many excursions in $\mathscr{E}_{\infty}(\delta)$.  Similarly to their continuum versions, we denote respectively by $\mathscr{E}_{\eta}$ and $\mathscr{E}_{\eta}(\delta)$ the two corresponding collections of \textbf{(discrete) excursions} for the loops of $\Gamma_{\eta}^{\lambda}$ intersecting $\partial\mathbb{D}^{\eta}$.

To study the convergence of $\gamma_{\eta}^{\lambda}:=\gamma_{\eta}^{\lambda}(\mathbb{D},a,b)$, we will use excursions in $\mathscr{E}_{\infty}(\delta)$ and $\mathscr{E}_{\eta}(\delta)$ respectively to construct $\delta$-cutoff versions of continuum and discrete exploration paths in the following.

The right boundary $\partial_{a,b}\mathbb{D}$ (resp. left boundary $\partial_{b,a}\mathbb{D}$) is the open arc of $\partial\mathbb{D}$ traversed counterclockwise from $a$ to $b$ (resp. $b$ to $a$).  Fix $\delta\in(0,1)$ and write $\mathscr{E}_{\infty}^{LR}(\delta)$ for the set of all the excursions of $\mathscr{E}_{\infty}$ starting at the left boundary and ending at the right, with all their endpoints at a distance larger than $\delta$ from both $a$ and $b$; write $\mathscr{E}_{\infty}^{RL}(\delta)$ for the analogous set of excursions from the right to the left boundary.  For a fixed $\delta>0$, there are a.s. only finitely many such excursions, which are ordered by the ordering of their intersection points with $\partial_{a,b}\mathbb{D}$.  Thus, these excursions alternate in direction.  Between an excursion $\gamma^{RL}$ in $\mathscr{E}_{\infty}^{RL}(\delta)$ and the next excursion $\gamma^{LR}$ in $\mathscr{E}_{\infty}^{LR}(\delta)$, there are a.s. finitely many excursions in $\mathscr{E}_{\infty}(\delta)$ with endpoints on $\partial_{b,a}\mathbb{D}$, and we choose the maximal ones in the sense that they are not ``contained inside" any other excursion of the same type (the arc on $\partial_{b,a}\mathbb{D}$ connecting the endpoints of such a maximal excursion is not contained in the arc on $\partial_{b,a}\mathbb{D}$ connecting the endpoints of any other excursion of the same type); note that each ``maximal excursion" and the arc on $\partial_{b,a}\mathbb{D}$ which connects the endpoints of the excursion produce a counterclockwise loop.  These maximal excursions are ordered by the ordering of their intersection points with the reversal of the arc $\partial_{b,a}\mathbb{D}$, and can be used to make a ``bridge" between $\gamma^{RL}$ and $\gamma^{LR}$ by concatenating $\gamma^{RL},\gamma^{LR}$, such ordered maximal excursions, and all the boundary arcs in the reversal of $\partial_{b,a}\mathbb{D}$ that start at the ending point of an excursion above and end at the starting point of the next one.  Between an excursion in $\mathscr{E}_{\infty}^{LR}(\delta)$ and the next excursion in $\mathscr{E}_{\infty}^{RL}(\delta)$, one can do a similar construction.

\begin{figure}
\centering
\includegraphics[height=0.45\textwidth]{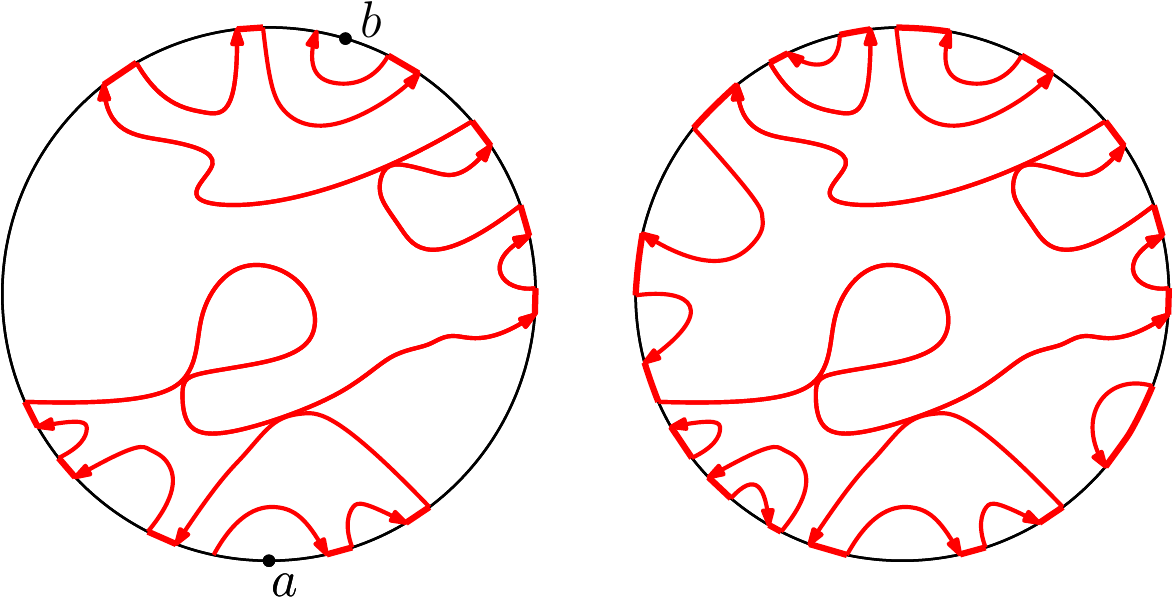}
\caption{\emph{Left}: The continuum $\delta$-exploration path $\gamma_{\infty}^{\lambda,\delta}$ is indicated by the concatenation of red curves.  The excursions and arcs on $\partial \mathbb{D}$ used in the construction of $\gamma_{\infty}^{\lambda,\delta}$ are represented by red arrows and heavy red lines, respectively. \emph{Right}: The continuum $\varepsilon$-loops in $\Gamma_{\infty,\partial \mathbb{D}}^{\lambda,\varepsilon}(\delta)$ are indicated by red loops.}\label{CLEcut}
\end{figure}

By doing this for each pair of successive excursions in $\mathscr{E}_{\infty}^{LR}(\delta)\cup\mathscr{E}_{\infty}^{RL}(\delta)$, we get a curve that connects two points on $\partial\mathbb{D}$, which is called the \textbf{continuum $\delta$-exploration path}, denoted by $\gamma_{\infty}^{\lambda,\delta}$.  See the left-hand side of Figure \ref{CLEcut}.  When $\mathscr{E}_{\infty}^{LR}(\delta)$ or $\mathscr{E}_{\infty}^{RL}(\delta)$ is empty, write $\gamma_{\infty}^{\lambda,\delta}:=\partial_{a,b}\mathbb{D}$.  The collection of excursions used in the construction of $\gamma_{\infty}^{\lambda,\delta}$ is denoted by $\mathscr{E}_{\infty,a,b}(\delta)$.  Note that $\mathscr{E}_{\infty,a,b}(\delta)\subset\mathscr{E}_{\infty}(\delta)$.  The \textbf{discrete $\delta$-exploration path} $\gamma_{\eta}^{\lambda,\delta}$ for $(\mathbb{D}^{\eta},a^{\eta},b^{\eta})$ is defined similarly, and we let $\mathscr{E}_{\eta,a,b}(\delta)$ denote the collection of discrete excursions used in the construction of $\gamma_{\eta}^{\lambda,\delta}$.

We want to show that $\mathscr{E}_{\eta,a,b}(\delta)$ converges in probability in the metric $dist$ to $\mathscr{E}_{\infty,a,b}(\delta)$, as $\eta\rightarrow 0$.  For this, we need the following facts.  First, as we have mentioned above, a.s. there are finitely many excursions in $\mathscr{E}_{\infty}(\delta)$.  Second, by classical arguments using the half-plane 3-arm exponent being 2, we know that for all small $\eta>0$, with high probability a loop of $\Gamma_{\eta}^{\lambda}$ with diameter larger than $\delta$ cannot come close to $\partial\mathbb{D}^{\eta}$ without nearby intersecting it (this fact was also used in the last paragraph of the proof of \cite[Theorem 12]{CCK19}), which, together with Theorem \ref{t1}, implies that as $\eta\rightarrow 0$, for each excursion in $\mathscr{E}_{\infty}(\delta)$ there exists an excursion in $\mathscr{E}_{\eta}(\delta)$ converging in probability to it in the metric $d$.  Third, using the 6-arm exponent being strictly larger than 2, it is standard to show that for all small $\eta>0$, with high probability there do not exist two distinct excursions of $\mathscr{E}_{\eta}(\delta)$ to be very close to each other in the metric $d$.  Then, the almost sure convergence of $\Gamma_{\eta}^{\lambda}$ to $\Gamma_{\infty}^{\lambda}$ in $D_{\mathcal{L}}$, combined with the above three facts, implies that $dist(\mathscr{E}_{\eta}(\delta),\mathscr{E}_{\infty}(\delta))\rightarrow 0$ in probability as $\eta\rightarrow 0$.  From this and the construction of $\mathscr{E}_{\infty,a,b}(\delta)$ and $\mathscr{E}_{\eta,a,b}(\delta)$, we deduce that
\begin{equation*}
dist(\mathscr{E}_{\eta,a,b}(\delta),\mathscr{E}_{\infty,a,b}(\delta))\rightarrow 0\quad\mbox{in probability as $\eta\rightarrow 0$}.
\end{equation*}
Combined with the construction of the continuum and discrete $\delta$-exploration paths, this implies that
\begin{equation}\label{e51}
d(\gamma_{\eta}^{\lambda,\delta},\gamma_{\infty}^{\lambda,\delta})\rightarrow 0\quad\mbox{in probability as $\eta\rightarrow 0$}.
\end{equation}

For a point $z\in\partial\mathbb{D}$, let $\ell_{\infty,z}^{\delta}$ (resp. $\ell_{\eta,z}^{\delta}$) be the minimal loop of $\Gamma_{\infty}^{\lambda}$ (resp. $\Gamma_{\eta}^{\lambda}$) which surrounds $z$ and has distance larger than $\delta$ from $z$.  By the definition of $\gamma_{\infty}^{\lambda,\delta}$, it is clear that for any $\delta_1,\delta_2\in(0,\delta]$, we have
\begin{equation}\label{e52}
d(\gamma_{\infty}^{\lambda,\delta_1},\gamma_{\infty}^{\lambda,\delta_2})\leq 2\delta+\diam(\ell_{\infty,a}^{\delta})+\diam(\ell_{\infty,b}^{\delta}).
\end{equation}
By Theorem \ref{t2}(\ref{largest}), a.s. each of $a$ and $b$ is surrounded by an infinite family of nested loops with diameters going to zero.
Thus, by (\ref{e52}), there exists a random curve $\gamma_{\infty}^{\lambda}$ from $a$ to $b$, which is measurable with respect to $\Gamma_{\infty}^{\lambda}$, such that
\begin{equation}\label{e54}
d(\gamma_{\infty}^{\lambda,\delta},\gamma_{\infty}^{\lambda})\rightarrow 0\quad\mbox{a.s. as $\delta\rightarrow 0$}.
\end{equation}
By the definitions of $\gamma_{\eta}^{\lambda,\delta}$ and $\gamma_{\eta}^{\lambda}$, it is clear that for all small enough $\eta>0$, we have
\begin{equation}\label{e75}
d(\gamma_{\eta}^{\lambda},\gamma_{\eta}^{\lambda,\delta})\leq 2\delta+\diam(\ell_{\eta,a}^{\delta})+\diam(\ell_{\eta,b}^{\delta}).
\end{equation}
Using RSW and FKG, it is easy to see that by letting $\delta\rightarrow 0$ and then $\eta\rightarrow 0$, both $\diam(\ell_{\eta,a}^{\delta})$ and $\diam(\ell_{\eta,b}^{\delta})$ go to zero in probability.  It follows from this and (\ref{e75}) that for any $\alpha>0$, there exists $\delta_0>0$, such that for any $\delta\in(0,\delta_0)$ and all sufficiently small $\eta>0$,
\begin{equation}\label{e53}
\mathbf{P}_{\eta}^{\lambda}[d(\gamma_{\eta}^{\lambda},\gamma_{\eta}^{\lambda,\delta})<\alpha]\geq 1-\alpha.
\end{equation}
Combining (\ref{e51}), (\ref{e54}) and (\ref{e53}), we derive the desired result: $d(\gamma_{\eta}^{\lambda},\gamma_{\infty}^{\lambda})\rightarrow 0$ in probability as $\eta\rightarrow 0$.

(\ref{cle}) By Theorem \ref{t1} and the fact that a.s. there is no non-trivial loop of $\Gamma_{\infty}^{\lambda}$ that is contained in $\overline{\mathbb{D}}$ and touches $\partial \mathbb{D}$ (which follows from the half-plane 3-arm exponent being 2), it follows that the collection of loops of $\Gamma_{\eta}^{\lambda}$ which are contained in the interior of $\mathbb{D}^{\eta}$ converges a.s. in $d_{\mathcal{L}}$ to the collection of loops of $\Gamma_{\infty}^{\lambda}$ which are contained in $\mathbb{D}$.  Thus, to prove item (\ref{cle}) we just need to consider the collection of loops of $\Gamma_{\eta,\mathbb{D}}^{\lambda}$ that touch $\partial\mathbb{D}^{\eta}$, called the collection of \textbf{boundary touching loops} of $\Gamma_{\eta,\mathbb{D}}^{\lambda}$ and denoted by $\Gamma_{\eta,\partial \mathbb{D}}^{\lambda}$.  We will prove the convergence of $\Gamma_{\eta,\partial\mathbb{D}}^{\lambda}$ by an argument similar to that for (\ref{sle}).

Recall that a.s. every continuum excursion has distinct starting and ending points, and distinct continuum excursions have no endpoints in common.  For an excursion $\gamma$, let $\underline{\gamma}$ be the segment of $\partial\mathbb{D}$ traversed counterclockwise from the ending point of $\gamma$ to its starting point.  Assume that $\gamma\in\mathscr{E}_{\infty}(\delta)$ and $0<\varepsilon<\delta/2<1/2$.  Consider the collection of all maximal excursions of $\mathscr{E}_{\infty}(\varepsilon)$ with endpoints in the interior of $\underline{\gamma}$, where ``maximal" means that the segment of $\underline{\gamma}$ which connects the endpoints of such an excursion does not lie in the segment of $\underline{\gamma}$ which connects the endpoints of any other excursion of the same type.  Observe that each maximal excursion, together with the segment of the reversal of $\underline{\gamma}$ which connects the endpoints of the excursion, produces a clockwise loop.  We order these maximal excursions by the ordering of their intersection points with $\underline{\gamma}$ to obtain a sequence $\gamma_1,\ldots,\gamma_k$.  Write $\gamma_0=\gamma_{k+1}=\gamma$.  For each $0\leq i\leq k$, let $\underline{\gamma}_i$ be the segment of $\underline{\gamma}$ that starts at the ending point of $\gamma_i$ and ends at the starting point of $\gamma_{i+1}$.  By concatenating $\gamma_0,\underline{\gamma}_0,\gamma_1,\underline{\gamma}_1,\dots,\gamma_k,\underline{\gamma}_k$, we get the \textbf{continuum $\varepsilon$-loop} $\ell_{\infty}^{\varepsilon}(\gamma)$ containing $\gamma$.  The collection of such loops constructed for all excursions in $\mathscr{E}_{\infty}(\delta)$ is denoted by $\Gamma_{\infty,\partial \mathbb{D}}^{\lambda,\varepsilon}(\delta)$.  See the right-hand side of Figure \ref{CLEcut}.  Similarly, for each $\gamma_{\eta}\in\mathscr{E}_{\eta}(\delta)$, one can construct the \textbf{discrete $\varepsilon$-loop} $\ell_{\eta}^{\varepsilon}(\gamma_{\eta})$ containing $\gamma_{\eta}$, and denote by $\Gamma_{\eta,\partial \mathbb{D}}^{\lambda,\varepsilon}(\delta)$ the collection of all such loops.

It is obvious that for any $0<\varepsilon<\delta/2<1/2$, $\varepsilon_1,\varepsilon_2\in(0,\varepsilon]$ and any $\gamma\in\mathscr{E}_{\infty}(\delta)$, we have
\begin{equation*}
d(\ell_{\infty}^{\varepsilon_1}(\gamma),\ell_{\infty}^{\varepsilon_2}(\gamma))\leq 2\varepsilon.
\end{equation*}
Then, since there are a.s. finitely many excursions in $\mathscr{E}_{\infty}(\delta)$, we get that for any fixed $\delta\in(0,1)$, a.s. for each $\gamma\in\mathscr{E}_{\infty}(\delta)$ there exists a random loop $\ell_{\infty}(\gamma)$ containing $\gamma$, which is measurable with respect to $\Gamma_{\infty}^{\lambda}$, such that
\begin{equation}\label{e55}
d(\ell_{\infty}^{\varepsilon}(\gamma),\ell_{\infty}(\gamma))\rightarrow 0\quad\mbox{uniformly for all $\gamma\in\mathscr{E}_{\infty}(\delta)$ as $\varepsilon\rightarrow 0$}.
\end{equation}
Denote by $\Gamma_{\infty,\partial D}^{\lambda}(\delta)$ the collection of all such limiting loops which contain excursions in $\mathscr{E}_{\infty}(\delta)$.
Note that it is possible that $\ell_{\infty}(\gamma)=\ell_{\infty}(\gamma')$ for distinct $\gamma,\gamma'$ in $\mathscr{E}_{\infty}(\delta)$.

As discussed in the proof of (\ref{sle}), for any fixed $\varepsilon\in(0,1/2)$ we have $dist(\mathscr{E}_{\eta}(\varepsilon),\mathscr{E}_{\infty}(\varepsilon))\rightarrow 0$ in probability as $\eta\rightarrow 0$.  From this and the construction of continuum and discrete $\varepsilon$-loops, we obtain
\begin{equation}\label{e57}
dist(\Gamma_{\infty,\partial\mathbb{D}}^{\lambda,\varepsilon}(\delta),\Gamma_{\eta,\partial\mathbb{D}}^{\lambda,\varepsilon}(\delta))\rightarrow 0\quad\mbox{in probability as $\eta\rightarrow 0$}.
\end{equation}

For $\gamma_{\eta}\in\mathscr{E}_{\eta}(\delta)$, let $\ell_{\eta}(\gamma_{\eta})$ be the loop in $\Gamma_{\eta,\partial\mathbb{D}}^{\lambda}$ that contains $\gamma_{\eta}$.  Denote by $\Gamma_{\eta,\partial\mathbb{D}}^{\lambda}(\delta)$ the collection of all the loops in $\Gamma_{\eta,\partial\mathbb{D}}^{\lambda}$ which contain excursions in $\mathscr{E}_{\eta}(\delta)$.  By the constructions, we know that for any $\gamma_{\eta}\in\mathscr{E}_{\eta}(\delta)$,
\begin{equation}\label{e56}
d(\ell_{\eta}^{\varepsilon}(\gamma_{\eta}),\ell_{\eta}(\gamma_{\eta}))\leq 2\varepsilon\quad\mbox{for all small enough $\eta>0$.}
\end{equation}
 Combining (\ref{e55}), (\ref{e57}) and (\ref{e56}), we obtain that for any $\delta\in(0,1)$,
\begin{equation}\label{e58}
dist(\Gamma_{\eta,\partial\mathbb{D}}^{\lambda}(\delta),\Gamma_{\infty,\partial\mathbb{D}}^{\lambda}(\delta))\rightarrow 0\quad\mbox{in probability as $\eta\rightarrow 0$}.
\end{equation}

By RSW and FKG, it is standard to show that with high probability each large boundary touching loop contains a large excursion, and each boundary point has a boundary touching loop close to it: For any fixed $\alpha,\delta'\in(0,1)$, there exists $\delta\in(0,\delta')$ such that for all $\eta>0$ small enough,
\begin{equation}\label{e60}
\mathbf{P}_{\eta}^{\lambda}\left[\begin{aligned}
&\mbox{all loops in $\Gamma_{\eta,\partial\mathbb{D}}^{\lambda}$ of diameter larger than $\delta'$ are contained in $\Gamma_{\eta,\partial\mathbb{D}}^{\lambda}(\delta)$;}\\
&\mbox{for each $z\in\partial\mathbb{D}$ there is a loop in $\Gamma_{\eta,\partial\mathbb{D}}^{\lambda}$ of diameter at most $\delta$ which}\\
&\mbox{has a distance at most $\delta$ from $z$}
\end{aligned}
\right]\geq 1-\alpha.
\end{equation}
Write $\Gamma_{\infty,\partial\mathbb{D}}^{\lambda}:=\{z\in\partial\mathbb{D}:\mbox{trivial loop at $z$}\}\bigcup_{\delta\in(0,1)}\Gamma_{\infty,\partial\mathbb{D}}^{\lambda}(\delta)$, which is measurable with respect to $\Gamma_{\infty}^{\lambda}$.
Combining (\ref{e58}) and (\ref{e60}), we obtain the convergence of discrete boundary-touching loop ensembles:
\begin{equation*}
d_{\mathcal{L}}(\Gamma_{\eta,\partial\mathbb{D}}^{\lambda},\Gamma_{\infty,\partial\mathbb{D}}^{\lambda})\rightarrow 0\quad\mbox{in probability as $\eta\rightarrow 0$}.
\end{equation*}
This, together with the argument at the beginning of the proof of (\ref{cle}), gives the desired result.
\end{proof}

\section*{Acknowledgements}
This work is supported by the National Natural Science Foundation of China (No. 12288201) and the Key Laboratory of Random Complex Structures and Data Science, CAS (No. 2008DP173182).

\subsection*{AI Use Statement}
AI was not used by the author at any stage of this work.

\end{document}